\documentclass[leqno,11pt]{article}
\usepackage{amsfonts}
\usepackage{amsmath}
\usepackage{amsthm}
\usepackage{hyperref}
\usepackage{geometry}
\usepackage{fancybox}
\usepackage{boxedminipage}
\usepackage{color}

\newtheorem{theorem}{Theorem}

\newtheorem{corollary}[theorem]{Corollary}

\newtheorem{lemma}[theorem]{Lemma}

\newtheorem{proposition}[theorem]{Proposition}

\theoremstyle{definition}

\newtheorem{remark}{Remark}

\newtheorem{assumption}{Assumption}

\def\R{\ensuremath{\mathbb R}}
\def\Q{\ensuremath{\mathbb Q}}

\def\C{\ensuremath{\mathbb C}}
\def\N{\ensuremath{\mathbb N}}
\def\R{\ensuremath{\mathbb R}}
\def\Q{\ensuremath{\mathbb Q}}

\def\C{\ensuremath{\mathbb C}}
\def\N{\ensuremath{\mathbb N}}

\def\Ds{\displaystyle}

\def\intdouble{\int \! \! \int}

\def\eps{\varepsilon}
\def\calA{\mathcal{A}}
\def\calC{\mathcal{C}}
\def\calN{\mathcal{N}}
\def\calJ{\mathcal{J}}
\def\calI{\mathcal{I}}
\def\calG{\mathcal{G}}
\def\calH{\mathcal{H}}
\def\calM{\mathcal{M}}
\def\calR{\mathcal{R}}
\def\calL{\mathcal{L}}
\def\calP{\mathcal{P}}
\def\calS{\mathcal{S}}
\def\calK{\mathcal{K}}
\def\calB{\mathcal{B}}
\def\calU{\mathcal{U}}
\def\calF{\mathcal{F}}

\def\frakC{\mathfrak{C}}

\newcommand{\psU}[2]{\left\langle #1, #2 \right\rangle_{\mathcal{U}_2}}
\newcommand{\normeU}[1]{\left\| #1 \right\|_{\mathcal{U}_2}}

\title{New results on biorthogonal families in cylindrical domains and controllability consequences}
\author{
\textsc{F.~Ammar Khodja}\footnote{Laboratoire de Mathématiques UMR 6623, Université de Franche-Comté, 16, route de Gray, 25030, Besançon cedex, France}, 
\textsc{A.~Benabdallah}\footnote{Aix Marseille Univ, CNRS, I2M, Marseille, France}, 
\textsc{M.~Gonz\'alez-Burgos}\footnote{Dpto. Ecuaciones Diferenciales y Análisis Numérico and, Instituto de Matemáticas de la Universidad de Sevilla, Universidad de Sevilla, Spain}, 
\textsc{M.~Morancey}\footnote{Aix Marseille Univ, CNRS, I2M, Marseille, France} 
\\
and \textsc{L.~de Teresa}\footnote{Instituto de Matemáticas, Universidad Nacional Autónoma de México, Circuito Exterior, C.U., C. P. 04510, Mexico City, Mexico} }
\date{}

\begin{document}

\maketitle

\begin{abstract}
In this article we consider moment problems equivalent to null controllability of some linear parabolic partial differential equations in space dimension higher than one. 
For these moment problems, we prove existence of an associated biorthogonal family and estimate its norm. 
The considered setting requires the space domain to be a cylinder and the evolution operator to be tensorized. 

Roughly speaking, we assume that the so-called Lebeau-Robbiano spectral inequality holds but only for the eigenvectors of the transverse operator. 
In the one dimensional tangent variable we assume the solvability of block moment problem as introduced in [Benabdallah, Boyer and Morancey - \textit{Ann. H. Lebesgue.} 3 (2020)]. 

We apply this abstract construction of biorthogonal families to the characterization of the minimal time for simultaneous null controllability of two heat-like equations in a cylindrical domain. To the best of our knowledge, this result is unattainable with other known techniques. 
\end{abstract}


\paragraph{Acknowledgment.}
This work owes a lot to the city of Chiclana de la Frontera where the authors spent an entire week without any disturbance. 

The second and fourth authors were partially supported by the ANR project TRECOS ANR-20-CE40-0009. 

The third author is supported by grant PID2020-114976GB-I00 funded by MICIU/ AEI/10.13039/501100011033 (Spain).


\section{Introduction}\label{s1}

\subsection{Biorthogonal families and moment method}

The first results on the boundary or internal null-controllability at a 
positive time $T$ of the heat equation were obtained in the 70's (see \cite{FR1:71}, \cite{FR2:74}, \cite{Fatto}) using the moment method.
This technique consists in writting a null controllability objective as a moment problem satisfied by the control.  
In those references, this moment problem is solved thanks to a biorthogonal family in $L^{2}(0,T \: ; \: \R)$ to $\left\{ t \mapsto e^{-\lambda
_{k}t} \right\}_{k \geq 1}$ where $\{-\lambda _{k} \}_{k\geq 1}$ is the sequence of eigenvalues of the Dirichlet Laplace operator. As proved in~\cite{Schwartz} a necessary and sufficient condition of existence of such biorthogonal family is the convergence of the series $\underset{k\geq 1}{\sum }\frac{1}{\lambda _{k}}$. 
Therefore, the Weyl's asymptotic restricts this approach to the one dimensional heat equation. 

In several space dimensions, other techniques were required to control to zero the heat equation. In particular, the use of Carleman inequalities (see \cite{LR:95} and \cite{FI:96} for the main references) generated a lot of results. 

Yet, in the last fifteen years the moment method was used again in the context of parabolic control problems. In particular, it allowed to solve control problems that seemed unattainable by Carleman's inequalities. One can cite for example, the boundary control of coupled parabolic equations \cite{FCGBdT}. 
It also allowed to deal with some parabolic control problems in which a positive minimal time or geometric conditions on the control region may be required for null controllability to hold (see for instance \cite{ABGD:14, AKBGBdT_JMAA16, Duprez_Tmin2017, Ou:19}). 
These are high-frequency phenomena due for instance to eigenvalue condensation and/or eigenvector localization which are well captured by the biorthogonal families. 

Let us also mention that this strategy also allowed to study null controllability for degenerate parabolic operators~\cite{CMV:20, CMV:21}. 

Recently, to take into account condensation of eigenvectors, the use of biorthogonal families to solve moment problems was replaced by the resolution of appropriate block moment problems (see~\cite{BBM, BM:23}) still under the assumption that the series of the inverse of the eigenvalues converges. This allows to consider situations in which the eigenvectors also condensate. 

The natural question is then to study the phenomenon in higher dimensions, \textit{i.e.}, when this series does not converge. 
Very few results are available in this direction. Let us mention the recent work~\cite{BO:24} where the authors prove null controllability in any time of two coupled heat equations in a rectangle with distinct diffusion speeds when the controls acts on two non parallel sides. Their proof relies on the moment method with a subtle decomposition of the moment problem into an infinite family of one dimensional moment problems. 
Their construction does not seem to be easily generalizable and the proved result is not contained in our study. 

Thus, the main question addressed in this paper is the construction of biorthogonal families associated with moment problems coming from parabolic control problems in space dimensions larger than one.

\subsection{Biorthogonal families in higher dimension and strategy of proof}

To expose more precisely the problem solved in this article, we introduce an abstract control problem 
\begin{equation} \label{PbAbstrait}
\left\{
\begin{aligned}
& y' + \calA y = \calB u
\\
&y(0) = y_0
\end{aligned}
\right.
\end{equation}
on a Hilbert space $H$. Assume that $\calA$ generates a $C^0-$semigroup and that $\calB$ ensures well-posedness for any $u \in L^2(0,T ; U)$ where $U$ is the Hilbert space of controls. 
The precise setting under study (with tensorized operators and cylindrical geometry) will be specified later on. 
We assume that the operator $\calA^*$ has a family of eigenvalues $\Lambda \subset (0,+\infty)$ and that the family of associated eigenvectors $\{ \phi_\lambda \}_{ \lambda \in \Lambda}$ forms a complete family in the state space. Then, the control $u \in L^2(0,T ;U)$ is such that $y(T) = 0$ if and only if 
\begin{equation} \label{PbMomentsAbstrait}
\int_0^T \left\langle u(T-t) , e^{-\lambda t} \calB^* \phi_\lambda \right\rangle_U dt = - e^{- \lambda T} \left\langle y_0 , \phi_\lambda \right\rangle, \qquad \forall \lambda \in \Lambda. 
\end{equation}
Thus, an appropriate generalization of biorthogonal families to the time exponentials is a family $\{ q_\mu \}_{\mu \in \Lambda} \subset L^2(0,T ; U)$ such that
\begin{equation} \label{BiorthoAbstrait}
\int_0^T \left\langle q_\mu(t) , e^{-\lambda t} \calB^* \phi_\lambda \right\rangle_U dt = \delta_{\lambda \mu}, \qquad
\forall \lambda, \mu \in \Lambda
\end{equation}
where $\delta_{\lambda,\mu}$ denotes the Kronecker delta function. 
As noticed for instance in~\cite{Morancey24}, when the family of eigenvectors $\{ \phi_\lambda \}_{\lambda \in \Lambda}$ forms a Hilbert basis of the state space, then~\eqref{PbMomentsAbstrait} gives that spectral null controllability in time $T$ (that is when the initial condition is any eigenvector) implies the existence of a biorthogonal family satisfying~\eqref{BiorthoAbstrait} and every bound on the control cost translates into bounds on this biorthogonal family. 
Thus, for example, if $\Omega$ is any smooth domain in $\R^n$, $\calA$ is the Laplace-Dirichlet operator and $\calB = \mathbf{1}_\omega$ with $\omega \subset \Omega$ an open set, then for any $T>0$ a biorthogonal family in the sense of~\eqref{BiorthoAbstrait} exists and satisfies 
\[
\| q_\lambda \|_{L^2((0,T) \times \omega ; \R)} \leq C e^{C \sqrt{\lambda}}, \qquad \forall \lambda \in \Lambda. 
\]
More generally, following the Lebeau-Robbiano iteration scheme (see~\cite{LR:95}) this holds in any setting where $\{ \phi_\lambda \}_{\lambda \in \Lambda}$ is a Hilbert basis and the following so-called Lebeau-Robbiano spectral inequality holds
\[
\left\| \sum_{\sqrt{\lambda} \leq N} a_\lambda \phi_\lambda \right\|_{L^2(\Omega ; \R)} 
\leq e^{C N} \left\| \sum_{\sqrt{\lambda} \leq N} a_\lambda \phi_\lambda \right\|_{L^2(\omega ; \R)} 
\]
for any $N \geq 1$ and $\{ a_\lambda \}_{\lambda \in \Lambda} \subset \R$.

\medskip
Our goal is thus to prove the existence of biorthogonal families as defined in~\eqref{BiorthoAbstrait} with suitable estimates but under weaker assumptions. Namely, we consider $\Omega = \Omega_1 \times (0,\pi)$ and the underlying evolution operator is assumed to be tensorized. A precise formulation of the assumptions is given in Section~\ref{ss21}. 
We still assume such a spectral inequality but only for the eigenvectors associated to the transverse operator. In the tangential variable, we will use results from~\cite{BBM}. Thus, our assumptions on the eigenvalues of the adjoint of the tangential operator include the summability of the series of their inverse as well as a weak-gap condition (see~\eqref{H3}).

Let us insist on the fact that our main result, Theorem~\ref{main}, is about biorthogonal families. Though the moment problem~\eqref{PbMomentsAbstrait} (and thus the definition of biorthogonal families in~\eqref{BiorthoAbstrait}) comes from the null controllability of system~\eqref{PbAbstrait}, the study of these biorthogonal families is of interest regardless of the controllability properties. 
For instance, at a given final time $T>0$, problem~\eqref{PbAbstrait} might not be null controllable whereas at the same time a biorthogonal family in the sense of~\eqref{BiorthoAbstrait} does exist.

\medskip
To avoid drowning the ideas into technicalities and notation, let us present our strategy of proof on the following example. 
Let $\Omega = (0,\pi)^2$ and
\begin{equation} \label{Ex:Dolecki2D}
\left\{
\begin{aligned}
& \partial_t y - \Delta	y = \delta_{x_0} \mathbf{1}_{\omega}(x') u(t,x',x), \quad (t,x',x) \in (0,T) \times \Omega,
\\
& y = 0, \quad \text{ on }   (0,T) \times \partial \Omega,
\\
& y(0,x',x) = y_0(x',x),  \quad (x',x) \in (0,\pi)^2. 
\end{aligned}
\right. 
\end{equation}
We emphasize that our study is not limited to this particular example but encompasses the abstract setting described in Section~\ref{s2}. 

In~\eqref{Ex:Dolecki2D}, the control has its support located on a segment parallel to one of the axes. This generalizes the study of~\cite{Dolecki:73} for the one dimensional system
\begin{equation} \label{Ex:Dolecki1D}
\left\{
\begin{aligned}
& \partial_t y - \partial_{xx}	y = \delta_{x_0} u(t,x), \quad (t,x) \in (0,T) \times (0, \pi),
\\
& y(t,0) = y(t,1) = 0, \quad t \in (0,T),
\\
& y(0,x) = y_0(x),  \quad x \in (0,\pi). 
\end{aligned}
\right. 
\end{equation}
There it is proved that the minimal time for null controllability in $H^{-1}(0, \pi ; \R)$ is
\[
T_0(x_0) = \limsup_{k \to + \infty} \frac{- \ln | \sin(k x_0) |}{k^2}.
\]
In~\cite{EHS:15}, the author proved that the $2D$ system~\eqref{Ex:Dolecki2D} is null controllable in any time $T>0$ under assumptions on $x_0$ that implies that $T_0(x_0) = 0$ and that the cost of null controllability in small time of~\eqref{Ex:Dolecki1D} is dominated by $e^{C/T}$. 

His strategy consists in proving first the null controllability of~\eqref{Ex:Dolecki2D} when $\omega = (0,\pi)$ using the null controllability of the associated one dimensional system~\eqref{Ex:Dolecki1D} and the fact that $\{ x' \mapsto \sin(mx') \}_{m \geq 1}$ is a Hilbert basis of $L^2(0, \pi ; \R)$. Then, using a Lebeau-Robbiano like strategy inspired by~\cite{BBGBO:2014}, this null controllability is transferred to~\eqref{Ex:Dolecki2D} with $\omega$ an open set in $(0,\pi)$. This step uses, in a crucial way, that the associated one dimensional problem~\eqref{Ex:Dolecki1D} is null controllable in any final time $T>0$, as well as the estimate on the cost of null controllability. 

The general construction of a biorthogonal family given by Theorem~\ref{main} applies to the moment problem associated with~\eqref{Ex:Dolecki2D} for any $x_0 \in (0, \pi)$ such that $\frac{x_0}{\pi} \not\in \Q$, which is a necessary and sufficient condition for the approximate controllability. 
The estimates on this biorthogonal family given in Theorem~\ref{main} imply that $T_0(x_0)$ is also the minimal null control time for system~\eqref{Ex:Dolecki2D} from $H^{-1}((0,\pi) \times (0,\pi) ; \R)$. 
To the best of our knowledge, such a result is not known and, at present, not attainable by other techniques than the moment method. 

Let us present our construction of a biorthogonal family associated to the problem~\eqref{Ex:Dolecki2D}. 

\begin{itemize}
\item[$\bullet$] \emph{Notion of biorthogonal family}. 

In this case the eigenvalues of (the adjoint of) the evolution operator are explicitly given by
\[
\Lambda = \{ k^2 + m^2  \: : \: k,m \geq 1 \}
\]
and for $k,m \geq 1$ an eigenvector associated to $k^2+m^2$ is given by 
\[
\varphi_{m,k} : (x',x) \in (0,\pi)^2 \mapsto  \sin(k x) \sin( mx'). 
\]
Thus, the moment problem~\eqref{PbMomentsAbstrait} reads as follows: the solution $y$ of~\eqref{Ex:Dolecki2D} satisfies $y(T)=0$ if and only if for any $k,m \geq 1$,
\begin{equation} \label{PbMoments_Dolecki2D}
\int_0^T \int_\omega u(T-t, x') e^{-(k^2 + m^2) t} \sin(m x') \sin(k x_0) dx' dt = -e^{-(k^2+m^2)T} \left\langle y_0, \varphi_{m,k} \right\rangle_{H^{-1},H^1_0}. 
\end{equation}
Thus, we look for a biorthogonal family $\{ Q_{m,k} \}_{k,m \geq 1} \subset L^2((0,T) \times \omega ; \R)$ in the sense that
\begin{equation} \label{Biortho_Dolecki2D}
\sin(k x_0) \int_0^T \int_\omega Q_{n,\ell}(t,x') e^{-(k^2 + m^2) t} \sin(m x') dx' dt = \delta_{k \ell} \delta_{m n},
\end{equation}
for any $k,\ell \geq 1$ and any $m, n \geq 1$.  

In the general setting we will look for a biorthogonal family to $F_{m,k}^{(j)}$ as defined by~\eqref{Fmk}.
\smallskip
\item[$\bullet$] \emph{A simpler problem}. 

In this article we were strongly inspired by~\cite{FR2:74}. There, the authors design biorthogonal families to $\left\{ t \mapsto e^{-\lambda t} \right\}_{\lambda \in \Lambda}$ in $L^2(0,T ; \R)$. Following~\cite{Schwartz}, their strategy first consists in solving the simpler problem to find a biorthogonal family in $L^2(0, +\infty ; \R)$ and then to deduce a biorthogonal family in $L^2(0,T ; \R)$ studying the properties of the restriction operator on appropriate spaces. We follow this idea but with a restriction in the $x'$ variable instead of the time variable. 

Thus, a first step is to design a biorthogonal family in the sense of~\eqref{Biortho_Dolecki2D} but in the simpler case where $\omega = (0,\pi)$. From previous results (for instance~\cite{BBM, GBO:19}), for any fixed $m\geq 1$, there exists $\left\{ \widetilde{q}_{m,k} \right\}_{k \geq 1} \subset L^2(0,T ; \R)$ such that 
\[
\sin(k x_0) \int_0^T \widetilde{q}_{m,\ell}(t) e^{-(k^2 + m^2) t} dt = \delta_{k \ell}, \quad \forall k, \ell \geq 1,
\]
and
\[
\| \widetilde{q}_{m,k} \|_{L^2(0,T ; \R)} \leq C \frac{e^{C \sqrt{k^2 + m^2}}}{|\sin(k x_0)|}, \quad \forall k, \ell \geq 1.
\]
This step crucially uses that the series of the inverse of the eigenvalues of (the adjoint of) the tangential operator converges. The general version of this result is Proposition~\ref{Prop:Annexe_biortho}. 

Then, we define 
\[
q_{m,k} : (t,x') \in (0,T) \times (0,\pi) \mapsto \widetilde{q}_{m,k}(t) \sin(m x'), \quad \forall k, m \geq 1.
\]
Thus, for any $k, \ell \geq 1$ and any $m, n \geq 1$, we have
\begin{align*}
&\sin(k x_0) \int_0^T \int_\omega q_{n,\ell}(t,x') e^{-(k^2 + m^2) t} \sin(m x') dx' dt 
\\
=&\sin(k x_0) \int_0^T \widetilde{q}_{n,\ell}(t) e^{-(k^2 + m^2) t} dt \int_0^\pi \sin(n x') \sin(m x') d x' 
\\
=&\delta_{mn} \sin(k x_0) \int_0^T \widetilde{q}_{n,\ell}(t) e^{-(k^2 + m^2) t} dt 
\\
=&\delta_{mn} \delta_{k\ell}
\end{align*}
and
\[
\| q_{m,k} \|_{L^2((0,T) \times (0,\pi) ; \R)} \leq C \frac{e^{C \sqrt{k^2 + m^2}}}{|\sin(k x_0)|}, \quad \forall k, \ell \geq 1.
\]
This step crucially uses orthogonality of the eigenvectors of the transverse operator which allows to consider a biorthogonal family to 
$\left\{ t \mapsto e^{-(k^2+m^2) t} \right\}_{k \geq 1}$ for every fixed $m \geq 1$.

The general version of this construction of a biorthogonal family in $L^2((0,T) \times \Omega_1 ; \calU_2)$ is given in Proposition~\ref{p2}.

\item[$\bullet$] \emph{The restriction operator}. 

Now, following the strategy developed in~\cite{FR2:74, Schwartz}, we prove that the restriction operator 
\[
\calR_\omega  : \varphi \mapsto \varphi_{| \omega}
\]
is an isomorphism between appropriate spaces. Having in mind integrated observability inequalities (see for instance~\cite[Section 3.3]{M:10}, which was the other great source of inspiration for the present paper), we introduce a weight function and prove that for $\alpha > 0$ sufficiently large we have
\begin{equation} \label{Inegalite_Pn_Dolecki2D}
\int_0^T \int_0^\pi e^{ - \frac{\alpha}{t} } \left| P_N(t,x') \right|^2 dx' dt 
\leq C \int_0^T \int_\omega e^{ - \frac{\alpha}{t} } \left| P_N(t,x') \right|^2 dx' dt,
\end{equation}
for any $N \geq 1$ and any $P_N$ given by 
\[
P_N(t,x') = \sum_{k,m \leq N} a_{m,k} e^{-(k^2 + m^2) t} \sin(m x').
\]
The weight function in the left-hand side of~\eqref{Inegalite_Pn_Dolecki2D} led us to modify the biorthogonal family designed in the previous step requiring that it vanishes near $t=0$. 

The proof of~\eqref{Inegalite_Pn_Dolecki2D} is too technical to be completely detailed in this introductory section but let us present the main ingredients. It relies on the fact that the eigenvectors in the transverse variable $\left\{ x' \mapsto \sin(m x') \right\}_{m \geq 1}$ satisfy the spectral inequality 
\[
\int_0^\pi \left| \sum_{m \leq N} b_m \sin(m x') \right|^2 dx' 
\leq C e^{ CN } \int_\omega \left| \sum_{m \leq N} b_m \sin(m x') \right|^2 dx' 
\]
and the identity
\[
\int_0^T \int_0^\pi P_N(t,x') q_{m,k}(t,x') dx' dt = a_{m,k}
\]
where $\left\{ q_{m,k} \right\}_{ k,m \geq 1}$ is the biorthogonal family designed at the previous step. The estimate of the norm of this biorthogonal family allows to estimate the coefficients $a_{m,k}$ with the norm of $P_N$ (see Lemma~\ref{l2} in the general setting). 
Then, the proof of~\eqref{Inegalite_Pn_Dolecki2D} amounts to estimate the rest of a converging series (see~\eqref{SigmaN}) which converges since the dissipation speed in the transverse variable is stronger than the cost coming from the spectral inequality. 
Hence, the proof of~\eqref{Inegalite_Pn_Dolecki2D} uses the same ingredients as the classical Lebeau-Robbiano strategy, especially from the point of view of observability as developed in~\cite{M:10}, but without using a partition of the time interval that usually requires controllability (or observability) in arbitrary small time. 

Then, inequality~\eqref{Inegalite_Pn_Dolecki2D} implies that the restriction operator
\begin{equation*}
\calR_\omega : \varphi \mapsto \varphi_{|\omega}
\end{equation*}
is an isomorphism between appropriate Hilbert spaces (see~\eqref{eta_alpha}). 
This gives, from $\left\{ q_{m,k} \right\}_{k,m \geq 1}$ (the biorthogonal family in $L^2((0,T) \times (0, \pi) ; \R)$) the sought biorthogonal family~\eqref{Biortho_Dolecki2D} satisfying
\[
\| Q_{m,k} \|_{L^2(0,T ; \R)} \leq C \frac{e^{C \sqrt{k^2 + m^2}}}{|\sin(k x_0)|}, \quad \forall k, \ell \geq 1.
\]

The general version of~\eqref{Inegalite_Pn_Dolecki2D} is given in Theorem~\ref{tPN} and the general version of the isomorphism property is given in Theorem~\ref{tmainrest}. 
\end{itemize}

\subsection{Structure of the article}

To end this introduction, let us present the structure of this article. 

In Section~\ref{s2}, we precisely state our assumptions and our main result (see Theorem~\ref{main}) concerning the existence and estimate of biorthogonal families. 

Section~\ref{s3} is devoted to the restriction operator in the variable $x'$. We state (see Theorem~\ref{tmainrest}) and prove the needed isomorphism property between appropriate spaces. 

Then, in Section~\ref{s4} we prove Theorem~\ref{main}: we design biorthogonal families in the simpler case $\omega= \Omega_1$ in Section~\ref{ss41} and detail how the isomorphism property of the restriction operator allows to conclude (see Proposition~\ref{Prop:restriction_biortho}). 

We provide in Section~\ref{s5} an application of this abstract construction of biorthogonal families to the characterization of the minimal time for simultaneous controllability of two linear parabolic partial differential equations. 

In Section~\ref{s6}, we provide an extension to the resolution of moment problem associated with operators with geometrically multiple eigenvalues. 

In Appendix~\ref{Annexe1} we recall the construction of biorthogonal families obtained in~\cite{BM:23}. Finally, in Appendix~\ref{Annexe2}, we revisit the classical Leabeau-Robbiano construction from the point of view of biorthogonal families. In particular, we prove that the obtained estimates on the restriction operator are sufficiently sharp to recover the bounds given by Miller in~\cite{M:10} on the cost of null-controllability of the heat equation in small time.


\section{Main results}\label{s2}
Let us fix $d\geq 2$, $T > 0$, $\Omega = \Omega_1 \times (0,\pi)\subset \mathbb{R}^{d}$, with $\Omega_1 \subset \R^{d-1}$ a bounded domain with boundary $\partial \Omega_1 \in C^1$, and $\omega \subset \Omega_1 $, an arbitrary non-empty open set of $\mathbb{R}^{d-1}$. 


Let us fix some general notations that will be used all along this work. First, we will write
	\begin{equation}\label{notation}
	\left\{
	\begin{array}{c}
\Ds ( x' ,x ) \in \R^d, \quad \hbox{with } x' = (x_1,\cdots,x_{d-1}) \in \Omega_1 \hbox{ and } x \in  ( 0,\pi  ), \\
	\noalign{\smallskip}
\Ds  Q_{T}:= (0,T) \times \Omega \quad \hbox{and} \quad \Sigma _{T} :=(0,T) \times \partial \Omega.
	\end{array}
	\right.
	\end{equation}

Secondly, if $ \calS \subset (0, \infty)$ is a sequence, we will use the notation $\calN_\calS $ for the counting function associated to $\calS$, i.e., for the function $\calN_\calS $ given by
	\begin{equation}\label{counting}
\calN_\calS (r) := \sharp \left\{ \lambda \in \calS : \lambda \leq r \right \}  , \quad r \in (0, \infty).
	\end{equation}

The main result of this paper establishes the existence of a biorthogonal family to an appropriate sequence of functions in $L^2 (Q_T)$. Before stating it, let us introduce the main hypotheses of this work.


\subsection{Assumptions} \label{ss21}
Let us consider two real non-decreasing sequences $\Lambda_1 \subset (0, \infty)$ and  $\Lambda_2 \subset (0, \infty) $ satisfying the following properties:

\begin{description}
\item[$\Lambda_1$:] There exist positive constants  $\kappa_1$ and $\theta_1$ such that
	\begin{equation}\label{H1}
\calN_{\Lambda_1}  ( r ) \leq \kappa_1 r^{\theta_1}, \quad \forall r \in (0, \infty)
	\end{equation}
where $ \calN_{\Lambda_1}$ is the counting function associated to $\Lambda_1$, see~\eqref{counting}. We denote the elements of $\Lambda_1$ by $\Lambda_1 := \left\{ \mu_m \right\}_{m \ge 1}$. 

\item[$\Lambda_2$:] There exist two constants $\kappa >0$ and $\theta \in (0,1)$ such that 
	\begin{equation}\label{H5}
\left| \calN_{\Lambda_2} (r_1) - \calN_{\Lambda_2} (r_2) \right| \le \kappa \left( 1 + \left| r_1 - r_2\right|^{\theta} \right), \quad \forall r_1 , r_2  \in (0, \infty),
	\end{equation}
where $\calN_{\Lambda_2}$ is the counting function associated to $\Lambda_2$.

Notice that~\eqref{H5} implies the weak-gap condition: for any $\rho >0$ and any $x>0$
\[
\sharp \left[ \Lambda_2 \cap (x - \rho/2 , x + \rho/2 ) \right] 
\le \calN_{\Lambda_2}\left( x+ \frac{\rho}{2} \right) - \calN_{\Lambda_2} \left( x- \frac{\rho}{2} \right) 
\leq \kappa ( 1 + \rho^\theta). 
\]
In all what follows we consider  $\rho >0$ fixed and $ p \in \N $ such that
	\begin{equation}\label{H3}
\sharp \left[ \Lambda_2 \cap (x - \rho/2 , x + \rho/2 ) \right] \le p , \qquad \forall x > 0.
	\end{equation}
Following~\cite[Proposition 2.2]{BM:23}, if the sequence $\Lambda_2 $ satisfies~\eqref{H3} there exists a countable family $\left\{ G_k \right\}_{k\geq 1}$ of disjoint subsets of $\Lambda_2$ satisfying
	\begin{gather}
\Ds  \Lambda_2 = \bigcup_{k \ge 1 }G_k, \quad G_k = \left\{ \lambda_k^{(1)}, \dots , \lambda_k^{(g_k)} \right\}, \quad  \lambda_k^{(1)} <  \lambda_k^{(2)} < \cdots < \lambda_k^{(g_k)} , \label{Gk} \\ 
	\noalign{\smallskip}
\quad \Ds  g_k \le p, \quad \max G_k - \min G_k \le \rho , \quad C(p, \rho) \le  \min G_{k+1} - \max G_k , \quad \forall k \ge 1, \label{Gk2}
	\end{gather}
with $C(p, \rho)$ a new positive constant only depending on $p$ and $\rho$. 

\end{description}

Let $\calB_1 := \left\{ \psi_m \right\}_{m \ge 1 }$ be an orthonormal basis of $L^2 (\Omega_1)$.

Let $\calU_2$ be a Hilbert space with inner product and associated norm respectively denoted by $(\cdot, \cdot)_{\calU_2}$ and $\| \cdot \|_{\calU_2}$.
Let us also consider an operator 
	$$
\frakC_2 \in \calL ( H^2 (0, \pi) \cap H_0^1 (0, \pi), \calU_2).
	$$
Let us assume that, associated to the sequence $\Lambda_2$ given by~\eqref{Gk}, we have $\calB_2 \subset L^2(0, \pi)$, a family of $L^2(0, \pi)$, given by
	\begin{equation}\label{B2}
\calB_2 := \bigcup_{k \ge 1} B_k , \quad B_k = \left\{ \phi_k^{(1)}, \dots , \phi_k^{(g_k)} \right\}  \subset H^2 (0, \pi) \cap H_0^1 (0, \pi) , \quad \forall k \ge 1 ,
	\end{equation}
and satisfying 
	\begin{equation}\label{H6}
\frakC_2 \phi_k^{(j)} \neq 0, \quad \forall k \ge 1 \hbox{ and } j : 1 \le j \le g_k.  
	\end{equation}
Finally, we assume that there exists $\beta > 0$ and $ \vartheta \in (0, 1) $ such that the following inequality
	\begin{equation}\label{H2}
\int_{\Omega_1 } \left\| \sum_{\mu_{m}^{\vartheta}  \leq \lambda } G_m  \psi _{m}(x' )  \right\|_{\calU_2}^{2}\, dx' 
\leq e^{\beta \lambda } \int_{\omega }\left \| \sum_{\mu_{m}^{\vartheta} \leq \lambda } G_m \psi _{m}(x') \right \|_{\calU_2} ^{2}\, dx' ,
	\end{equation}
holds for any $ \lambda \in (0, \infty)$ and any $G_m \in \operatorname{Span} \left\{  \frakC_2 \phi_k^{(j)} \: : \: k \geq 1, \: 1 \leq j \leq g_k \right\}$ for any $m \geq 1$. 
%

\begin{remark} \label{Rem:Inegalite_Spectrale}
For applications to the study of null controllability for parabolic problems the main settings we have in mind are boundary control or distributed control that is to say, respectively,
\[
\calU_2 = \R \qquad \text{ and } \qquad \frakC_2 \phi_k^{(j)} = -\left( \phi_k^{(j)} \right)'(0)
\]
or
\[
\calU_2 = L^2(0,\pi)  \qquad \text{ and } \qquad \frakC_2 \phi_k^{(j)} = \mathbf{1}_{(a,b)} \phi_k^{(j)}
\]
with $0 \leq a < b \leq \pi$. Let us mention that, in these two settings, the validity of inequality~\eqref{H2} is a direct consequence of the following so-called spectral inequality associated with $\calB_1$:
	\begin{equation}\label{H2'}
\int_{\Omega_1 }\left \vert \sum_{\mu_{m}^{\vartheta} \leq \lambda }b_{m}\psi _{m}(x' )\right \vert ^{2}\, dx' \leq e^{\beta \lambda } \int_{\omega }\left \vert \sum_{\mu_{m}^{\vartheta} \leq \lambda }b_{m}\psi _{m}(x' )\right \vert ^{2}\, dx' ,
	\end{equation}
for any $ \lambda \in (0, \infty)$ and $ \left \{ b_m \right\}_{m \ge 1} \in \ell^2$. This will be detailed on actual examples in Section~\ref{s5}. 
\end{remark}

The above hypotheses can be written more concisely as follows:

\begin{assumption}  \label{A1}
We have two positive real sequences $\Lambda_1 $ and $\Lambda_2$, an orthonormal basis $\calB_1$ of $L^2 (\Omega_1)$, a sequence $\calB_2$ in $L^2 (0, \pi)$, a Hilbert space $\calU_2$, and an operator $\frakC_2 \in \calL ( H^2 (0, \pi) \cap H_0^1 (0, \pi), \calU_2)$ such that
	\begin{equation*}
	\left\{
	\begin{array}{l}
\Ds  \Lambda_1  \hbox{ satisfies~\eqref{H1} with } \kappa_1, \theta_1 > 0; \\
	\noalign{\smallskip}
 \Ds \Lambda_2 \hbox{ satisfies~\eqref{H5} and~\eqref{H3} with } p \in \N, \rho, \kappa > 0 \hbox{ and } \theta \in (0, 1); \\
	\noalign{\smallskip}
\Ds \calB_2 \hbox{ is given by~\eqref{B2} and satisfies~\eqref{H6} };  \\
	\noalign{\smallskip}
\hbox{the spectral inequality~\eqref{H2} holds with } \beta >0 \hbox{ and } \vartheta \in (0, 1).
	\end{array}
	\right.
	\end{equation*}
The sequence $\Lambda_2$ is labeled accordingly to the grouping~\eqref{Gk}, \eqref{Gk2}.
\end{assumption}

\begin{remark}
Notice that assumption~\eqref{H5} also implies
	\begin{equation}\label{H4}
\calN_{\Lambda_2} (r) \le 2\kappa r^{\theta} , \quad \forall r \in (1, \infty).
	\end{equation}
From Weyl's law, this explains why the tangential variable $x$ is one dimensional in our study. 
To fit into the framework of~\cite{BM:23}, we will write $\Lambda_2 \in \calL (p, \rho, \theta , \kappa)$, where
	\begin{equation}\label{Lclass}
\calL (p, \rho, \theta , \kappa) := \left\{ \Lambda : \Lambda \subset (0, \infty) \mbox{ is a sequence satisfying \eqref{H3} and \eqref{H5}} \right\}. 
	\end{equation}
\end{remark}

\subsection{The main result}\label{ss22}

With the notations of Section~\ref{ss21} (see Assumption~\ref{A1}), let us consider
	\begin{equation}\label{ekm}
	\left\{
	\begin{array}{l}
e_k^{(j)} (t) = e^{- \lambda_k^{(j)} t} \frakC_2 \phi_k^{(j)} \in \calU_2 , \quad \forall k \ge 1 \hbox{ and } j:  1\le j \le g_k,     \\
  	\noalign{\smallskip}
e_{m, k}^{(j)}  (t) = e^{- \lambda_{m, k}^{( j)} t} \frakC_2 \phi_k^{(j)} \in \calU_2, \quad \forall m, k \geq 1 \hbox{ and } j:  1\le j \le g_k, 
	\end{array}
	\right.
	\end{equation}
where $t \in (0,T)$ and $\lambda _{m,k}^{ ( j ) }$ is given by
	\begin{equation}\label{spec}
\lambda_{m, k}^{(j)} = \mu_m + \lambda_k^{( j )} ,  \quad (m, k )\in \N^2, \quad 1 \le j \le g_k. 
	\end{equation}
We will also use the sequence $ \calF := \left\{ F_{m,k}^{(j)}\right\}_{\substack{ m, k \ge 1\\  1 \le j \le g_k } }$ of elements of $\calU_2$ given by
	\begin{equation}\label{Fmk}
F_{m,k}^{(j)} (t, x') := e_{m, k}^{(j)} (t) \psi_{m} ( x') = e^{-\lambda _{m,k}^{ ( j ) }t} \psi_{m} ( x') \frakC_2 \phi_k^{(j)} , \quad ( t,x'   ) \in (0,T) \times \Omega_1 , 
	\end{equation}
for any $m,k\geq 1$ and $j : 1 \le j \le g_k $, where $\lambda _{m,k}^{ ( j ) }$ is given by~\eqref{spec}. 
For any $k \geq 1$, we define the matrix
\begin{equation} \label{Def:M_k}
M_k = \sum_{\ell=1}^{g_k} \operatorname{Gram}_{\calU_2} \left(
\delta_{k,\ell}^1 \frakC_2 \phi_k^{(1)}, \dots, \delta_{k,\ell}^{g_k} \frakC_2 \phi_k^{(g_k)}
\right)
\end{equation}
where
\begin{equation} \label{Def:Delta_k}
\left\{
\begin{aligned}
\delta_{k, 1}^j &= 1,  \qquad \forall 1 \leq j \leq g_k,
\\
\delta_{k, \ell}^j &= \prod_{i= 1}^{\ell-1} \left( \lambda_{k}^{(j)} - \lambda_{k}^{(i)} \right).
\end{aligned}
\right.
\end{equation}
From~\cite[Proposition 13]{BM:23}, we have that assumption~\eqref{H6} implies that the matrix $M_k$ is invertible. 

The main result of this paper establishes the existence of a biorthogonal family associated to the sequence $\left\{ F_{m,k}^{(j)}\right\}_{\substack{ m, k \ge 1 \\  1 \le j \le g_k } }$ in $L^{2}( (0, T) \times  \omega ; \calU_2) $ and provides an estimate of the norm of its elements. It reads as follows:


\begin{theorem}\label{main}
Let us assume that $(\Lambda_1, \calB_1, \Lambda_2, \calB_2, \calU_2 , \frakC_2)$ satisfies Assumption~\ref{A1}. Then, there exists a constant $\calC >0$, only depending on $p$, $\rho$, $\theta$, $\kappa$, $\beta$, $\vartheta$, $\theta_1$ and $\kappa_1$, such that for any $T >0$, the sequence $\left\{ F_{m,k}^{(j)}\right\}_{\substack{ m, k \ge 1\\  1 \le j \le g_k } }$ (see~\eqref{Fmk}) admits a biorthogonal family $\left\{ Q_{m,k}^{(j)}\right\}_{\substack{ m, k \ge 1\\  1 \le j \le g_k } }$ in $L^{2}( ( 0,T ) \times \omega ; \calU_2)$, \textit{i.e.}, such that for any $ m,n \ge 1$, any $k , \ell \geq 1$, any $j: 1 \le j \le g_k $ and any $i: 1 \le i \le g_\ell$ we have 
\[
\int_0^T \int_\omega \psU{ Q_{m,k}^{(j)}(t,x') }{ F_{n,\ell}^{(i)}(t,x') } dx' dt = \delta_{mn} \delta_{k \ell} \delta_{ji},
\] 
that satisfies
	\begin{equation}\label{main-esti}
\Ds \left\| Q_{m,k}^{(j)} \right\|^2_{L^2 ((0, T) \times \omega; \calU_2) } \le \calC \exp \left(\frac{ \calC }{T^b } + \frac{ \calC} {T^{ \theta '} } \right) \exp \left( \calC \left[ \lambda_{m, k}^{(1)} \right]^{ \frac{ b}{1 + b}} + \calC \left[ \lambda_{m , k }^{(1)} \right]^\theta \right)  \left( M_k^{-1} \right)_{j,j},
	\end{equation}
for any $ m, k \ge 1$ and $j: 1 \le j \le g_k $, where $M_k$ is the matrix defined in~\eqref{Def:M_k} and $\theta '$ is given by 
	\begin{equation}\label{theta}
 \theta ' = \frac \theta {1 - \theta} \in (0, \infty),
	\end{equation}
and $b$ is given by and
	\begin{equation}\label{bcoef}
b := \vartheta \max \left\{ \frac 1{1 - \vartheta} , \frac 1{1 - \theta} \right\} . 
	\end{equation}
\end{theorem}


The proof of Theorem~\ref{main} will be done in the next two sections.
First, the idea consists in proving that the sequence $\calF$ (see~\eqref{Fmk}) has a biorthogonal family in $ L^{2}( ( 0,T ) \times \Omega_{1}; \calU_2)$. In this step we will use that the set $\left\{ \psi_m \right\}_{m \ge 1}$ is an orthonormal basis of $L^2(\Omega_1)$ (see Section~\ref{ss41}). 

Then, the main argument in the proof of Theorem~\ref{main} is the following one: we define the restriction operator from the closed subspace of $L^{2}_\rho ((0,T) \times \Omega_1 ; \calU_2)$ (with an appropriate weight function which blows up near $t= 0$) spanned by $\calF$ into $E^\omega$, the subspace of $L^{2}((0,T)\times \omega ; \calU_2) $ spanned by the restriction to $\omega $ of the elements of $\calF$. We prove that this operator is a bi-continuous bijection between the two spaces (see Section~\ref{s3}) which allows to deduce that $\calF$ has a biorthogonal family in $ L^{2}( ( 0,T ) \times \omega; \calU_2)$ (see Section~\ref{ss42}). In fact, this biorthogonal family $\left\{ Q_{m,k}^{(j)}\right\}_{\substack{ m, k \ge 1\\  1 \le j \le g_k } }$ belongs to the space $E^\omega$ and, in consequence, is unique and optimal.


\begin{remark}
From the expression of the constant $b$ (see~\eqref{bcoef}), we deduce:
\begin{itemize}
\item If $\vartheta \le \theta$ then $b = \frac \vartheta{1- \theta} \le \theta '$ and $\frac{ b}{1 + b} \le \theta$. In this case, inequality~\eqref{main-esti} becomes
	\begin{equation*}
\Ds \left\| Q_{m,k}^{(j)} \right\|^2_{L^2 ((0, T) \times \omega; \calU_2) } \le \calC \exp \left( \frac{ \calC} {T^{ \theta '} } \right) \exp \left( \calC \left[ \lambda_{m , k }^{(1)} \right]^\theta \right)  \left( M_k^{-1} \right)_{j,j},
	\end{equation*}
for a new positive constant $\calC $ and any $ m, k \ge 1$ and $j: 1 \le j \le g_k $.
\item If $\vartheta \ge \theta$ then $b = \frac \vartheta{1- \vartheta}  \ge \theta '$ and $ \theta \le \frac{ b}{1 + b} $. In this case, inequality~\eqref{main-esti} becomes
	\begin{equation*}
\Ds \left\| Q_{m,k}^{(j)} \right\|^2_{L^2 ((0, T) \times \omega; \calU_2) } \le \calC \exp \left(\frac{ \calC }{T^b } \right) \exp \left( \calC \left[ \lambda_{m, k}^{(1)} \right]^{ \frac{ b}{1 + b}} \right) \left( M_k^{-1} \right)_{j,j},
	\end{equation*}
for a new positive constant $\calC $ and any $ m, k \ge 1$ and $j: 1 \le j \le g_k $.
\end{itemize}
\end{remark}
\begin{remark}
The formulation of assumption~\eqref{H2} can be compared with~\cite[Assumption (6)]{M:10}. Roughly stated, in~\cite{M:10}, the author proves that 
\begin{itemize}
\item[$\bullet$] if observability holds in any time for a reference operator
\item[$\bullet$] and the considered observation operator and this reference operator satisfy a spectral inequality~\cite[Assumption (6)]{M:10}  
\end{itemize}
then observability holds in any time for the considered operator. 

Theorem~\ref{main} can somehow be seen as the analogous of~\cite{M:10} for biorthogonal families. Indeed our assumptions on $\Lambda_2$ implies the existence of biorthogonal families in $L^2((0,T) \times \Omega_1 ; \calU_2)$ and the spectral inequality~\eqref{H2} allows to transfer it to a biorthogonal families in $L^2((0,T) \times \omega ; \calU_2)$. 

The geometrical setting we consider is less general than the one of~\cite{M:10} but our assumptions are weaker since the spectral inequality is only assumed for the eigenvectors of the transverse operator. Also our results allow to study null controllability for a fixed given final time $T>0$ and is not limited to situations where null controllability holds at any time. Thus our results are not contained in~\cite{M:10} and conversely they do not completely cover the setting of~\cite{M:10}. Nevertheless, our method allows us to obtain some of the results of~\cite{M:10}. This is illustrated in Appendix~\ref{Annexe2} where, for simplicity, we only treat a particular case. The method can be extended to more general cases.
\end{remark}

\subsection{Notation} \label{ss23}

We gather in this section some notation that will be used throughout this article. 

\begin{description}
\item[Divided differences.]
\leavevmode\par 
In all this manuscript the notation $f[x_1, \dots, x_n]$ denotes divided differences. 
For pairwise distinct $x_1, \dots, x_n \in \R$ and $f_1, \dots, f_n$ in a real vector space, the divided differences are defined by
\[
f[x_i]=f_i, \qquad \forall i \in \{1 , \dots, n\}
\]
and then recursively 
for any $k \in \{2, \dots, n \}$, for any pairwise distinct $i_1, \dots, i_k \in \{1, \dots,n\}$, by 
\[
f[x_{i_1}, \dots, x_{i_k}] = \frac{f[x_{i_1}, \dots, x_{i_{k-1}}]-f[x_{i_2}, \dots, x_{i_k}]}{x_{i_1}-x_{i_k}}.
\]
We will use the mean value theorem for divided differences. It states that if $f$ is a $n$ times differentiable function then, for any $k \in \{2, \dots, n \}$, for any pairwise distinct $i_1, \dots, i_k \in \{1, \dots,n\}$, there exists $z \in \operatorname{Conv}\left\{ x_{i_1}, \dots, x_{i_k} \right\}$ such that
\[
f[x_{i_1}, \dots, x_{i_k}] = \frac{f^{(k-1)}(z)}{(k-1)!}
\]
where $f_i = f(x_i)$. 

\item[Linear combination.]
\leavevmode\par 
Recall that $F_{m,k}^{(j)}$ is defined in~\eqref{Fmk} by
\[
F_{m,k}^{(j)} (t, x') = e^{-\lambda _{m,k}^{ ( j ) }t} \psi_{m} ( x') \frakC_2 \phi_k^{(j)} \in \calU_2,
\]
for $(t,x') \in (0,T) \times \Omega_1$, $m,k\geq 1$ and $j : 1 \le j \le g_k $. In this article we will often deal with elements of
\[
\mathrm{span} \, \left \{ F_{m, k}^{(j)} : k,m\geq 1, \ j : 1 \le j \le g_k \right\}.
\]
These finite combinations of the functions $F_{m,k}^{(j)}$ will be denoted by
	\begin{equation} \label{PN}
P_{N} (t, x')  := \sum_{\mu _{m}^\vartheta, k \leq N} \sum_{j= 1 }^{g_k}  a_{m,k}^{ (j , N ) } F_{m,k}^{(j)} (t, x')  = \sum_{\mu _{m}^\vartheta \leq N} G_{m}^{(N)} (t) \psi_m (x '), 
	\end{equation}
where $(t, x' ) \in (0, T) \times \Omega_1$, $\lambda _{m,k}^{ ( j ) }$ is given in~\eqref{spec}, 
	\begin{equation}\label{G_m}
G_{m}^{(N) } (t) : = \sum_{k=1}^{N} g_{m,k}^{ ( N ) } (t),  \quad  g_{m,k}^{ ( N ) } (t) := \sum_{j= 1 }^{g_k} a_{m,k}^{ ( j, N ) } e^{-\lambda _{m,k}^{ ( j ) }t} \frakC_2 \phi_k^{(j)} \in \calU_2  ,
	\end{equation}
and $a_{m,k}^{ (j , N ) } \in \R$, for any $k, m \ge 1:  k, \mu_m^\vartheta \le N$ and $1 \le j \le g_k$. 
\item[Weights and functional spaces.]
\leavevmode\par
For any $\alpha > 0$, let us introduce the function
	\begin{equation}\label{eta}
\eta _{\alpha } ( x ' ) =\left \{ 
	\begin{array}{ll}
0 & \mathrm{if} \ x' \in  \omega \\ 
	\noalign{\smallskip}
\alpha \beta & \mathrm{if} \ x' \in \Omega_1 \backslash  \omega,
	\end{array}
	\right.
	\end{equation}
where $\beta >0$ is the constant appearing in~\eqref{H2}. 
We also introduce the Hilbert spaces 
	\begin{equation} \label{L2etaalpha}
L_{\eta _{\alpha }}^ {2}\left( ( 0,T ) \times \Omega_1  ; \calU_2\right) := \left\{ f : \int_{0}^T \!\! \int_{\Omega_1} e^{-\frac{\eta _{\alpha } ( x' ) }{t^b} }  \normeU{ f ( t, x' ) }^2 \,dt \, dx' < \infty \right\},
	\end{equation}
with $b >0$ given by~\eqref{bcoef}. This space is equipped with the scalar product 
	$$
\left(  f , g \right)_{L^2_{\eta_{\alpha }}((0,T) \times \Omega_1 ; \calU_2)} = \int_{0}^{T}\! \! \int_{\Omega_1 } e^{-\frac{\eta _{\alpha } ( x' ) }{t^b } } \psU{f}{g} \, dt \, dx',  \quad \forall f, g \in L_{\eta _{\alpha }}^{2} ( ( 0,T ) \times \Omega_1 ; \calU_2). 
	$$
We can now define the Hilbert spaces: 
	\begin{equation}\label{eta_alpha}
	\left\{
	\begin{array}{c}
\Ds E_{\eta _{\alpha }} =\overline{\mathrm{span} \, \left \{ F_{m, k}^{(j)} : k,m\geq 1, \ j : 1 \le j \le g_k \right\}  }^{L_{\eta _{\alpha }}^{2} ((0, T) \times \Omega_1 ; \calU_2)} , \\
	\noalign{\smallskip}
\Ds E^\omega =\overline{\mathrm{span}\, \left \{ F_{m, k}^{(j)} |_\omega : k,m\geq 1, \ j : 1 \le j \le g_k  \right\}  }^{L^{2} ((0, T) \times \omega ; \calU_2) } ,
	\end{array}
	\right.
	\end{equation}
where $F_{m, k}^{(j)}$ is given in~\eqref{Fmk}.
\end{description}


\section{The restriction operator}\label{s3}
In this section we provide the main idea used in the proof of Theorem~\ref{main}. We will prove (see Theorem~\ref{tmainrest}) that the restriction operator $1_\omega$ is a bi-continuous bijection between the closed subspace of $L^{2}_{\eta_\alpha} (( 0,T) \times  \Omega_1 ; \calU_2)$ spanned by $\left\{ F_{m,k}^{(j)}\right\}_{\substack{ m, k \ge 1 \\  1 \le j \le g_k } }$ and the subspace of $L^{2}( ( 0,T) \times  \omega ; \calU_2) $ spanned by $\left\{ F_{m,k}^{(j)} |_\omega \right\}_{\substack{ m, k \ge 1 \\  1 \le j \le g_k } }$. 
Recall that we are using the notations and assumptions of Section~\ref{ss21} for the sequences $\Lambda_1$ and $\Lambda_2$.


The key point to study the restriction operator is the following result:


\begin{theorem}\label{tPN}
Let us assume that $(\Lambda_1, \calB_1, \Lambda_2, \calB_2, \calU_2 , \frakC_2)$ satisfies Assumption~\ref{A1}. Then, there exist $\alpha > 0$ (only depending on $p$, $\rho$, $\theta$, $\kappa$, $\beta$, $\vartheta$ and $\theta_1$) and $ \tau_ 0 \in (0, 1]$ (only depending on $p$, $\rho$, $\theta$, $\kappa$, $\beta$, $\vartheta$, $\theta_1$ and $\kappa_1$) such that
	\begin{equation}\label{f7}
\Ds  \int_{0}^T  \!\! \int_{ \Omega_1 }  e^{-\frac{\alpha \beta }{t^{ b} }} \normeU{ P_{N} (t,x') }^{2} \, dx'\,dt  \leq  6 \int_0^T   \!\! \int_{ \omega } \normeU{ P_{N}(t,x') }^{2} \, dx'\,dt , 
	\end{equation}
for any $ T \in (0, \tau_0]$, any $N \ge 1$ and any $P_N $ given by~\eqref{PN}, $\beta $ is the constant in~\eqref{H2} and $b$ is given by~\eqref{bcoef}.
\end{theorem}


The proof of Theorem~\ref{tPN} is postponed to Section~\ref{ss32}.


\begin{remark}\label{r1}
When $T >0$ is arbitrary, it is possible to prove a slightly different version of Theorem~\ref{tPN}. For the details, see Remark~\ref{r5} and Theorem~\ref{tPNp}.
\end{remark}

The weight appearing in the left-hand side of inequality~\eqref{f7} motivates the definition of the function $\eta_\alpha $ (see~\eqref{eta}). From this definition, it is clear that 
	$$
L^{2} ( (0, T) \times \Omega_1 ; \calU_2) \hookrightarrow L_{\eta _{\alpha }}^{2} ( (0, T) \times \Omega_1 ; \calU_2)
	$$
with continuous injection. On the other hand, if $\varphi \in L_{\eta _{\alpha }}^{2} ( (0, T) \times \Omega_1 ; \calU_2) $, we then have $\varphi |_\omega \in L^{2} ( (0, T) \times \omega ; \calU_2) $. Therefore, we can define the restriction operator 
%
	\begin{align} \label{Def:restriction}
\calR_{\omega} : L_{\eta _{\alpha }}^{2} ( (0, T) \times \Omega_1 ; \calU_2) &\to L^{2} ( (0, T) \times \omega ; \calU_2)
\\
\varphi &\mapsto \calR_{\omega}( \varphi  ) =\varphi |_\omega
	\notag
	\end{align}
which satisfies $\calR_{\omega} \in \calL \left( L_{\eta _{\alpha }}^{2} ( (0, T) \times \Omega_1 ; \calU_2) , L^{2} ( (0, T) \times \omega ; \calU_2) \right)$.
The main result of this section reads as follows:


\begin{theorem}\label{tmainrest}
Let us assume that $(\Lambda_1, \calB_1, \Lambda_2, \calB_2, \calU_2 , \frakC_2)$ satisfies Assumption~\ref{A1} and consider $ \tau_0 , \alpha > 0$, the constants provided by Theorem~\ref{tPN}, $b$ given by~\eqref{bcoef} and $\eta_\alpha$, the function defined in~\eqref{eta} with $\beta > 0$ given in~\eqref{H2}. Then, if $T \in (0, \tau_0]$, the operator $\calR_{\omega} $ satisfies
	\begin{equation}\label{dir_indir_ineq}
\left\| \varphi \right\|_{L^2_{\eta_{\alpha }}((0,T) \times \Omega_1 ; \calU_2)}^2 
\le 7 \left\| \calR_{\omega} (\varphi ) \right\|_{L^{2} ( (0, T) \times \omega ; \calU_2) }^2 
\le 7 \left\| \varphi \right\|_{L^2_{\eta_{\alpha }}((0,T) \times \Omega_1 ; \calU_2)}^2, \quad \forall \varphi \in E_{\eta _{\alpha }}  .
	\end{equation}
Moreover, $\calR_{\omega} \left( E_{\eta _{\alpha }}  \right) = E^\omega  $ and, therefore, $\calR_{\omega} \in \calL \left( E_{\eta _{\alpha }} , E^\omega  \right)$ is an isomorphism. 
\end{theorem}


\begin{proof}
Let us take $\alpha , \tau_0 >0$ provided by Theorem~\ref{tPN} and $b >0$ given by~\eqref{bcoef}. With these constants, it is possible to apply Theorem~\ref{tPN}. As a consequence, let us first prove that
	\begin{equation}\label{dir_indir_ineq0}
\left\| P_N \right\|_{L^2_{\eta_{\alpha }}((0,T) \times \Omega_1 ; \calU_2)}^2 
\le 7 \left\| \calR_{\omega} (P_N ) \right\|_{L^{2} ( (0, T) \times \omega ; \calU_2)}^2 
\le 7 \left\| P_N \right\|_{L^2_{\eta_{\alpha }}((0,T) \times \Omega_1 ; \calU_2)}^2, 
	\end{equation}
for any $N \ge 1 $ and $P_N$ given by~\eqref{PN}, with $a_{m,k}^{ (j , N ) } \in \R$, for any $k, m \ge 1:  k, \mu_m^\vartheta \le N$ and $1 \le j \le g_k$. Using Theorem~\ref{tPN}, we have
\begin{align*}
\Ds \left\| P_N \right\|_{L^2_{\eta_{\alpha }}((0,T) \times \Omega_1 ; \calU_2)}^2 
&= \int_{0}^T \!\! \int_{ \left( \Omega _{1} \setminus \omega \right) }  e^{-\frac{\alpha \beta }{t^b }} \normeU{ P_{N} (t,x') }^{2}  \, dx' \,dt + \int_{0}^T  \!\! \int_{\omega } \normeU{ P_{N}(t,x') }^{2}  \, dx' \,dt
\\
&\leq \int_{0}^T \!\! \int_{\Omega_1 } e^{-\frac{\alpha \beta }{t^b } } \normeU{P_{N} ( t,x' )}^{2} \, dx'\,dt +  \int_{0}^T  \!\! \int_{\omega } \normeU{ P_{N} (t,x')}^{2}  \, dx' \,dt
\\
&\leq 7  \int_{0}^T  \!\! \int_{\omega } \normeU{ P_{N}(t,x')}^{2}  \, dx'\,dt = 7 \left\| \calR_{\omega} (P_N ) \right\|_{L^{2} ( (0, T) \times \omega ; \calU_2) }^2.
\end{align*}
On the other hand, one has
	\begin{align} \notag
\Ds \left\| \varphi \right\|_{L^2_{\eta_{\alpha }}((0,T) \times \Omega_1 ; \calU_2)}^2 
&= \int_{0}^T \!\! \int_{ \left( \Omega _{1} \setminus \omega \right) }  e^{-\frac{\alpha \beta }{t^b }} \normeU{ \varphi(t,x')}^{2}  \, dx'\,dt + \int_{0}^T  \!\! \int_{\omega } \normeU{ \varphi(t,x') }^{2}  \, dx' \,dt
\\
&\ge \int_{0}^T  \!\! \int_{\omega} \normeU{ \varphi(t,x')}^{2}  \, dx'\,dt .
\label{f8}
	\end{align}
for any $\varphi \in  E_{\eta _{\alpha }} $. Thus,~\eqref{dir_indir_ineq0} holds for any $P_N$ given by~\eqref{PN}.

Let us prove that $\calR_{\omega} \left( E_{\eta _{\alpha }} \right) \subseteq E^\omega $. Indeed, first, we have 
	$$
\calR_{\omega} \in \calL \left( E_{\eta _{\alpha }} , L^{2} ( (0, T) \times \omega ; \calU_2) \right). 
	$$ 
Secondly, we also have
	$$
E^\omega =\overline{\mathrm{span}\, \left \{ \calR_{\omega} \left( F_{m, k}^{(j)} \right) : k,m\geq 1, \ j : 1 \le j \le g_k\right \} }^{L^{2} ((0, T) \times \omega ; \calU_2)}.
	$$
Thus, if $\varphi \in E_{\eta _{\alpha }} $ there exists a sequence $\left\{ P_N \right\}_{N \ge 1}$ ($P_N$ given by~\eqref{PN}) such that $ P_N \to \varphi$ in $L^2_{\eta_\alpha } ( (0,T) \times \Omega_1 ; \calU_2)$. In particular, $\calR_{\omega} \left(P_N \right) \in E^\omega $ and, from~\eqref{f8}, $\calR_{\omega} \left(P_N \right) \to \calR_{\omega}  (\varphi )$ in $L^{2} ((0, T) \times \omega ; \calU_2)$. We deduce therefore that $ \calR_{\omega}  (\varphi ) \in E^\omega$. 

Let us now prove the inclusion $ E^\omega \subseteq  \calR_{\omega} \left( E_{\eta _{\alpha }}  \right)$. To this end, let us consider $\psi \in E^\omega $. 
From the definition of this space, again, there exists a sequence $\left\{ P_N \right\}_{N \ge 1}$ ($P_N$ given by~\eqref{PN}) such that $ \calR_\omega (P_N ) \to \psi$ in $L^2 ( (0,T) \times \omega ; \calU_2)$. 
This implies that $ \left\{ \calR_\omega (P_N ) \right\}_{N \ge 1}$ is a Cauchy sequence in $L^2 ( (0,T) \times \omega ; \calU_2)$ and, from~\eqref{dir_indir_ineq0}, $ \left\{ P_N  \right\}_{N \ge 1}$ is also a Cauchy sequence in $L^2_{\eta_\alpha} ( (0,T) \times \Omega_1 ; \calU_2)$. 
Thus, there exists $\varphi \in E_{\eta _{\alpha }} $ such that $ P_N \to \varphi$ in $L^2_{\eta_\alpha } ( (0,T) \times \Omega_1 ; \calU_2)$.
Since $ \calR_\omega (P_N ) \to \calR_\omega (\varphi )$ in $L^2 ( (0,T) \times \omega ; \calU_2)$ (see~\eqref{f8}), we infer that $\psi = \calR_\omega (\varphi ) $. This completes the proof of the inclusion. 

Finally, inequality~\eqref{dir_indir_ineq} is a direct consequence of~\eqref{dir_indir_ineq0}. This ends the proof of the result.  
\end{proof}


\begin{remark}
As in Remark~\ref{r1}, it is possible to prove a new version of Theorem~\ref{tmainrest} valid for any $T > 0$. In this case the positive constant $ \alpha > 0$ depends on $p$, $\rho$, $\theta$, $\kappa$, $\beta$, $\vartheta$, $\theta_1$, $\kappa_1$ and $T$. The proof can be deduced using Theorem~\ref{tPNp} instead of Theorem~\ref{tPN} and following the same argument as in the proof of Theorem~\ref{tmainrest}. 
\end{remark}


The remaining part of this section is dedicated to prove Theorem~\ref{tPN}. The proof will use the preliminary results stated in the following section.

 
\subsection{Preliminary results} \label{ss31}
%


Let us start by stating and proving a technical result that will be used below:

\begin{lemma}\label{l1}
Let us consider $ \calS := \{ \nu_m \}_{m\geq 1} \subset (0, \infty)$ an increasing sequence satisfying
	\begin{equation*}
\calN_\calS (r) \le  \widetilde \calC r^q, \quad \forall r \in (0, \infty),
	\end{equation*}
for two constants $\widetilde \calC >0 $ and $q > 0$ (the counting function $\calN_\calS $ is defined in~\eqref{counting}). Then, there exists a positive constant $\widehat \calC$ (only depending on $\widetilde \calC$ and $q$) such that
	\begin{equation}\label{sum}
\Ds \sum_{\gamma < \nu_{m}} e^{-\sigma \nu _{m} } \le \widehat \calC \, \frac{ 1 + \left( \sigma \gamma \right)^{ q} }{\sigma^{ q }} e^{- \sigma \gamma} ,
	\end{equation}
for any $\sigma, \gamma > 0$.
\end{lemma}


\begin{proof}
Given $\sigma, \gamma > 0$, we can write
	\begin{equation*}
	\begin{array}{l}
\Ds \sum_{\gamma < \nu_{m}} e^{-\sigma \nu _{m} } = \int_\gamma^{\infty }e^{- \sigma r }\, d\calN_\calS  ( r )  
= \left[ e^{- \sigma r }\calN_\calS  ( r ) \right]_\gamma^{\infty } + \sigma \int_\gamma^{\infty } e^{-\sigma r }\calN_\calS (r) dr \\
	\noalign{\smallskip}
\Ds \phantom{\sum_{\gamma < \nu_{m}} e^{-\sigma \nu _{m} }} \leq \widetilde \calC \sigma \int_\gamma^{\infty } r^{ q } e^{-\sigma r}\, dr.
	\end{array}
	\end{equation*}
With the change of variable $s = \sigma (r - \gamma)$ in this last integral, we get:
	\begin{equation*}
\Ds \sum_{\gamma < \nu_{m}} e^{-\sigma \nu _{m} } \le \frac{\widetilde \calC }{\sigma^{ q }} e^{- \sigma \gamma } \int_0^\infty  \left( s + \sigma \gamma \right)^{ q} e^{- s }\, d s = \frac{\widetilde \calC }{\sigma^{ q }} e^{- \sigma \gamma } \left( \calJ_1 + \calJ_2 \right) , 
	\end{equation*}
where
	\begin{equation*}
\calJ_1 =  \int_0^1  \left( s + \sigma \gamma \right)^{ q} e^{- s }\, d s \le \left( 1 + \sigma \gamma \right)^{ q} \int_0^1 e^{- s }\, d s \le \left( 1 + \sigma \gamma \right)^{ q} ,
	\end{equation*}
and
	\begin{equation*}
\calJ_2 = \int_1^\infty  \left( s + \sigma \gamma \right)^{ q} e^{- s }\, d s \le \left( 1 + \sigma \gamma \right)^{ q}  \int_1^\infty   s^{ q} e^{- s }\, d s \le \left( 1 + \sigma \gamma \right)^{ q} \Gamma (q + 1),
	\end{equation*}
where $\Gamma (z)$ is the gamma function. Therefore,
	\begin{equation*}
\Ds \sum_{\gamma < \nu_{m}} e^{-\sigma \nu _{m} } \le {\widetilde \calC }  \left( 1 + \Gamma (q + 1) \right) \frac{\left( 1 + \sigma \gamma \right)^{ q}}{\sigma^{ q }} e^{- \sigma \gamma }. 
	\end{equation*}

Finally, taking into account the inequality
	$$
 ( 1 + x  )^{ q} \le \max \left\{1, 2^{q - 1} \right\}  \left( 1 + x^{ q} \right) , \quad \forall x \in [0, \infty),
	$$
we deduce the existence of a new constant $\widehat \calC$, only depending on $\widetilde \calC$ and $q$, for which~\eqref{sum} holds. This finishes the proof of the result. 
\end{proof}


As an intermediate tool, we will use the resolution of some block moment problems as developed in~\cite{BBM}. The refined estimates given in the following theorem are proved in~\cite[Appendix A, Theorem 46]{BM:23}. 

\begin{theorem}\label{tblock}
Let us fix $p \in \N$, $\rho, \kappa > 0$ and $\theta \in (0, 1)$. Then, there exists a constant $\widetilde C_0 > 0$ (only depending on $p$, $\rho$, $\theta$ and $\kappa$) such that for any $T >0$, any $\Lambda_2 \in \calL (p, \rho, \theta , \kappa)$ (which we will assume is given by~\eqref{Gk}) and any $\left\{ f_{k}^{(j)} \right\}_{\substack{ k \ge 1 \\  1 \le j \le g_k } } \subset \R$, there exists a family $\left\{ r_{k} \right\}_{k \geq 1} \subset L^2(0,T)$ satisfying
	\begin{equation*}
	\left \{ 
	\begin{array}{l}
\Ds \int_0^T e^{-  \lambda_{k}^{(j)}  t} r_{\ell } (t) \, dt = 0, \quad \forall k , \ell \ge 1 \hbox{ and } 1 \le j \le g_k : k \not= \ell , \\ 
	\noalign{\smallskip}
\Ds \int_0^T e^{-  \lambda_{k}^{(j)}  t} r_{k } (t) \, dt = f_{k}^{(j)} , \quad \forall k \ge 1 \hbox{ and } 1 \le j \le g_k ,
	\end{array}
	\right.
	\end{equation*}
and
	\begin{equation*}
\left\| r_{k}  \right\|_{L^2(0,T)} \leq \widetilde C_0 \exp \left( \frac{\widetilde C_0} {T^{ \theta '} } \right) e^{\widetilde C_0 \left[ \lambda_{k}^{(1)} \right]^\theta } \max_{1 \le j \le g_k } \left\{ \left| f \left[ \lambda_{k}^{(1)} , \dots , \lambda_{k}^{(j)}  \right] \right|  \right\}, \quad \forall k \ge 1,
	\end{equation*}
with $f \left[\lambda_{k}^{(j)} \right] = f_{k}^{(j)} $, for any $ k \ge 1$ and $j:1 \le j \le g_k$, and $\theta ' $ given in~\eqref{theta}. 
\end{theorem}

Actually, we will use the following moment problem. It will allow us to deal with the blow-up of the weight $t \mapsto \exp \left( \frac{\alpha \beta}{t^b} \right)$ near $t=0$ and to obtain uniform estimates with respect to $m$. 


\begin{corollary}\label{cblock}
Let us fix $p \in \N$, $\rho, \kappa > 0$ and $\theta \in (0, 1)$. Then, there exists a positive constant $\calC_0$ (only depending on $p$, $\rho$, $\theta$ and $\kappa$) such that for any $T >0$, $\eps \in (0, T/4)$, $m \ge 1$, $\Lambda_2 \in \calL (p, \rho, \theta , \kappa)$ (given by~\eqref{Gk}) and $\left\{ f_{m,k}^{(j)} \right\}_{\substack{ k \ge 1 \\  1 \le j \le g_k } } \subset \R$, there exists a family $\left\{ r_{m,k}^\eps \right\}_{k \geq 1} \subset L^2(0,T)$ satisfying
	\begin{equation}\label{fblock}
	\left \{ 
	\begin{array}{l}
\Ds \int_0^T e^{- \lambda_{m, k}^{(j)}  t} r_{m, \ell }^\eps (t) \, dt = 0, \quad \forall k , \ell, m \ge 1 \hbox{ and } 1 \le j \le g_k : k \not= \ell , \\ 
	\noalign{\smallskip}
\Ds \int_0^T e^{- \lambda_{m, k}^{(j)}  t} r_{m, k }^\eps (t) \, dt = f_{m,k}^{(j)} , \quad \forall k,m \ge 1 \hbox{ and } 1 \le j \le g_k,
	\end{array}
	\right.
	\end{equation}
($\lambda_{m, k}^{(j)} $ is given in~\eqref{spec}) and 
	\begin{equation}\label{f2block}
	\left\{
	\begin{array}{l}
\Ds r^{\eps}_{m,k}  \equiv 0 \quad \mbox{in} \quad  ( 0,\varepsilon ), \\
	\noalign{\smallskip}
\Ds \left\| r_{m,k}^\eps  \right\|_{L^2(0,T)} \leq \calC_0 \exp \left( \frac{\calC_0} {T^{ \theta '} } \right) e^{\calC_0 \left[ \lambda_{m,k}^{(1)} \right]^\theta } e^{\eps \mu_m} \calK_{m,k}^\eps,  
	\end{array}
	\right.
	\end{equation}
for any $(m, k) \in \N^2 $, where
	\begin{equation}\label{f3block}
\calK_{m,k}^\eps = \max_{1 \le j \le g_k } \left\{ \left| f_{m,\eps} \left[ \lambda_{k}^{(1)} , \dots , \lambda_{k}^{(j)}  \right] \right|  \right\}, \quad \forall k \ge 1,
	\end{equation}
with $f_{m,\eps} \left[ \lambda_{k}^{(j)} \right] = e^{\eps \lambda_{k}^{(j)} } f_{m, k}^{(j)} $, for any $ k \ge 1$ and $j:1 \le j \le g_k$, and $\theta ' $ given in~\eqref{theta}. 
\end{corollary}



\begin{proof}
Let us consider the positive constant $\widetilde C_0 $ associated to $p \in \N$, $\rho, \kappa > 0$ and $\theta \in (0, 1)$ provided by Theorem~\ref{tblock}. Let us also consider $\Lambda_2 \in \calL (p, \rho, \theta , \kappa)$ (given by~\eqref{Gk}) and $\mu_m \in (0, \infty)$. Thus, the sequence
	$$
\Lambda^{(m )} = \mu_m + \Lambda_2 = \left\{ \lambda_{m,k}^{(j)}\right\}_{ k \ge 1,  1 \le j \le g_k } 
	$$
satisfies $\Lambda^{(m )} \in \calL (p, \rho, \theta , \kappa)$. Indeed, condition~\eqref{H3} holds for the parameters $\rho >0$ and $p \in \N$. On the other hand, it is not difficult to check 
	$$
\calN_{\Lambda^{(m )}} (r) = 
	\left\{
	\begin{array}{ll}
\Ds 0 & \hbox{ if } r \in (0, \mu_m], \\
	\noalign{\smallskip}
\Ds  \calN_{\Lambda_2} (r - \mu_m )& \hbox{ if } r \in (\mu_m, \infty),
	\end{array}
	\right.
	$$
and, therefore, $\Lambda^{(m )}$ satisfies~\eqref{H4} for the parameters $\kappa$ and $\theta \in (0,1)$. Finally, let us check condition~\eqref{H5} for $\Lambda^{(m )} $. This condition is direct if $r_1, r_2 \in (0, \mu_m]$. If $r_1, r_2 \in (0, \infty) $ are such that $r_1 \le \mu_m < r_2 $, then
	$$
\left| \calN_{\Lambda^{(m)}} (r_1) - \calN_{\Lambda^{(m)}} (r_2) \right| = \calN_{\Lambda_2} (r_2 - \mu_m ) \le \kappa (r_2 - \mu_m )^\theta \le \kappa \left( 1 + \left( r_2 - r_1 \right)^{\theta} \right).
	$$
Now, if $r_1,  r_2 \in (\mu_m , \infty)  $ one has:
	$$
\left| \calN_{\Lambda^{(m)}} (r_1) - \calN_{\Lambda^{(m)}} (r_2) \right| = \left| \calN_{\Lambda_2} (r_1 - \mu_m) - \calN_{\Lambda_2} (r_2- \mu_m ) \right|  \le \kappa \left( 1 + \left| r_2 - r_1 \right|^{\theta} \right).
	$$

If we fix $T >0$, $\eps \in (0, T/4)$ and a sequence $\left\{ f_{m,k}^{(j)} \right\}_{\substack{ k \ge 1 \\  1 \le j \le g_k } } \subset \R$, we can apply Theorem~\eqref{tblock} to the sequence $\Lambda^{(m)}$ and obtain the existence of a family $\left\{ \widehat r_{m,k}^\eps \right\}_{k \geq 1} \subset L^2(0,T - \eps)$ satisfying
	\begin{equation*}
	\left \{ 
	\begin{array}{l}
\Ds \int_0^{T-\eps} e^{- \lambda_{m, k}^{(j)}  t} \widehat r_{m, \ell }^\eps (t) \, dt = 0, \quad \forall k , \ell, m \ge 1 \hbox{ and } 1 \le j \le g_k : k \not= \ell , \\ 
	\noalign{\smallskip}
\Ds \int_0^{T-\eps} e^{- \lambda_{m, k}^{(j)}  t} \widehat r_{m, k }^\eps (t) \, dt = e^{\eps \mu_m }e^{\eps \lambda_k ^{(j)}} f_{m,k}^{(j)} , \quad \forall k,m \ge 1 \hbox{ and } 1 \le j \le g_k,
	\end{array}
	\right.
	\end{equation*}
and
	\begin{equation*}
	\begin{split}
\left\| \widehat r_{m,k}^\eps  \right\|_{L^2(0,T -\eps )} & \leq \widetilde C_0 \exp \left( \frac{\widetilde C_0} {(T-\eps)^{ \theta '} } \right) e^{\widetilde C_0 \left[ \lambda_{m,k}^{(1)} \right]^\theta } e^{\eps \mu_m} \calK_{m,k}^\eps \\
	\noalign{\smallskip}
&\leq \widetilde C_0 \exp \left( \frac{ 4^{ \theta '} \widetilde C_0} {(3T)^{ \theta '} } \right) e^{\widetilde C_0 \left[ \lambda_{m,k}^{(1)} \right]^\theta } e^{\eps \mu_m} \calK_{m,k}^\eps,\quad \forall (m,k) \in \N^2,
	\end{split}
	\end{equation*}
with $\calK_{m,k}^\eps$ given in~\eqref{f3block}. Finally, it is not difficult to check that the function $r_{m, k}^\eps$ given by:
	$$
r_{m, k}^\eps ( t ) =
	\left\{
	\begin{array}{ll}
\Ds 0, & \mbox{if } t \in (0, \eps], \\
	\noalign{\smallskip}
\Ds \widehat r_{m,k}^\eps (t - \eps) , & \mbox{if } t \in (\eps , T),
	\end{array}
	\right.
	$$
satisfies~\eqref{fblock} and~\eqref{f2block} for $\calC_0 = \left( \frac 43 \right)^{ \theta '} \widetilde C_0 $. This ends the proof. 
\end{proof}


Using the resolution of block moment problems given by Corollary~\ref{cblock} we obtain the following estimate. 


\begin{lemma}\label{l2}
Let us consider $(\Lambda_1, \calB_1, \Lambda_2, \calB_2, \calU_2 , \frakC_2)$ satisfying Assumption~\ref{A1}. Let us also consider $\calC_0 >0$ and $\beta >0$ the constants provided, resp., by Corollary~\ref{cblock} and inequality~\eqref{H2}. Then, for any $N\geq 1$, $T > 0$, $\eps \in (0, T/4)$, $\tau \in (\eps , T)$, $\alpha, b > 0 $ and $(m,k) \in \N^2: \mu_{m}^\vartheta, k \leq N$, if $P_{N}$ is given by~\eqref{PN}, we have:
	\begin{equation}\label{fl2}
\normeU{ g_{m,k}^{ ( N ) } ( \tau ) } \leq \widetilde \calH ( \eps, T) e^{\calC_0 \left[ \lambda_{m,k}^{(1)} \right]^\theta }   e^{- (\tau - \eps) \lambda_{m,k}^{( 1 )} }  \left( \int_{0}^T  \!\! \int_{\Omega _{1}} e^{-\frac{\alpha \beta }{t^{ b }} } \normeU{ P_{N} (t, x') }^{2} \, dx'\,dt \right)^{1/2},
	\end{equation}
where  $\lambda_{m,k}^{( 1 )}$, $g_{m,k}^{ ( N ) }$, and $\theta '$  are respectively given in~\eqref{spec},~\eqref{G_m} and~\eqref{theta}, and
	\begin{equation*}
\widetilde \calH (\eps, T) :=  \calC_0  \max \{ 1, T^{p - 1}\}  \exp \left( {\frac{\alpha \beta }{2 \varepsilon^{ b } }} + \frac{\calC_0} {T^{ \theta '} } \right).  
	\end{equation*}
\end{lemma}



\begin{proof}
Fixed $T > 0$, $\eps \in (0, T/4 )$ and $\tau \in (\eps , T)$, let us consider the constant $ \calC_0 > 0$ and the family $\left\{ r_{m,k}^{\tau , \eps} \right\}_{k \geq 1} \subset L^2(0,T)$ provided by Corollary~\ref{cblock} associated with $ f_{m,k}^{(j)} = e^{ - \lambda_{m,k}^{(j)} \tau }$. On the other hand, we define:
	\begin{equation*}
R_{m,k}^{\tau , \eps} (t, x' ) := r_{m,k}^{\tau , \eps} (t) \psi_m (x'), \quad \forall (m, k) \in \N^2 , \quad (t, x') \in (0, T) \times \Omega_1. 
	\end{equation*}
Let us first note that, from the properties of $ r_{m,k}^{\tau , \eps}$ (see Corollary~\ref{cblock}), we deduce:
	\begin{equation*}
	\left\{
	\begin{array}{l}
\Ds R_{m, k }^{\tau , \eps}  \equiv 0 \quad \mbox{in} \quad  ( 0,\varepsilon ) \times \Omega_1, \\
	\noalign{\smallskip}
\Ds \left\| R_{m, k }^{\tau , \eps} \right\|_{L^2((0, T) \times \Omega_1)} \leq \calC_0 \exp \left( \frac{\calC_0} {T^{ \theta '} } \right) e^{\calC_0 \left[ \lambda_{m,k}^{(1)} \right]^\theta } e^{\eps \mu_m} \calK_{k}^{\tau, \eps},  
	\end{array}
	\right.
	\end{equation*}
for any $(m, k) \in \N^2 $, where
	\begin{equation*}
\calK_{m, k}^{\tau, \eps} = \max_{1 \le j \le g_k } \left\{ \left| f_{m,\eps}^\tau \left[ \lambda_{k}^{(1)} , \dots , \lambda_{k}^{(j)}  \right] \right|  \right\}, \quad \forall k \ge 1,
	\end{equation*}
with $f_{m,\eps}^\tau \left[ \lambda_{k}^{(j)} \right] = e^{\eps \lambda_{k}^{(j)} } e^{- \lambda_{m, k}^{(j)} \tau }$, for any $ k \ge 1$ and $j:1 \le j \le g_k$. Using the expression~\eqref{spec}, we can write
	\begin{equation*}
\calK_{m, k}^{\tau, \eps} = e^{- \mu_m \tau} \max_{1 \le j \le g_k } \left\{ \left| f_{\tau,\eps} \left[ \lambda_{k}^{(1)} , \dots , \lambda_{k}^{(j)}  \right] \right|  \right\}, \quad \forall k \ge 1,
	\end{equation*}
with $f_{\tau,\eps} \left[ \lambda_{k}^{(j)} \right] = e^{- ( \tau - \eps) \lambda_{k}^{(j)} } $, for any $ k \ge 1$ and $j:1 \le j \le g_k$. If we introduce the function $f_{\tau,\eps} : x \mapsto e^{- ( \tau - \eps) x }$, then, from the mean value theorem for divided differences, for any $j:1 \le j \le g_k$ there exists $\xi_j \in \left( \lambda_{k}^{(1)}, \lambda_{k}^{(j)}\right) $ such that
	$$
\left| f_{\tau,\eps} \left[ \lambda_{k}^{(1)} , \dots , \lambda_{k}^{(j)}  \right] \right| 
= \left| \frac{d^{j - 1}f_{\tau,\eps}}{d x^{j - 1} }  (\xi_j)  \right| \le T^{j - 1} e^{- (\tau - \eps) \lambda_k^{(1)}}. 
	$$

Taking into account~\eqref{Gk2}, we deduce
	\begin{equation*}
\calK_{m, k}^{\tau, \eps} \le \max \{ 1, T^{p - 1}\} e^{- \mu_m \tau} e^{- (\tau - \eps) \lambda_k^{(1)}}, \quad \forall k \ge 1,
	\end{equation*}
and
	\begin{equation}\label{aux2}
\left\| R_{m, k }^{\tau , \eps} \right\|_{L^2((0, T) \times \Omega_1)} \leq \calC_0  \max \{ 1, T^{p - 1}\}  \exp \left( \frac{\calC_0} {T^{ \theta '} } \right) e^{\calC_0 \left[ \lambda_{m,k}^{(1)} \right]^\theta } e^{- (\tau - \eps) \lambda_{m, k}^{(1)} } ,
	\end{equation}
for any $(m, k) \in \N^2 $.

Let us now demonstrate the result. To do so, let us consider $N \ge 1$, $P_N$ given by~\eqref{PN} and $\alpha, b > 0$. 
Using successively the orthonormality of the sequence $\calB_1 = \left\{ \psi_m \right\}_{m \ge 1}$ in $L^2(\Omega_1)$ and the block moment problem~\eqref{fblock} (remember that $ f_{m,k}^{(j)} = e^{ - \lambda_{m,k}^{(j)} \tau }$) we get
\begin{align*}
\int_0^T \int_{ \Omega_1} R_{m,k}^{\tau, \eps} (t,x') P_N(t,x')   \, dx'\, dt
&= \sum_{k' = 1}^N \sum_{\mu_{m'}^\vartheta \leq N} \int_{\Omega_1} \psi_{m'}(x') \psi_m(x') \, dx' \int_0^T r_{m,k}^{\tau, \eps}(t) g_{m',k'}^{(N)}(t) \, dt
\\
&= \sum_{k'=1}^N \sum_{j=1}^{g_{k'}} a_{m,k}^{(j,N)} \int_0^T r_{m,k}^{\tau, \eps}(t) e^{-\lambda_{m,k'}^{(j)} t} \, dt \frakC_2 \phi_k^{(j)}
\\
&= \sum_{j=1}^{g_k} a_{m,k}^{j,N} e^{-\lambda_{m,k}^{(j)} \tau} \frakC_2 \phi_k^{(j)}
\\
&= g_{m,k}^{(N)} (\tau). 
\end{align*}
Then, recalling that $R_{m, k }^{\tau , \eps}  \equiv 0$ in  $( 0,\varepsilon ) \times \Omega_1$ we obtain
	\begin{align*}
\Ds \normeU{g_{m,k}^{ ( N ) } ( \tau )} 
& \le 
\left( \int_\eps^T \int_{ \Omega_1} e^{\frac{\alpha \beta}{t^b}}  \left| R_{m,k}^{\tau, \eps} (t,x') \right|^2  \, dx'\, dt \right)^{\frac{1}{2}}
\left( \int_0^T \int_{ \Omega_1} e^{-\frac{\alpha \beta}{t^b}}  \normeU{ P_N(t,x')}^2   \, dx'\, dt \right)^{\frac{1}{2}}
\\
	\noalign{\smallskip}
\Ds & \leq e^{\frac{\alpha \beta }{2 \varepsilon^{ b  } }}\left( \int_0^T\!\! \int_{ \Omega _{1}}e^{-\frac{\alpha \beta }{t^{ b }} } \normeU{P_{N}(t,x')}^{2} \, dx'\, dt \right)^{\frac{1}{2}} \left\| R_{m, k }^{\tau , \eps} \right\|_{L^{2}((0,T) \times \Omega_1)} ,
	\end{align*}
for any $(m,k) \in \N^2: \mu _{m}^\vartheta , k \leq N$. This inequality together with~\eqref{aux2} provide~\eqref{fl2} and the proof of the result. 
\end{proof}


 
\subsection{Proof of Theorem~\ref{tPN}} \label{ss32}
In what follows, we will prove Theorem~\ref{tPN}. To do so, we will work with $\alpha, b > 0 $ (to be determined below). We will first assume that $\alpha $, $b$, $N$ and $T$ satisfy 
	\begin{equation*}
\left( \frac{\alpha }{N} \right)^{{\frac 1{b} }} < T ,
	\end{equation*}
and we will divide the proof into two steps. See Remark~\ref{r4} for the case $T \le (\alpha / N)^{\frac  1{b} } $. 


\paragraph{First step: working on $\left( 0, \left( \alpha / N \right)^{\frac  1{b} } \right)$.}
This step is devoted to the proof of the following result:


\begin{lemma}\label{lrest1}
Under the conditions of Theorem~\ref{tPN}, for any $\alpha, b >0$ and $T > (\alpha / N)^{\frac  1{b} } $, we have: 
	\begin{equation}\label{frest1}
\int_{0}^{\left(\frac \alpha  N \right)^{\frac  1{b} } }  \!\! \int_{ \Omega_1 }  e^{-\frac{\alpha \beta }{t^{ b} }}  \normeU{ P_{N}(t,x')}^{2} \, dx'  \,dt
\leq \int_{0}^{\left(\frac \alpha  N \right)^{\frac  1{b} } } \!\! \int_{ \omega } \normeU{ P_{N}(t,x')}^{2}  \, dx' \,dt ,
	\end{equation}
where $P_{N}$ is given by~\eqref{PN}.
\end{lemma}


\begin{proof}
Let $t \in \left( 0, \left( \alpha / N \right)^{\frac  1{b} } \right)$ be fixed. This implies that $N < \frac{\alpha}{t^b}$.

Using the expression~\eqref{PN}, we can apply to $P_{N}$ the spectral inequality~\eqref{H2} with $b_{m,k}^{(j)} = a_{m,k}^{(j,N)} e^{-\lambda_{m,k}^{(j)} t}$ and $\lambda = N$ to obtain
%
	\begin{equation*}
\Ds \int_{\Omega _{1}} \normeU{ P_{N} (t, x')}^{2} \, dx' 
\leq  e^{\beta N} \int_{\omega } \normeU{ P_{N} (t, x')}^{2} \, dx'  
\leq \exp \left( \frac{\alpha \beta }{t^{ b } }\right) \int_{\omega } \normeU{ P_{N} (t, x')}^{2} \, dx'.
	\end{equation*}
Thus,
	\begin{equation}\label{aux1}
\exp \left( -\frac{\alpha \beta }{t^{ b } } \right)  \int_{\Omega _{1}} \normeU{ P_{N} (t, x') }^{2}\, dx' 
\leq \int_{\omega } \normeU{ P_{N} (t, x')}^{2} \, dx' ,\quad \forall t\in \left( 0,\left(\frac \alpha  N \right)^{\frac  1{b} } \right).
	\end{equation}
Integrating with respect to $t$ in $\left( 0,(\alpha / N)^{\frac  1{b} }\right)$, we deduce inequality~\eqref{frest1}. This ends the proof.
\end{proof}


\begin{remark}\label{r4}
It is interesting to note that if $T \le \left(  \alpha / N \right)^{\frac  1{b} } $ then the proof of Theorem~\ref{tPN} is straightforward. Indeed, if $t \in (0, T)$, in particular, $t < \left(  \alpha / N \right)^{\frac  1{b} } $ and inequality~\eqref{aux1} holds for any $t \in (0,T)$. Again, integrating with respect to $t$ in $( 0, T )$, we deduce the proof of Theorem~\ref{tPN}.
\end{remark}


\paragraph{Second step and main argument: working on $\left( \left(  \alpha / N \right)^{\frac  1{b} } ,T\right) $.}
In this part we will work with $t$ in the interval $\left( \left( \alpha / N \right)^{\frac  1{b} } , T \right)$ and prove an estimate similar to~\eqref{frest1} but in the open set  $\left(  \left(  \alpha / N \right)^{\frac  1{b} }, T \right)  \times \Omega_1$. One has:


\begin{lemma}\label{lrest2}
Under the conditions of Theorem~\ref{tPN}, let us consider the positive constant $b$ given by~\eqref{bcoef}. Then, there exist positive constants $\calC_1$ (only depending on $p$, $\rho$, $\theta$, $\kappa$, $\theta_1$ and $\kappa_1$) and $\alpha_1 $ (only depending on $p$, $\rho$, $\theta$, $\kappa$, $\beta$, $\vartheta$ and $\theta_1$) such that for any $T \in (0, 1]$ and $\alpha \ge \alpha_1$ satisfying $\left(  \alpha / N \right)^{\frac  1{b} } < T$, one has
	\begin{multline}\label{frest2}
\Ds \int_{\left(\frac \alpha  N \right)^{\frac  1b } }^T  \!\! \int_{ \Omega_1 } e^{-\frac{\alpha \beta }{t^b} } \normeU{ P_{N}(t,x') }^{2} \, dx'\,dt \leq 3 \int_{\left(\frac \alpha  N \right)^{\frac  1b } }^T   \!\! \int_{ \omega} \normeU{ P_{N}(t,x') }^{2}  \, dx' \,dt  
\\
\Ds + \calC_1 e^{ \frac{2 \calC_0}{T^{\theta '}}}  \calH(\alpha, \beta, T) \int_{0}^T  \!\! \int_{ \Omega_1 } e^{-\frac{\alpha \beta }{t^b} } \normeU{ P_{N} (t,x')}^{2}  \, dx'\,dt,
	\end{multline}
for any $P_{N}$ given by~\eqref{PN}. In~\eqref{frest2}, $\calC_0 >0$, $\beta >0$ and $\theta '$ are the constants provided, resp., by Corollary~\ref{cblock}, inequality~\eqref{H2} and~\eqref{theta}, and $\calH(\alpha, \beta, T)$ is given by 
	\begin{equation}\label{calH}
\calH(\alpha, \beta, T) = \left( T^{ a \theta_1 } + \alpha^{ \frac{\theta_1} \vartheta} \right) \frac 1{ \calC_\alpha} \exp \left( - \calC_\alpha \right)  \exp \left( - \frac{ \calC_\alpha }{ T^a} \right), \ a = \max \left\{ \frac 1{1 - \vartheta} , \frac 1{1 - \theta} \right\} - 1, 
	\end{equation}
and
	\begin{equation}\label{Calpha}
\calC_\alpha := \frac 12 \left( \alpha^{1/ \vartheta} -4^{ b } \alpha \beta - 2(1 - \theta) \calC_0^{\frac 1 {1 - \theta}} \left( 4 \theta \right)^{ \theta '} \right) \ge \frac 12. 
	\end{equation}
\end{lemma}


\begin{proof}
Fix $t\in \left(0 ,T \right)$. Using orthogonality of the sequence $\calB_1= \{ \psi_m \}_{m \geq 1}$ in $L^2(\Omega_1)$, we can write (see~\eqref{PN} and~\eqref{G_m}):
	\begin{equation} 
	\begin{aligned} \notag 
\int_{\Omega_ 1} \normeU{ P_{N} ( t, x')}^{2}\, dx' 
= &\int_{\Omega_1} \normeU{ \sum_{\mu _{m}^\vartheta \leq \lambda } G_{m}^{ ( N ) } ( t ) \psi_m(x')}^{2} dx' 
\\
+ &\int_{\Omega_1} \normeU{ \sum_{\lambda < \mu _{m}^\vartheta \leq N} G_{m}^{ ( N ) } ( t ) \psi_m(x')}^{2} dx',
	\label{Eq_1}
	\end{aligned}
	\end{equation}
for any $\lambda \in (0, N)$ and $t\in \left( 0 ,T\right) $. The first sum in~\eqref{Eq_1} is estimated by applying the spectral inequality~\eqref{H2}. 
It follows that:
	\begin{equation*}
\int_{\Omega_1} \normeU{ \sum_{\mu _{m}^\vartheta \leq \lambda } G_{m}^{ ( N ) } ( t ) \psi_m(x')}^{2} dx'
\leq e^{\beta \lambda }\underset{:=I}{\underbrace{\int_{\omega } \normeU{ \sum_{\mu _{m}^\vartheta \leq \lambda } G_{m}^{ ( N ) } ( t) \psi _{m}(x' ) }^{2}\, dx' ,}} 
	\end{equation*}
for any $t \in \left( 0 ,T\right)$. Using orthonormality of the sequence $\calB_1= \{ \psi_m \}_{m \geq 1}$ in $L^2(\Omega_1)$, the second term  in~\eqref{Eq_1} writes
\begin{equation*}
\int_{\Omega_1} \normeU{ \sum_{\lambda < \mu _{m}^\vartheta \leq N} G_{m}^{ ( N ) } ( t ) \psi_m(x')}^{2} dx' 
= \sum_{\lambda < \mu _{m}^\vartheta \leq N} \normeU{ G_{m}^{ ( N ) } ( t )}^{2}.
\end{equation*}
But,
	\begin{equation*}
	\begin{array}{l}
\Ds I  =   \int_{\omega } \normeU{ \left( \sum_{\mu _{m}^\vartheta \leq N} - \sum_{\lambda < \mu _{m}^\vartheta \leq N} \right) G_{m}^{ ( N  ) } ( t ) \psi _{m}(x' )}^{2}\, dx' \\ 
	\noalign{\smallskip}
 \Ds \phantom{I} \leq  2 \int_{\omega } \normeU{\sum_{\mu _{m}^\vartheta \leq N} G_{m}^{ ( N ) } ( t) \psi _{m}(x' )}^{2}\, dx' + 2 \int_{\Omega _{1}} \normeU{ \sum_{\lambda <\mu _{m}^\vartheta \leq N} G_{m}^{ ( N ) } (t) \psi _{m}(x' )}^{2}\, dx' \\ 
	\noalign{\smallskip}
 \Ds \phantom{I}\leq  2 \int_{\omega } \normeU{ P_{N} ( t, x')}^{2} \, dx' + 2 \sum_{\lambda <\mu _{m}^\vartheta \leq N} \normeU{
G_{m}^{ ( N ) } (t) }^{2}.
	\end{array}
	\end{equation*}
Thus, one obtains:
	\begin{align*}
\int_{\Omega_1} \normeU{ \sum_{\mu _{m}^\vartheta \leq \lambda } G_{m}^{ ( N ) } ( t ) \psi_m(x')}^{2} dx'
\leq \, &2e^{\beta \lambda }\int_{\omega } \normeU{P_{N} ( t, x')}^{2}\, dx' 
\\
+ &2 e^{\beta \lambda }\sum_{\lambda < \mu _{m}^\vartheta \leq N} \normeU{ G_{m}^{ ( N ) } (t)}^{2},
	\end{align*}
for any $\lambda \in (0, N)$ and $t\in \left( 0 ,T\right) $. Inserting this last inequality in~\eqref{Eq_1} and dividing by $e^{\beta \lambda }$ we get: 
	\begin{equation*}
e^{-\beta \lambda }\int_{\Omega _{1}} \normeU{ P_{N} ( t, x') }^{2}\, dx' \leq 3 \left( \int_{\omega } \normeU{P_{N} ( t, x')}^{2}\, dx' + \sum_{\lambda <\mu _{m}^\vartheta  \leq N} \normeU{ G_{m}^{ ( N ) } ( t)}^{2}\right) ,
	\end{equation*}
for any $\lambda \in (0, N)$ and $t \in \left( 0 ,T\right)$.

From now on we restrict ourselves to the case $t\in \left( \left(  \alpha / N \right)^{\frac  1{b} } ,T\right) $ (recall that $b $ is given by~\eqref{bcoef}). 
If in the last inequality we take $\lambda =\frac{\alpha }{t^{b }}$, then $\lambda \in (0, N)$. Introducing the notation
	\begin{equation}\label{SigmaN}
\Sigma^{ ( N ) } (t) := \sum_{\frac{\alpha }{t^{b}} < \mu_{m}^\vartheta \leq N} \normeU{ G_{m}^{ ( N ) } (t)}^{2}, \qquad  t \in \left( \left( \frac{\alpha }{N} \right)^{ \frac 1{b }} ,T\right), 
	\end{equation}
this becomes
	\begin{equation}\label{Eq_2}
e^{-\frac{\alpha \beta }{t^{b}}} \int_{\Omega _{1}} \normeU{ P_{N} ( t, x')}^{2}\, dx' \leq 3 \left( \int_{\omega } \normeU{ P_{N} ( t, x')}^{2}\, dx' + \Sigma^{ ( N ) } (t) \right) , 
	\end{equation}
for any $t \in \left( \left(  \alpha / N \right)^{\frac  1{b} } ,T\right) $.

Our next objective is to provide an estimate of $\Sigma ^{ ( N ) } ( t ) $ for $t \in \left( \left(  \alpha / N \right)^{\frac  1{b} } ,T\right) $. Let us start with $ G_{m}^{ ( N ) } $. Using~\eqref{G_m} and estimate~\eqref{fl2} in Lemma~\ref{l2} for $\alpha >0$, $b$ given in~\eqref{bcoef}, $\eps = t/4 \in (0, T/4)$ and $\tau = t \in (t/4, T) $, we infer
	$$
\Ds \normeU{ G_{m}^{ ( N ) } (t)} \le \calC_0 \max \{ 1, T^{p - 1}\}  \exp \left( \frac{4^{ b} \alpha \beta }{ 2 t^{ b} } + \frac{\calC_0} {T^{ \theta '} } \right) \mathcal{A}_m \left( \int_{0}^T  \!\! \int_{\Omega _{1}} e^{-\frac{\alpha \beta }{t^{ b}} } \normeU{ P_{N}(t,x')}^{2} \, dx'\,dt \right)^{\frac{1}{2}},
	$$	
for any $t \in \left( \left(  \alpha / N \right)^{\frac  1{b} } ,T\right) $ where
	\begin{equation*}
\mathcal{A}_m = \sum_{k=1}^{N} \exp \left( - \frac 34 \lambda_{m, k}^{(1)} t + \calC_0 \left[ \lambda_{m,k}^{(1)} \right]^\theta \right)
	\end{equation*}
and the constants $\calC_0 >0$, $\beta >0$ and $\theta '$ are provided, resp., by Corollary~\ref{cblock}, inequality~\eqref{H2} and~\eqref{theta}. 
If we now use the inequality
	$$
\calC_0 \left[ \lambda_{m,k}^{(1)} \right]^\theta \leq \frac{t}{4} \lambda_{m,k}^{(1)}  + (1 - \theta) \calC_0^{\frac 1 {1 - \theta}} \left( 4 \theta \right)^{ \theta '} \frac 1 {t^{ \theta '} }, 
	$$
valid for any $t  > 0$, we deduce 
	\begin{equation*}
\Ds \normeU{ G_{m}^{ ( N ) } ( t )} \le \calC_0 \max \{ 1, T^{p - 1}\}  e^{\frac 12 \calH_\alpha } e^{- \frac 12  \mu _{m}  t  } \calS_2 (t) \left( \int_{0}^T  \!\! \int_{\Omega _{1}} e^{-\frac{\alpha \beta }{t^{ b}} } \normeU{P_{N}(t,x')}^{2}  \, dx'\,dt \right)^{\frac{1}{2}} ,	
	\end{equation*}
for any $\alpha > 0$ and $t \in \left( \left(  \alpha / N \right)^{\frac  1{b} } ,T\right) $, where
	\begin{equation*}
\calH_\alpha := \frac{4^{ b} \alpha \beta}{ t^{b} } + \frac{ 2 (1 - \theta) \calC_0^{\frac 1 {1 - \theta}} \left( 4 \theta \right)^{ \theta '} }{ t^{\theta '} } + \frac{2 \calC_0}{T^{\theta '}} \hbox{ and } \calS_2 (t) := \sum_{k \ge 1} e^{- \frac 12 \lambda_k^{(1)} t  }  . 
	\end{equation*}

We can now estimate $\Sigma^{ ( N ) } ( t ) $ for any $t  \in \left( \left(  \alpha / N \right)^{\frac  1{b} } ,T\right)$ and $N \ge 1$ (see \eqref{SigmaN}):
	\begin{equation*}
\Sigma^{ ( N ) } ( t ) \leq \calC_0^2 \max \{ 1, T^{p - 1}\}^2 e^{ \calH_\alpha }  \calS_1 (t) \calS_2(t)^{2}  \int_{0}^T  \!\! \int_{\Omega _{1}} e^{-\frac{\alpha \beta }{t^{ b}} } \normeU{ P_{N}(t,x')}^{2} \, dx'\,dt,
	\end{equation*}
where
	\begin{equation*}
\calS_1 (t) := \sum _{\frac{ \alpha }{t^{b}} <\mu _{m}^\vartheta  } e^{-\mu _{m}t}
	\end{equation*}

We can bound the series appearing in the previous inequality using Lemma~\ref{l1}. Indeed, the sequence $\left\{ \mu_m \right\}_{m \ge 1}$ satisfies~\eqref{H1}. From inequality~\eqref{sum} applied to $q = \theta_1$, $\gamma = \frac{ \alpha^{1/ \vartheta} }{t^{b / \vartheta } } > 0$ and $\sigma = t >0 $ ($b$ is given in~\eqref{theta}), we can write
	$$
\Ds\calS_1 (t)  = \sum _{\frac{ \alpha^{1/ \vartheta} }{t^{b / \vartheta } } <{\mu _{m}} } e^{-\mu _{m}t} \le \widehat \calC_1 \, \frac{ t^{ a \theta_1 } + \alpha^{ \frac{\theta_1} \vartheta} }{t^{ (a + 1) \theta_1}} \exp \left( - \frac{\alpha^{1/ \vartheta} } {t^{ a }} \right), \quad \forall t  \in \left( \left( \frac{\alpha }{N} \right)^{ \frac 1{b }} ,T\right), 
	$$
where $ \alpha > 0$, 
	\begin{equation*}
a = \frac{b}\vartheta - 1 = \max \left\{ \frac 1{1 - \vartheta} , \frac 1{1 - \theta} \right\} - 1, 
	\end{equation*}
and $\widehat \calC_1 = \widehat \calC_1 (\kappa_1, \theta_1) > 0$ is the constant provided by Lemma~\ref{l1}.

On the other hand, the sequence $\Lambda_2$ satisfies~\eqref{H4}. Thus, again using Lemma~\eqref{l1} and inequality~\eqref{sum}, we get
	\begin{equation*}
\calS_2 (t) = \sum_{ k\geq 1} e^{- \frac 12 \lambda_k^{(1)} t  } \le \sum_{\lambda \in \Lambda_2 } e^{- \frac 12 \lambda t  } \leq \frac{ \widehat \calC_2}{  t^\theta }, \quad \forall t > 0,
	\end{equation*}
where $\widehat \calC_2$ is a positive constant depending on $\kappa$ and $\theta$.

Coming back to the estimate of $\Sigma^{ ( N ) } ( t )$, we deduce
	$$
	\begin{array}{l}
\Ds \Sigma^{ ( N ) } ( t ) 
\le \frac{\calC_1}3 \max \{ 1, T^{p - 1}\}^2 e^{ \calH_\alpha }   \frac{ t^{ a \theta_1 } + \alpha^{ \frac{\theta_1} \vartheta} }{t^{ (a + 1) \theta_1 + 2 \theta} } \exp \left( - \frac{\alpha^{1/ \vartheta} } {t^{a }} \right) \int_{0}^T  \!\! \int_{\Omega _{1}} e^{-\frac{\alpha \beta }{t^{ b}} } \normeU{ P_{N}(t,x')}^{2} \, dx'\,dt 
\\
	\noalign{\smallskip}
\Ds \phantom{ \Sigma^{ ( N ) } ( t )  } = \frac{\calC_1}3 \max \{ 1, T^{p - 1}\}^2 e^{ \frac{2 \calC_0}{T^{\theta '}}} h (t) \int_{0}^T  \!\! \int_{\Omega _{1}} e^{-\frac{\alpha \beta }{t^{ b}} } \normeU{ P_{N}(t,x')}^{2} \, dx'\,dt ,
	\end{array}
	$$
for any $\alpha > 0$, $t \in \left( \left(  \alpha / N \right)^{\frac  1{b} } ,T\right) $ and $x \in (0, \pi) $, where $a$ and $b$ are given in~\eqref{calH} and~\eqref{bcoef}, $\calC_1 $ is a new positive constant, only depending on $p$, $\rho$, $\theta$, $\kappa$, $\theta_1$ and $\kappa_1$, and $h$ is the function
	\begin{equation}\label{h}
h (t ) := \frac{ t^{ a \theta_1 } + \alpha^{ \frac{\theta_1} \vartheta} }{t^{ (a + 1) \theta_1 + 2 \theta} }  \, \exp \left( \frac{ A }{ t^b } + \frac{ B }{ t^{ \theta '} } - \frac{\alpha^{1/ \vartheta} }{ t^{ a } } \right) , \quad  A= 4^{ b } \alpha \beta, \quad B =   2(1 - \theta) \calC_0^{\frac 1 {1 - \theta}} \left( 4 \theta \right)^{ \theta '} ,  
	\end{equation}
with $t \in (0, \infty)$.

Going back to~\eqref{Eq_2}, we get:
	\begin{multline*}
\Ds e^{-\frac{\alpha \beta }{t^b}} \int_{\Omega _{1}} \normeU{ P_{N} ( t, x') }^{2}\, dx' 
\leq 3 \int_{\omega } \normeU{ P_{N} ( t, x')}^{2}\, dx'   
\\
+ \calC_1 \max \{ 1, T^{p - 1}\}^2 e^{ \frac{2 \calC_0}{T^{\theta '}}} h (t) \int_{0}^T  \!\! \int_{\Omega _{1}} e^{-\frac{\alpha \beta }{t^{ b }} } \normeU{ P_{N}(t,x')}^{2} \, dx'\,dt ,
	\end{multline*}
for any $\alpha > 0$ and $t \in \left( \left(  \alpha / N \right)^{\frac  1{\theta'} } ,T\right) $.

At this point, as $1/ \vartheta > 1$, we can choose $\alpha_0 = \alpha_0 (\beta, \theta, \vartheta, \calC_0) > 0$ ($\calC_0 > 0$ is the constant provided by Corollary~\ref{cblock}) sufficiently large such that $ {\alpha^{1/ \vartheta} } - A - B = 2 \calC_\alpha \ge 1 $ for any $\alpha \ge \alpha_0$ (see~\eqref{Calpha}). Using the expression~\eqref{h}, we deduce that $h \in L^1(0, T)$. Integrating the previous inequality with respect to $t$ in $\left( \left(  \alpha / N \right)^{\frac  1 b } ,T\right)$, we obtain
	\begin{multline}\label{f1}
\Ds \int_{\left(\frac \alpha  N \right)^{\frac  1b } }^T  \!\! \int_{\Omega _{1}} e^{-\frac{\alpha \beta }{t^b }} \normeU{ P_{N}(t,x')}^{2} \, dx' \,dt
\leq 3 \int_{\left(\frac \alpha  N \right)^{\frac  1b } }^T  \!\! \int_{\omega } \normeU{ P_{N}(t,x')}^{2} \, dx' \,dt  
\\
+ \calC_1 \max \{ 1, T^{p - 1}\}^2 e^{ \frac{2 \calC_0}{T^{\theta '}}} \calI  \int_{0}^T  \!\! \int_{\Omega _{1}} e^{-\frac{\alpha \beta }{t^b } } \normeU{ P_{N}(t,x')}^{2} \, dx' \,dt,
	\end{multline}
where 
	$$
\calI := \int_0^T h (t)\, dt ,
	$$
and $ h $, $A$ and $B$ are given in~\eqref{h} and $ b > 0$ in~\eqref{bcoef}. Let us estimate $\calI$ and, to this end, let us assume that $T \le 1$. Thanks to the expressions of $a$ and $b$ (see~\eqref{calH} and~\eqref{bcoef}) one has $a \ge b $ and $a \ge \theta '$ (see~\eqref{theta}). Thus, we deduce
	\begin{equation*}
h (t) \le  \frac{ T^{ a \theta_1 } + \alpha^{ \frac{\theta_1} \vartheta} }{t^{ (a + 1) \theta_1 + 2 \theta} } \, \exp \left( - \frac{2 \calC_\alpha }{ t^{ a } } \right) , \quad \forall t \in (0,T), \quad \forall \alpha \ge \alpha_0 ,
	\end{equation*}
($\calC_\alpha \ge 1 $ is given in~\eqref{Calpha}), and, 
	\begin{align*}
\Ds \calI \le &\int_0^T \frac{ T^{ a \theta_1 } + \alpha^{ \frac{\theta_1} \vartheta} }{t^{ (a + 1) \theta_1 + 2 \theta} } \, \exp \left( - \frac{2 \calC_\alpha }{ t^{ a } } \right) \, dt 
\\
\le \, &\calM \left( T^{ a \theta_1 } + \alpha^{ \frac{\theta_1} \vartheta} \right) \int_0^T \frac{1}{t^{a + 1}} \exp \left( - \frac{ \calC_\alpha }{t^a} \right) \, dt  
\\
=  \, &\calM \left( T^{ a \theta_1 } + \alpha^{ \frac{\theta_1} \vartheta} \right) \frac 1{a \calC_\alpha} \exp \left( - \frac{ \calC_\alpha }{T^a} \right) ,
	\end{align*}
where
	$$
\calM = \sup_{t \in (0, 1)} \frac{1}{t^{ (a + 1)( \theta_1 - 1) + 2 \theta} } \exp \left( - \frac{ \calC_\alpha }{ t^a} \right) = \exp \left( - \calC_\alpha \right) , \hbox{ if } \calC_\alpha \ge \frac 1a \left(  (a + 1)( \theta_1 - 1) + 2 \theta \right).
	$$

Summarizing, there exists a new constant $\alpha_1 >0$ (only depending on $p$, $\rho$, $\theta$, $\kappa$, $\beta$, $\vartheta$ and $\theta_1$) such that for any $T \le 1$ and $\alpha \ge \alpha_1$ one has
	\begin{equation*}
\calI \le \left( T^{ a \theta_1 } + \alpha^{ \frac{\theta_1} \vartheta} \right) \frac 1{a \calC_\alpha} \exp \left( - \calC_\alpha \right)  \exp \left( - \frac{ \calC_\alpha }{ T^a} \right). 
	\end{equation*}

Coming back to~\eqref{f1}, we obtain~\eqref{frest2} for a new positive constant $\calC_1 $, only depending on $p$, $\rho$, $\theta$, $\kappa$, $\theta_1$ and $\kappa_1$. This completes the proof of the result.
\end{proof}



\begin{remark}\label{r2}
In the proof of Lemma~\ref{lrest2} we have assumed that $T \in (0, 1]$. In this case, the constant $\alpha_1 > 0 $ provided by this result is independent of $T$. In the general case $T > 1$, it is possible to prove a slightly different version of this result. Indeed, let us consider inequality~\eqref{f1} with $\alpha \ge \alpha_0 $ (recall that if $\alpha \ge \alpha_0$, one has that $\calC_\alpha \ge 1/2$, see~\eqref{Calpha}). In this case, $h \in L^1 (0, T)$ (see~\eqref{h}). The goal is again to obtain an estimate of $\calI$. To do so, we rewrite 
	\begin{equation*}
h(t) \le  \frac{ T^{ a \theta_1 } + \alpha^{ \frac{\theta_1} \vartheta} }{t^{ (a + 1) \theta_1 + 2 \theta} }  \, \exp \left( \frac{ A }{ t^b } + \frac{ B }{ t^{ \theta '} } - \frac{\alpha^{1/ \vartheta} }{ t^{ a } } \right) = \left( T^{ a \theta_1 } + \alpha^{ \frac{\theta_1} \vartheta} \right) h_1 (t) h_2 (t) , \quad \forall t \in (0,T)  , 
	\end{equation*}
with 
	\begin{equation*}
	\begin{array}{l}
\Ds h_1(t ) = \frac{1}{t^{\frac{a + 1}{2}}} \exp \left(- \frac{ \calC_\alpha }{t^a } \right) , \quad  \calC_\alpha = \frac 12 \left( \alpha^{1/ \vartheta} - A - B \right), \\
	\noalign{\smallskip}
\Ds h_2(t ) = \frac{1}{t^\gamma } \, \exp \left(\frac{ \calC_\alpha }{t^a } \right) \, \exp \left( \frac{ A }{ t^b } + \frac{ B }{ t^{ \theta '} } - \frac{\alpha^{1/ \vartheta} }{ t^{ a } } \right), \quad \gamma = (a +1)\left(\theta_1 - \frac 12 \right) + 2 \theta . 
	\end{array}
	\end{equation*}
Thus,
	\begin{equation}\label{I}
	\begin{array}{l}
\Ds \calI = \int_0^T h(t) \, dt \le \left( T^{ a \theta_1 } + \alpha^{ \frac{\theta_1} \vartheta} \right) \left( \int_0^T h_1(t)^2 \, dt \right)^{1/2} \left( \int_0^T h_2(t)^2 \, dt \right)^{1/2}  \\
  	\noalign{\smallskip}
\Ds \phantom{\calI} = \left( T^{ a \theta_1 } + \alpha^{ \frac{\theta_1} \vartheta} \right) \calI_1 F (\alpha)^{1/2} . 
	\end{array}
	\end{equation}
On one hand,
	\begin{equation*}
\calI_1 = \frac 1{\sqrt{2a \calC_\alpha } } \exp \left( - \frac{ \calC_\alpha }{T^a} \right), \quad \forall \alpha \ge \alpha_0. 
	\end{equation*}
On the other hand, $F : \alpha \in [ \alpha_0, \infty) \mapsto F (\alpha) \in \R$ is the function
	\begin{equation*}
F (\alpha) = \int_0^T h_2(t)^2 \, dt ,
	\end{equation*}
which is well defined when $\alpha \ge \alpha_0$. Then, $F$ is differentiable
in any compact set of $[\alpha_0, \infty)$ and
	\begin{equation*}
	\begin{array}{l}
\Ds F' (\alpha )= \int_0^T \frac{1}{t^{ a } } \, h_2(t)^2 \left( 2\cdot 4^b \beta t^{a- b } - 4^b \beta - \frac 1\vartheta \alpha^{ \frac 1 \vartheta - 1} \right) \, dt \\
	\noalign{\smallskip}
\Ds \phantom{ F' (\alpha )} \le \left( 2\cdot 4^b \beta T^{a- b } - 4^b \beta - \frac 1\vartheta \alpha^{ \frac 1 \vartheta - 1} \right) \int_0^T \frac{1}{t^{ a } } \, h_2 (t)^2 \, dt \le 0, \quad \forall \alpha \in \left[ \alpha_1 (T), \infty \right),
	\end{array}
	\end{equation*}
with 
	\begin{equation*}
\widehat \alpha_1 (T) = \max \left\{ \alpha_0 ,   \left[4^b \beta \vartheta \left( 2 T^{a- b } - 1 \right) \right]^{\frac{\vartheta}{1 - \vartheta}} \right\}. 
	\end{equation*}
%
%
Coming back to~\eqref{I}, we deduce therefore
	\begin{equation*}
\calI \le \left( T^{ a \theta_1 } + \alpha^{ \frac{\theta_1} \vartheta} \right)  F (\widetilde \alpha_1 (T) )^{1/2} \frac 1{\sqrt{2a \calC_\alpha } } \exp \left( - \frac{ \calC_\alpha }{T^a} \right) \quad \forall \alpha \ge \widehat \alpha_1 (T),
	\end{equation*}
when $T > 1$. 

From~\eqref{f1}, we infer the existence of two constants $\widehat C_1 >0$, only depending on $p$, $\rho$, $\theta$, $\kappa$, $\beta$, $\vartheta$, $\theta_1$, $\kappa_1$ and $T$, and $\widehat \alpha_1 > 0 $, only depending on $p$, $\rho$, $\theta$, $\kappa$, $\beta$, $\vartheta$ and $T$, such that for any $\alpha \ge \widehat \alpha_1 $ satisfying $\left(  \alpha / N \right)^{\frac  1{b} } < T$, one has

	\begin{multline}\label{frest2p}
\Ds \int_{\left(\frac \alpha  N \right)^{\frac  1b } }^T  \!\! \int_{\Omega_1} e^{-\frac{\alpha \beta }{t^b} } \normeU{ P_{N}(t,x')}^{2} \, dx'\,dt 
\leq 3\int_{\left(\frac \alpha  N \right)^{\frac  1b } }^T   \!\! \int_{ \omega} \normeU{ P_{N}(t,x')}^{2}  \, dx' \,dt   \\
	\noalign{\smallskip}
+ \widehat \calC_1 T^{p-1} e^{ \frac{2 \calC_0}{T^{\theta '}}} \widehat \calH(\alpha, \beta, T) \int_{0}^T  \!\! \int_{\Omega_1} e^{-\frac{\alpha \beta }{t^b} } \normeU{ P_{N}(t,x')}^{2}  \, dx' \,dt,
	\end{multline}
with $P_{N}$ given by~\eqref{PN}. In~\eqref{frest2p}, $\calC_0 >0$, $\beta >0$ and $\theta '$ are the constants provided, resp., by Corollary~\ref{cblock}, inequality~\eqref{H2} and~\eqref{theta}, and $\widehat \calH(\alpha, \beta, T)$ is given by 
	\begin{equation}\label{calHp}
\widehat \calH(\alpha, \beta, T) = \left( T^{ a \theta_1 } + \alpha^{ \frac{\theta_1} \vartheta} \right) \frac 1{ \sqrt {\calC_\alpha } }  \exp \left( - \frac{ \calC_\alpha }{ T^a} \right).
	\end{equation}
and $\calC_\alpha \ge 1/2$ by~\eqref{Calpha}. 
\end{remark}


We are now in a position to prove Theorem~\ref{tPN}:


\begin{proof}[Proof of Theorem~\ref{tPN}]
We will prove this result as a consequence of Lemmas~\ref{lrest1} and~\ref{lrest2}. Let us first take a final time $T \in (0, 1]$ and $\alpha > 0$ satisfying~\eqref{Calpha}. We prove the result for any $N \ge 1$ and $P_N$ given by~\eqref{PN}, with $a_{m,k}^{ (j , N ) } \in \R$, for any $k, m \ge 1:  k, \mu_m^a \le N$ and $1 \le j \le g_k$. 

If $T \le \alpha/N $ then the proof of Theorem~\ref{tPN} can be deduced reasoning as in Remark~\ref{r4}. 

Let us now assume that $ T > \alpha/N $ and consider the constants $\calC_1 , \alpha_1 > 0$ provided by Lemma~\ref{lrest2}. Combining inequalities~\eqref{frest1} and~\eqref{frest2}, for any $\alpha \ge \alpha_1 $, we get:
	\begin{multline}\label{f5}
\Ds \int_{0}^T  \!\! \int_{ \Omega_1 }  e^{-\frac{\alpha \beta }{t^{ b} }} \normeU{ P_{N}(t,x')}^{2}  \, dx' \,dt 
\leq 3 \int_0^T   \!\! \int_{ \omega } \normeU{ P_{N}(t,x')}^{2}  \, dx' \,dt  \\
	\noalign{\smallskip}
\Ds + \calC_1 e^{ \frac{2 \calC_0}{T^{\theta '}}}  \calH(\alpha, \beta, T) \int_{0}^T  \!\! \int_{ \Omega_1 } e^{-\frac{\alpha \beta }{t^b} }\normeU{ P_{N}(t,x')}^{2}  \, dx'\,dt,
	\end{multline}
where $\calC_0 >0$, $\beta >0$ $\theta '$, $b$ and $\calC_\alpha$ are the constants provided, resp., by Corollary~\ref{cblock}, inequality~\eqref{H2},~\eqref{theta}, \eqref{bcoef} and~\eqref{calH}, and $\calH(\alpha, \beta, T)$ is given by~\eqref{calH}. 

Let us remember that $\calC_\alpha \ge 1/2$ for any $\alpha \ge \alpha_1 $. On the other hand, there exists a new positive constant $\alpha_2 \ge \alpha_1$ (only depending on $p$, $\rho$, $\theta$, $\kappa$, $\beta$, $\vartheta$ and $\theta_1$) such that $\calC_\alpha \ge 2 \calC_0 + 1$ for any $\alpha \ge \alpha_2$. In particular, if we take into account that $0 < T \le 1$ and $\theta' \le a $ (see~\eqref{theta} and~\eqref{calH}) and take $\alpha = \alpha_2$, from the expression of $\calH(\alpha, \beta, T)$, we deduce
	\begin{multline*}
\Ds \int_{0}^T  \!\! \int_{ \Omega_1 }  e^{-\frac{\alpha \beta }{t^{ b} }} \normeU{ P_{N}(t,x')}^{2}  \, dx'  \,dt
\leq 3 \int_0^T   \!\! \int_{ \omega} \normeU{ P_{N}(t,x')}^{2}  \, dx' \,dt
\\
\Ds +   \calH_0( T) \int_{0}^T  \!\! \int_{ \Omega_1 } e^{-\frac{\alpha \beta }{t^b} } \normeU{ P_{N}(t,x')}^{2}  \, dx'\,dt,
	\end{multline*}
where 
	$$
\calH_0( T) = \calC_2 \left( T^{ a \theta_1 } + \alpha^{ \frac{\theta_1} \vartheta} \right)  \exp \left( - \frac{ 1 }{ T^a} \right), 
	$$
and $\calC_2 > 0 $ is a new constant only depending on $p$, $\rho$, $\theta$, $\kappa$, $\beta$, $\vartheta$, $\theta_1$ and $\kappa_1$.

From the previous expression, it is clear that there exists $ \tau_ 0 \in (0, 1]$, only depending on $p$, $\rho$, $\theta$, $\kappa$, $\beta$, $\vartheta$, $\theta_1$ and $\kappa_1$, such that 
	$$
\calH_0 (T) \le \frac 12, \quad \forall T \in (0, \tau_0]. 
	$$
This implies that the following inequality 
	\begin{multline}\label{f6}
\Ds  \int_{0}^T  \!\! \int_{ \Omega_1 }  e^{-\frac{\alpha \beta }{t^{ b} }} \normeU{ P_{N}(t,x')}^{2} \, dx' \,dt
\leq  3 \int_0^T   \!\! \int_{ \omega} \normeU{ P_{N}(t,x')}^{2}  \, dx' \,dt 
\\
+  \frac{1}{2} \int_{0}^T  \!\! \int_{ \Omega_1}  e^{-\frac{\alpha \beta }{t^{ b} }} \normeU{P_{N}(t,x')}^{2}  \, dx' \,dt, 
	\end{multline}
holds for any $T \in (0, \tau_0]$ with $\alpha = \alpha_2$. This concludes the proof of Theorem~\ref{tPN}. 
\end{proof}



\begin{remark}\label{r5}
Taking into account Remark~\ref{r2}, it is possible to prove a version of Theorem~\ref{tPN} valid for any $ T > 0$. Indeed, let us take $T >0$, $N \ge 1$ and $P_N$ given by~\eqref{PN}, with $a_{m,k}^{ (j , N ) } \in \R$, for any $k, m \ge 1:  k, \mu_m^a \le N$ and $1 \le j \le g_k$.  As in the proof of Theorem~\ref{tPN}, if $T \le \alpha/N $ then we can reason as in Remark~\ref{r4}. If $T > \alpha/N $, we take $ \alpha \ge \widehat \alpha_1$ and we combine inequalities~\eqref{frest1} and~\eqref{frest2p}. We obtain an inequality as~\eqref{f5} with the constant
	\begin{equation*}
\widehat \calC_1 T^{p-1} e^{ \frac{2 \calC_0}{T^{\theta '}}} \widehat \calH(\alpha, \beta, T) ,
	\end{equation*}
instead of $\calC_1 e^{ \frac{2 \calC_0}{T^{\theta '}}}  \calH(\alpha, \beta, T)$ where $\widehat \calH(\alpha, \beta, T)$ is given in~\eqref{calHp}. We can argue as in the proof of Theorem~\ref{tPN} and obtain a new positive constant $ \widehat \alpha_2 \ge  \widehat \alpha_1$ (only depending on $p$, $\rho$, $\theta$, $\kappa$, $\beta$, $\vartheta$, $\theta_1$, $\kappa_1$ and $T$) such that
	\begin{equation*}
\widehat \calC_1 T^{p-1} e^{ \frac{2 \calC_0}{T^{\theta '}}} \widehat \calH(\alpha, \beta, T) \le \frac 12 , \quad \forall \alpha \ge  \widehat \alpha_2 .
	\end{equation*}
In particular, taking $\alpha = \widehat \alpha_2 $ we deduce the following version of Theorem~\ref{tPN} which is valid for any $T > 0$:

\begin{theorem}\label{tPNp}
Let us assume that $(\Lambda_1, \calB_1, \Lambda_2, \calB_2, \calU_2 , \frakC_2)$ satisfies Assumption~\ref{A1}. Then, there exists $\alpha > 0$ (only depending on $p$, $\rho$, $\theta$, $\kappa$, $\beta$, $\vartheta$, $\theta_1$, $\kappa_1$ and $T$)  such that
	\begin{equation*}
\Ds  \int_{0}^T  \!\! \int_{ \Omega_1 }  e^{-\frac{\alpha \beta }{t^{ b} }} \normeU{ P_{N}(t,x')}^{2}  \, dx' \,dt
\leq  6 \int_0^T   \!\! \int_{ \omega} \normeU{ P_{N}(t,x')}^{2}  \, dx'\,dt, 
	\end{equation*}
for any $ T > 0 $, any $N \ge 1$ and any $P_N $, where $P_N$ is given by~\eqref{PN}, $\beta $ is the constant in~\eqref{H2} and $b$ is given by~\eqref{bcoef}.
\end{theorem}


\end{remark}



\section{Proof of the main result: Existence of biorthogonal families in $L^2 ((0,T) \times \omega; \calU_2)$} \label{s4}
We devote this section to prove Theorem~\ref{main}. 
Its proof is done in two steps. First, we design a suitable biorthogonal family to $\left\{ F_{m,k}^{(j)}\right\}_{\substack{ m, k \ge 1\\  1 \le j \le g_k } }$ (see~\eqref{Fmk}) in $L^2 ((0,T) \times \Omega_1; \calU_2)$ using in a fundamental way the biorthogonal family coming from Appendix~\ref{Annexe1}. This is done in Section~\ref{ss41} (see Proposition~\ref{p2}). 

Then, in Section~\ref{ss42} we deduce the existence of a biorthogonal family to $\left\{ F_{m,k}^{(j)}\right\}_{\substack{ m, k \ge 1\\  1 \le j \le g_k } }$ in $L^2 ((0,T) \times \omega; \calU_2)$ and end the proof of Theorem~\ref{main} thanks to the following consequence of Theorem~\ref{tmainrest}. 

\begin{proposition}\label{Prop:restriction_biortho}
Assume that $(\Lambda_1, \calB_1, \Lambda_2, \calB_2, \calU_2 , \frakC_2)$ satisfies Assumption~\ref{A1} and consider $ \tau_0 , \alpha > 0$, the constants provided by Theorem~\ref{tPN}, $b$ given by~\eqref{bcoef} and $\eta_\alpha$, the function defined in~\eqref{eta} with $\beta > 0$ given in~\eqref{H2}. Assume that $T \in (0, \tau_0]$. 

For any $q \in L^2((0,T) \times \Omega_1; \calU_2)$ such that 
\begin{equation}\label{PropRestriction_Hypo}
\int_0^T \int_{\Omega_1} e^{\frac{\eta_\alpha(x')}{t^b}} \normeU{q(t,x')}^2 dx' dt < +\infty,
\end{equation}
there exists $Q \in E^\omega \subset L^2((0,T) \times \omega; \calU_2)$ (see~\eqref{eta_alpha}) such that, for any $(m,k) \in \N^2$ and $j : 1 \leq j \leq g_k$,
\begin{equation} \label{PropRestriction_conservation_ps}
\int_0^T \int_\omega \psU{F_{m,k}^{(j)}(t,x')}{Q(t,x')} dx' dt 
= \int_0^T \int_{\Omega_1} \psU{F_{m,k}^{(j)}(t,x')}{q(t,x')} dx' dt 
\end{equation}
and
\begin{equation} \label{PropRestriction_estimee}
\| Q \|_{ L^2((0,T)\times\omega; \calU_2)}^2 \leq 7  \int_0^T \int_{\Omega_1} e^{\frac{\eta_\alpha(x')}{t^b}} \normeU{q(t,x')}^2 dx' dt.
\end{equation}
\end{proposition}

\begin{proof}
We proved, in Theorem~\ref{tmainrest}, that the restriction operator $\calR_\omega \in \calL \left( E_{\eta _{\alpha }}  , E^\omega \right)$ is an isomorphism satisfying~\eqref{dir_indir_ineq} or, equivalently, 
	\begin{equation}\label{f0}
1 \le \left\| \calR_\omega^{-1} \right\|_{\calL \left( E^\omega  , E_{\eta _{\alpha }}   \right) }^2 \le 7. 
	\end{equation}
Recall that the Hilbert spaces $E_{\eta _{\alpha }}$ and $E^\omega $ are given in~\eqref{eta_alpha} and consider
	$$
\calP_{\eta_\alpha} : L_{\eta _{\alpha }}^{2} ((0, T) \times \Omega_1 ; \calU_2) \longrightarrow E_{\eta _{\alpha }}   ,
	$$
the orthogonal projection operator.
 
Let $q \in L^2((0,T) \times \Omega_1; \calU_2)$ and define 
\[
\widetilde q : (t,x') \in (0,T) \times \Omega_1 \mapsto e^{\frac{\eta_\alpha ( x') }{t^b} } q(t,x').
\]
From~\eqref{PropRestriction_Hypo}, we have that $\widetilde q \in L_{\eta _{\alpha }}^{2} ((0, T) \times \Omega_1 ; \calU_2)$.
Finally, we set
	\begin{equation*}
\Psi := \calP_{\eta _{\alpha }} \left( \widetilde q  \right) \quad \mathrm{and} \quad Q := \left( \calR_\omega^{-1}\right)^* \left( \Psi \right) \in E^\omega . 
	\end{equation*}
Then,~\eqref{f0} implies
\begin{align*}
\| Q \|_{ L^2((0,T)\times\omega; \calU_2)}^2
= \left\| \left( \calR_\omega^{-1}\right)^* \calP_{\eta _{\alpha }} \left( \widetilde q  \right) \right\|_{ L^2((0,T)\times\omega; \calU_2)}^2
&\leq 7 \left\| \calP_{\eta _{\alpha }} \left( \widetilde q  \right) \right\|_{ L^2_{\eta_\alpha}((0,T)\times\Omega_1; \calU_2)}^2
\\
&\leq 7 \left\| \widetilde q \right\|_{ L^2_{\eta_\alpha}((0,T)\times\Omega_1; \calU_2)}^2
\end{align*}
which proves~\eqref{PropRestriction_estimee}. We now prove~\eqref{PropRestriction_conservation_ps}. For any $(m,k) \in \N^2$ and $j : 1 \leq j \leq g_k$,
\begin{align*}
\int_0^T \! \! \int_{\Omega_1} \psU{F_{m,k}^{(j)}(t,x')}{q(t,x')} dx' dt 
&= \int_0^T \! \! \int_{\Omega_1} e^{-\frac{\eta_\alpha ( x') }{t^b} } \psU{F_{m,k}^{(j)}(t,x')}{e^{\frac{\eta_\alpha ( x') }{t^b} } q(t,x')} dx' dt 
\\
&= \left( F_{m, k}^{(j)} , \widetilde q  \right)_{L^2_{\eta_{\alpha }}((0,T) \times \Omega_1 ; \calU_2)}   
\\
&= \left( F_{m, k}^{(j)} , \calP_{\eta _{\alpha }} \widetilde q\right)_{L^2_{\eta_{\alpha }}((0,T) \times \Omega_1 ; \calU_2)}  
\\
&= \left( \calR_\omega^{-1} \calR_\omega \left( F_{m, k}^{(j)} \right) , \Psi \right)_{L^2_{\eta_{\alpha }}((0,T) \times \Omega_1 ; \calU_2)} 
\\
&=  \left(  \calR_\omega \left( F_{m, k}^{(j)} \right) , \left( \calR_\omega^{-1}\right)^* \left( \Psi \right) \right)_{L^2 ((0, T) \times \omega; \calU_2)}  
\\
&=\int_0^T \! \! \int_\omega \psU{F_{m,k}^{(j)}(t,x')}{Q(t,x')} dx' dt 
\end{align*}
which ends the proof of Proposition~\ref{Prop:restriction_biortho}. 
\end{proof}


\subsection{Existence of biorthogonal families in $L^2 ((0,T) \times \Omega_1; \calU_2)$}\label{ss41}
Let us consider the sequence $\left\{ F_{m,k}^{(j)}\right\}_{\substack{ m, k \ge 1\\  1 \le j \le g_k } }$ given by~\eqref{Fmk}. Our first objective will be to prove the existence of a biorthogonal family to it in $L^2 ((0,T) \times \Omega_1; \calU_2)$ and give a bound on its norm. 
More precisely, to deal with the weighted norm appearing in the assumption~\eqref{PropRestriction_Hypo}, we prove the following stronger result.

\begin{proposition}\label{p2}
Let us consider $(\Lambda_1, \calB_1, \Lambda_2, \calB_2, \calU_2 , \frakC_2)$ satisfying Assumption~\ref{A1}. Then, there exists a positive constant $ \calC_0 $ such that for any $T >0$ and $\eps \in (0, T/4)$, the family $\left\{ F_{m,k}^{(j)} \right\}_{\substack{ m, k \ge 1\\  1 \le j \le g_k } }$ admits a biorthogonal family  $ \left\{ q^{\eps, (j)}_{m,k} \right\}_{\substack{ m, k \ge 1\\  1 \le j \le g_k } }$ in $L^2 ((0,T) \times \Omega_1; \calU_2)$ satisfying
	\begin{equation}\label{2biort}
	\left\{
	\begin{array}{l}
\Ds q^{\eps, (j)}_{m,k} (t, \cdot) \equiv 0, \quad \forall t \in ( 0,\varepsilon ), \\
	\noalign{\smallskip}
\Ds \left\| q^{\eps, (j)}_{m,k}  \right\|_{L^2 ((0,T) \times \Omega_1; \calU_2)}^2 \leq \calC_0 \exp \left( \frac{ \calC_0} {T^{ \theta '} } \right) e^{ \calC_0 \left[ \lambda_{m, k}^{(1)} \right]^\theta } e^{2 \eps \lambda_{m, k}^{(j)}  } \left( M_k^{-1} \right)_{j,j} , 
	\end{array}
	\right.
	\end{equation}
for any $(m, k) \in \N^2 $ and $ j: 1 \le j \le g_k$ where $\theta'$ is defined by~\eqref{theta} and $M_k$ is the matrix given by~\eqref{Def:M_k}.
\end{proposition}


\begin{proof}
Let $\widetilde C_1 >0$ be the constant given by Proposition~\ref{Prop:Annexe_biortho}. Let us fix $T >0$, $\eps \in (0, T/4)$ and $m \ge 1$. 
Applying Proposition~\ref{Prop:Annexe_biortho} on the time interval $(0,T-\eps)$ we deduce the existence of a biorthogonal family $\left\{\widetilde q_{\eps , m, k}^{(j)} \right\}_{\substack{ k \ge 1 \\  1 \le j \le g_k } } $ to $\left\{ e_{m, k}^{(j)} \right\}_{\substack{ k \ge 1 \\  1 \le j \le g_k } } $ (see~\eqref{ekm}) in $L^2(0,T-\eps; \calU_2) $ satisfying
	\begin{align*}
\Ds \left\| \widetilde q_{\eps , m, k}^{(j)}  \right\|_{L^2(0,T-\eps; \calU_2)}^2 
&\leq \widetilde C_1 \exp \left( \frac{ \widetilde C_1} {(T - \eps)^{ \theta '} } \right) e^{ \widetilde C_1 \left[ \lambda_{m, k}^{(1)} \right]^\theta } \left( M_k^{-1} \right)_{j,j}
\\
&\leq \calC_0 \exp \left( \frac{ \calC_0 } {T^{ \theta '} } \right) e^{ \calC_0 \left[ \lambda_{m, k}^{(1)} \right]^\theta } \left( M_k^{-1} \right)_{j,j}, 
	\end{align*}
for any $k \geq 1$ and $ j: 1 \le j \le g_k$, where $\calC_0 = (\frac43 )^{\theta '} \widetilde C_1$.

Let us consider the function $\widetilde q_{m, k}^{\eps, (j)} $ given by:
	$$
\widetilde q_{m, k}^{\eps, (j)} (t) =
	\left\{
	\begin{array}{ll}
\Ds 0, & \mbox{if } t \in (0, \eps], \\
	\noalign{\smallskip}
\Ds e^{\eps \lambda_{m, k}^{(j)} } \widetilde q_{\eps , m, k}^{(j)} (t - \eps) , & \mbox{if } t \in (\eps , T).
	\end{array}
	\right.
	$$
Then,
	\begin{align*}
\Ds \delta_{k \ell} \delta_{i j} 
&= \int_0^{T - \eps} \!\! \psU{ e_{m, k}^{(j)} (t)}{ \widetilde q_{\eps , m, \ell }^{(i)} (t)} \, dt
 = \int_\eps^{T} \!\! e^{\eps \lambda_{m, k}^{(j)} } \psU{e_{m, k}^{(j)} (t)}{ \widetilde{q}_{\eps, m, \ell}^{(i)} (t - \eps)} \, dt
 \\
	\noalign{\smallskip}
& = \int_0^{T} \!\! \psU{ e_{m, k}^{(j)} (t)}{ \widetilde q_{m, \ell}^{\eps, (i)} ( t)} \, dt,
	\end{align*}
for any $k, \ell \in \N$ and $i,j : 1 \le j \le g_k, \ 1 \le i \le g_\ell$. 
The previous formula proves that, for any $m \in \N$, $ \left\{ \widetilde q_{m, k}^{\eps, (j)} \right\}_{\substack{ k \ge 1 \\  1 \le j \le g_k } }$ is a biorthogonal family to $\left\{ e_{m, k}^{(j)} \right\}_{\substack{ k \ge 1 \\  1 \le j \le g_k } } $ in $L^2(0,T; \calU_2) $. Moreover, $\widetilde q_{m, k}^{\eps, (j)} (t) \equiv 0$ for any $t \in (0, \eps)$ and
	\begin{align}\notag
\Ds \left\| \widetilde q_{m, k}^{\eps, (j)} \right\|_{L^2(0,T; \calU_2)}^2 
&= e^{2 \eps \lambda_{m, k}^{(j)} } \left\| \widetilde q_{\eps , m, k}^{(j)}  \right\|_{L^2(0,T-\eps; \calU_2)} 
\\
&\leq \calC_0 \exp \left( \frac{\calC_0} {T^{ \theta '} } \right) e^{ \calC_0 \left[ \lambda_{m, k}^{(1)} \right]^\theta } e^{2 \eps \lambda_{m, k}^{(j)} } \left( M_k^{-1} \right)_{j,j}, 
\label{aux}
	\end{align}
for any $(m, k) \in \N^2 $ and $ j: 1 \le j \le g_k$. 

Finally, let us define
	\begin{equation*}
q^{\eps, (j)}_{m,k} (t, x') = \widetilde q_{m, k}^{\eps, (j)} (t) \psi_m (x'), \quad \forall (t, x') \in (0,T) \times \Omega_1,
	\end{equation*}
with $(m, k) \in \N^2 $ and $ j: 1 \le j \le g_k$. 
Using that $\calB_1 := \left\{ \psi_m \right\}_{m \ge 1 }$ is an orthonormal basis  of $L^2 (\Omega_1)$ and the expressions of $e_{m, k}^{(j)} $ and $F_{m,k}^{(j)} $ (see~\eqref{ekm} and~\eqref{Fmk}), one has:
	\begin{align*}
&\Ds \int_0^T\!\! \int_{\Omega_1} \psU{F_{m,k}^{(j)}(t,x')}{  q^{\eps, (i)}_{n , \ell}(t,x')}   \, dx' \,dt 
\\
=& \left( \int_0^T \!\! \psU{e_{m, k}^{(j)} (t)}{ \widetilde q^{\eps, (i)}_{n , \ell} (t)}  \,dt \right) \left( \int_{\Omega_1} \psi_{m} (x') \psi_{n} (x')  \, dx' \right) 
\Ds = \delta_{mn} \delta_{k \ell} \delta_{i j}  ,
	\end{align*}
for any $(m,k), (n , \ell ) \in \N^2$ and $i,j : 1 \le j \le g_k, \ 1 \le i \le g_\ell$. On the other hand, property~\eqref{2biort} is a direct consequence of~\eqref{aux}. This ends the proof of Proposition~\ref{p2}.  
\end{proof}


\subsection{Existence of biorthogonal families in $L^2 ((0,T) \times \omega ; \calU_2)$}\label{ss42}

We are now in position to prove Theorem~\ref{main}. 

\begin{proof}[Proof of Theorem~\ref{main}]
Let us first assume that $T \in (0, \tau_0]$, with $\tau_0 > 0$ the final time provided by Theorem~\ref{tPN}. Let us also consider $ \alpha > 0$, the constant provided by Theorem~\ref{tPN}, and $\eta_\alpha$, the function defined in~\eqref{eta}.

Let us construct a family $\left\{ Q_{n,\ell}^{(i)}\right\}_{\substack{ n, \ell \ge 1\\  1 \le i \le g_\ell } }$ biorthogonal to $\left\{ F_{m,k}^{(j)}\right\}_{\substack{ m, k \ge 1\\  1 \le j \le g_k } }$ (see~\eqref{Fmk}) in $L^{2}( ( 0,T ) \times \omega ; \calU_2) $. To this end, we are going to apply Proposition~\ref{p2} in order to construct an appropriate sequence $ \left\{ q^{\eps, (j)}_{m,k} \right\}_{\substack{ m, k \ge 1\\  1 \le j \le g_k } }$ (depending on a parameter $\varepsilon \in (0, T/4)$) biorthogonal to $\left\{ F_{m,k}^{(j)} \right\}_{\substack{ m, k \ge 1\\  1 \le j \le g_k } }$ in $L^2 ((0,T) \times \Omega_1; \calU_2)$. From this sequence and as a consequence of Proposition~\ref{Prop:restriction_biortho}, we will construct the family $\left\{ Q_{m,k}^{(j)}\right\}_{\substack{ m, k \ge 1\\  1 \le j \le g_k } } \subset E^\omega$, biorthogonal to $\left\{ F_{m,k}^{(j)} \right\}_{\substack{ m, k \ge 1\\  1 \le j \le g_k } }$ in $L^2 ((0,T) \times \omega; \calU_2)$. Observe that this family belongs to $E^\omega \subset L^2 ((0,T) \times \omega; \calU_2)$ (for the definition of $E^\omega$, see~\eqref{eta_alpha}) and, therefore, is unique and optimal in the following sense: if $\left\{ \widetilde Q_{m,k}^{(j)}\right\}_{\substack{ m, k \ge 1\\  1 \le j \le g_k } } $ is biorthogonal to $\left\{ F_{m,k}^{(j)} \right\}_{\substack{ m, k \ge 1\\  1 \le j \le g_k } }$ in $L^2 ((0,T) \times \omega; \calU_2)$, then
	$$
\left\| Q_{m,k}^{(j)} \right\|^2_{L^2 ((0, T) \times \omega; \calU_2) } \le \left\| \widetilde Q_{m,k}^{(j)} \right\|^2_{L^2 ((0, T) \times \omega; \calU_2) } , \quad \forall m, k \ge 1 \hbox{ and } j: 1 \le j \le g_k. 
	$$
Thus, the family  $\left\{ Q_{n,\ell}^{(i)}\right\}_{\substack{ n, \ell \ge 1\\  1 \le i \le g_\ell } } \subset E^\omega$ is independent of $\eps$, but not its estimate which comes from Proposition~\ref{p2}. To optimize this estimate, we will make a convenient choice of the parameter $\eps$. More precisely,  let us fix $(n, \ell ) \in \N^2 $ and $i : 1 \le i \le g_\ell$, We begin by applying Proposition~\ref{p2} for $\eps = \eps_{n\ell}^{(i)} \in (0, T/4)$ given by
	\begin{equation}\label{eps}
\eps_{n\ell}^{(i)} = \left\{
	\begin{array}{ll}
\Ds \frac T8,  & \mathrm{if} \ \Ds \frac T4 \le \left( \lambda_{n, \ell}^{(i)} \right)^{ \frac{ - 1}{1 + b}}, \\
	\noalign{\medskip}
\Ds\left( \lambda_{n, \ell}^{(i)} \right)^{ \frac{ - 1}{1 + b}} , & \mathrm{if} \ \Ds \frac T4 > \left( \lambda_{n, \ell}^{(i)} \right)^{ \frac{ - 1}{1 + b}} , 
	\end{array}
	\right.
	\end{equation}
where $b$ is given in~\eqref{bcoef}.
We deduce then the existence of a biorthogonal family 
	$$ 
\left\{ q^{\eps_{n\ell}^{(i)}, (j)}_{m,k} \right\}_{\substack{ m, k \ge 1\\  1 \le j \le g_k } }
	$$ 
in $L^2 ((0,T) \times \Omega_1; \calU_2)$  to $\left\{ F_{m,k}^{(j)} \right\}_{\substack{ m, k \ge 1\\  1 \le j \le g_k } }$ that satisfies~\eqref{2biort} for $\eps_{n \ell}^{(i)}$. 

Since $q^{\eps_{n\ell}^{(i)}, (i)}_{n,\ell} = 0$ for any $t \in (0, \eps_{n\ell}^{(i)})$ we obtain that $q^{\eps_{n\ell}^{(i)}, (i)}_{n,\ell}$ satisfies the assumption~\eqref{PropRestriction_Hypo}. Thus, applying Proposition~\ref{Prop:restriction_biortho} to $q^{\eps_{n\ell}^{(i)}, (i)}_{n,\ell}$, it comes that there exists $Q_{n, \ell}^{(i)} \in E^\omega \subset L^2((0,T)\times\omega; \calU_2)$ such that, for any $ (m, k) \in \N^2$ and $j : 1 \le j \le g_k$, one has
	\begin{align} \notag
\int_0^T \!\! \int_\omega \psU{F_{m, k}^{(j)} (t, x')}{Q_{n, \ell}^{(i)} (t, x')}  \, dx' \, dt
&= \int_0^T \!\! \int_{\Omega_1} \psU{F_{m, k}^{(j)}(t,x')}{q^{\eps_{n\ell}^{(i)}, (i)}_{n,\ell}(t,x')}  \, dx' \, dt
\\
&=\delta_{mn} \delta_{k\ell } \delta_{ij}.
\label{f2}
	\end{align}
Hence, we conclude that the family $\left\{ Q_{m,k}^{(j)} \right\}_{\substack{ m, k \ge 1\\  1 \le j \le g_k } }  $ is biorthogonal to $\left\{ F_{m,k}^{(j)} \right\}_{\substack{ m, k \ge 1\\  1 \le j \le g_k } }$ in $L^{2} ( ( 0,T ) \times \omega ; \calU_2)$. 
To finish the proof, let us check that $\left\{ Q_{m,k}^{(j)} \right\}_{\substack{ m, k \ge 1\\  1 \le j \le g_k } }  $ satisfies~\eqref{main-esti}. From estimate~\eqref{PropRestriction_estimee}, it comes that
\[
\Ds \left\| Q_{m,k}^{(j)} \right\|^2_{L^2 ((0, T) \times \omega; \calU_2) } 
\le 7 \left\| e^{\frac{\eta_\alpha ( \cdot ) } {t^b} } q^{\eps_{mk}^{(j)} , (j)}_{m,k} \right\|^2_{L^2_{\eta_{\alpha }}((0,T) \times \Omega_1 ; \calU_2)}
\]
where $\eta_\alpha $ and $b$ are given in~\eqref{eta} and~\eqref{bcoef}. 
From~\eqref{2biort} we obtain
	\begin{align*}
&\displaystyle \left\| e^{\frac{\eta_\alpha (\cdot ) }{t^b}} q^{\eps_{mk}^{(j)}, (j)}_{m,k} \right\|^2_{L^2_{\eta_{\alpha }}((0,T) \times \Omega_1 ; \calU_2) }
= \int_{\eps_{m, k}^{(j)}}^{T} \! \int_{\Omega_1 } e^{\frac{\eta _{\alpha } ( x' ) }{t^b} } \normeU{q^{\eps_{mk}^{(j)}, (j)}_{m,k}}^2   \, dx'  \, dt
\\ 
\le &\exp \left(\frac{ \alpha \beta }{ \left(\eps_{mk}^{(j)}\right)^b }\right)   \left \Vert q^{\eps_{mk}^{(j)}, (j)}_{m,k} \right \Vert_{L^2((0,T) \times \Omega_1 ; \calU_2)}^2 \\
\le &\calC_0 \exp \left(\frac{ \alpha \beta }{\left(\eps_{mk}^{(j)}\right)^b }\right)   \exp \left( \frac{ \calC_0} {T^{ \theta '} } \right) e^{ \calC_0 \left[ \lambda_{m,k}^{(1)} \right]^\theta } e^{2 \eps_{mk}^{(j)} \lambda_{m,k}^{(j)}  } \left( M_k^{-1} \right)_{j,j}  ,
	\end{align*}
for $ n, \ell \geq 1$ and $i: 1 \le i \le g_\ell$. 

When $  T/4 \le \left( \lambda_{m, k}^{(j)} \right)^{ \frac{ - 1}{1 + b}} $, the previous inequality provides (see~\eqref{eps})
	\begin{equation*}
\Ds \left\| Q_{m,k}^{(j)} \right\|^2_{L^2 ((0, T) \times \omega; \calU_2) } \le  7 \calC_0 \exp \left(\frac{ 8^b \alpha \beta }{T^b }\right) \exp \left( \frac{ \calC_0} {T^{ \theta '} } \right) e^{ \calC_0 \left[ \lambda_{m , k }^{(1)} \right]^\theta } e^{\left[ \lambda_{m , k }^{(j)} \right]^{ \frac{ b }{1 + b}}  } \left( M_k^{-1} \right)_{j,j} , 
	\end{equation*}
for any $ (m, k) \in \N^2$ and $j : 1 \le j \le g_k$. 
On the other hand, when $  T/4 > \left( \lambda_{m, k}^{(j)} \right)^{ \frac{ - 1}{1 + b}}$, we infer 
	\begin{equation*}
\Ds \left\| Q_{m,k}^{(j)} \right\|^2_{L^2 ((0, T) \times \omega; \calU_2) } \le 7 \calC_0 \exp \left( \alpha \beta  \left[ \lambda_{m, k}^{(j)} \right]^{ \frac{ b}{1 + b}}\right) \exp \left( \frac{ \calC_0} {T^{ \theta '} } \right) e^{ \calC_0 \left[ \lambda_{m , k }^{(1)} \right]^\theta } e^{2 \left[ \lambda_{m , k }^{(j)} \right]^{ \frac{ b}{1 + b}}  } \left( M_k^{-1} \right)_{j,j},
	\end{equation*}
for any $ m, k \ge 1$ and $j: 1 \le j \le g_k $. Taking into account~\eqref{Gk2}, from the two previous inequalities we deduce inequality~\eqref{main-esti} when $T \in (0, \tau_0]$.

\medskip

The case $ T > \tau_0$ can be easily deduced reasoning as follows: We consider a family $\left\{ Q_{m,k}^{(j)} \right\}_{\substack{ m, k \ge 1\\  1 \le j \le g_k } }$ biorthogonal to $\left\{ F_{m,k}^{(j)} \right\}_{\substack{ m, k \ge 1\\  1 \le j \le g_k } }$ in $L^2 ((0, \tau_0) \times \omega; \calU_2) $ satisfying~\eqref{main-esti} for a positive constant $\calC$. It is clear that the extension by zeros of $Q_{m,k}^{(j)} $:
	\begin{equation*}
\widetilde Q_{m, k}^{(j)} (t, \cdot ) = 
	\left\{
	\begin{array}{cl}
Q_{m, k}^{(j)} (t, \cdot )   & \mathrm{if} \quad t \le \tau_0,     \\
	\noalign{\smallskip}
 0 & \mathrm{if} \quad t > \tau_0,
 	\end{array}
	\right.
	\end{equation*}
is a biorthogonal family to $\left\{ F_{m,k}^{(j)} \right\}_{\substack{ m, k \ge 1\\  1 \le j \le g_k } }$ in $L^2 ((0, T) \times \omega; \calU_2) $ that also satisfies~\eqref{main-esti}. This ends the proof of Theorem~\ref{main}. 
\end{proof}


\section{Application: Null controllability of a coupled parabolic system}
\label{s5}

In this section we use the biorthogonal family designed in Theorem~\ref{main} to study simultaneous null controllability for a system of two parabolic equations. This system is the extension in higher dimension of the system studied in~\cite{Ou:19}-\cite{Ou:19-b}. 

Let $T > 0$ and $\Omega = \Omega_1 \times (0,\pi)\subset \mathbb{R}^{d}$ with $\Omega_1 \subset \R^{d-1}$ ($ d \ge 2$), a bounded domain with boundary $\partial \Omega_1 \in C^1$. 
In Section~\ref{Subsec:Ex_bord}, we address the boundary null controllability problem for the following diagonal system of two coupled parabolic equations
	\begin{equation}\label{boundary_ap}
	\left \{ 
	\begin{array}{ll}
\Ds \partial _{t}y + \calA y = 0  & \mathrm{in}\ Q_{T} , \\ 
	\noalign{\smallskip}
\Ds y = \mathfrak{b} u1_{\omega \times \{ 0\}} & \mathrm{on}\ \Sigma_{T}  , \\ 
	\noalign{\smallskip}
\Ds y_{\mid t=0}=y_0,   & \mathrm{in}\ \Omega ,
	\end{array}
	\right. 
	\end{equation}
where $\omega \subset \Omega_1 $ is an arbitrary non-empty open set of $\mathbb{R}^{d-1}$, $y = \left( y_{1},y_{2}\right) $ is the state and 
	\begin{equation}\label{Lbcoupled}
\mathfrak{b}=\left( 
	\begin{array}{c}
\mathfrak{b}_{1} \\ 
\mathfrak{b}_{2}
	\end{array}
\right) \in \mathbb{R}^{2}, \quad \calA =\left( 
	\begin{array}{cc}
- \Delta  & 0 \\ 
0 & - \Delta + q
	\end{array}
	\right),
	\end{equation}
with $q \in L^{2}(\Omega )$ satisfying 
	\begin{equation}\label{potential}
q ( x' ,x ) = q ( x ), \quad \hbox{ a.e.~in } \Omega .  
	\end{equation}
For a given initial condition $y_0 \in H^{-1} (\Omega; \R^2) $ the question is the possibility of finding a control $u \in L^2((0, T) \times \omega )$ such that $y(T) = 0$. 

Similarly, in Section~\ref{Subsec:Ex_interne}, we address the internal null controllability problem for the following diagonal system of two coupled parabolic equations
	\begin{equation}\label{internal_ap}
	\left \{ 
	\begin{array}{ll}
\Ds \partial _{t}y + \calA y = \mathfrak{b} u 1_{\omega \times (a,b)}  & \mathrm{in}\ Q_{T} , \\ 
	\noalign{\smallskip}
\Ds y = 0 & \mathrm{on}\ \Sigma_{T}  , \\ 
	\noalign{\smallskip}
\Ds y_{\mid t=0}=y_0,   & \mathrm{in}\ \Omega ,
	\end{array}
	\right. 
	\end{equation}
where $0 \leq a < b \leq \pi$ and $y_0 \in L^2(\Omega ; \R^2)$. 

First, in Section~\ref{Subsec:Ex_prelim}, we prove that these systems fit into the framework of Theorem~\ref{main}.


\subsection{Moment problem and spectral assumptions} \label{Subsec:Ex_prelim}

\paragraph{Boundary control moment problem.}

First, let us deal with the boundary control problem~\eqref{boundary_ap}. 

The operator $-\calA$ with domain $D (\calA ) = H^2(\Omega; \R^2 )\cap H^1_0(\Omega; \R^2)$ is self-adjoint and generates a $C^0-$semigroup on $L^2(\Omega; \R^2)$. Thus, given $y_0 \in H^{-1} (\Omega; \R^2 )$ and  $u \in L^2((0, T) \times \omega)$, the initial value problem~\eqref{boundary_ap} has a unique solution $y \in C^0 \left([0,T] ; H^{-1} (\Omega; \R^2 ) \right) $ that satisfies
 \[
\left\langle y(T), z \right\rangle_{H^{-1},H^1_0} - \left\langle y_0 , e^{-T \calA} z \right\rangle_{H^{-1},H^1_0} 
= \int_0^T \int_\omega u(t,x') \left( \begin{pmatrix} \mathfrak{b}_1 \partial_x \\ \mathfrak{b}_2 \partial_x \end{pmatrix} \cdot e^{-(T-t) \calA} z \right)_{\mid x=0} (x') \, dx' dt
\]
for any $z \in H^1_0(\Omega;\R^2)$. 
Thus, $y(T)=0$ if and only if $v := u(T- \cdot)$ satisfies
\begin{equation} \label{Ex_PbMoments_bord_abstrait}
-\left\langle y_0 , e^{-T \calA} z \right\rangle_{H^{-1},H^1_0} 
= \int_0^T \int_\omega v(t,x') \left( \begin{pmatrix} \mathfrak{b}_1 \partial_x \\ \mathfrak{b}_2 \partial_x \end{pmatrix} \cdot e^{-t \calA} z \right)_{\mid x=0} (x') \, dx' dt
\end{equation}
for any $z \in H^1_0(\Omega;\R^2)$.

\medskip
\paragraph{Spectral analysis.}

Let us introduce some notations. First,  $A_{1}: L^{2} ( \Omega _{1} ) \rightarrow L^{2} ( \Omega_1 ) $  is the Dirichlet-Laplace operator on the open set $\Omega_1 \subset \R^{d - 1}$ defined by:
	\begin{equation*}
A_{1}=-\Delta _{1}=-\sum_{k=1}^{d-1}\frac{\partial ^{2}}{\partial x_{k}^{2}} , \quad D ( A_{1} ) = H^{2} ( \Omega_{1} ) \cap H_{0}^{1}( \Omega_{1} )  .  
	\end{equation*}
We denote its spectrum as $\sigma (A_1) = \Lambda_1 := \{ \mu_m \}_{m\geq 1}$ and $\calB_1 := \{ \psi_m \}_{m\geq 1}$ is the associated sequence of normalized eigenfunctions in $L^2(\Omega_1)$. Secondly, we will also consider 
	\begin{equation*}
\Lambda_2 = \Lambda_2^{(1)} \cup \Lambda_2^{(2)}
	\end{equation*}
where $\Lambda_2^{(1)} = \left\{ \nu_k^{( 1 )} \right\}_{k \ge 1}$ and $\Lambda_2^{(2)}= \left\{ \nu_k^{( 2 )}  \right\}_{k \ge 1}$ are, resp., the sequences of eigenvalues of the operators $- \frac {\partial^2}{\partial x^2} $ and $- \frac {\partial^2}{\partial x^2} + q $ in $(0, \pi) $ with homogenous Dirichlet boundary conditions. 
To fit in the framework studied in this article we assume in all what follows that
	\begin{equation}\label{distinct}
\nu_k^{(1)} \not= \nu_\ell^{(2)}, \quad \forall k, \ell \ge 1,
	\end{equation}
We will denote 
	\begin{equation*}
\calB_2 =\left\{ \phi_k^{( 1 )} , \phi_k^{( 2 )}\right\}_{k \ge 1} \subset H^2 (0, \pi) \cap H_0^1 (0, \pi) ,
	\end{equation*}
the corresponding eigenfunctions associated to the previous operators satisfying
	\begin{equation}\label{H6bis}
\left( \phi_k^{(j)} \right)' (0) = 1, \quad \hbox{for }  j=1,2.
	\end{equation}

In fact, 
	\begin{equation}\label{asymp0}
 \nu_k^{( 1 )}  = k^2, \qquad \phi_k^{( 1 )} (x) = \frac 1k  \sin ( k x ), \quad x \in (0, \pi), \quad \forall k \ge 1,  
	\end{equation}
and, from~\cite{K:95}, Theorem~4.11, page 135, one has:
%
	\begin{equation}\label{asymp}
\Ds \nu_k^{(2)} = \nu_k^{( 1 )} + \overline q + \xi_k, \quad \forall k \ge 1, \quad \overline q = \frac 1\pi \int_{0}^{\pi } q (x)\,dx, 
	\end{equation}
where $\left\{ \xi_k \right\}_{k \ge 1} \in \ell^2$.

\begin{remark}\label{not2}
In order to study the controllability properties of system~\eqref{boundary_ap}, we can assume, without loss of generality, that $ \widehat \nu := \inf_{k \ge 1} \nu_k^{(2)} > 0 $. Indeed, if $\widehat \nu \le 0$, we can perform in problem~\eqref{boundary_ap} the change $\widetilde{y} = e^{-c t} y$ with $c = - \widehat \nu + 1 > 0$. This change transforms~\eqref{boundary_ap} into the equivalent null controllability problem for $\widetilde{y}$:
	\begin{equation*}
	\left \{ 
	\begin{array}{ll}
\Ds \partial _{t} \widetilde{y} + \left(\calA + c \right) \widetilde{y} = 0  & \mathrm{in}\ Q_{T} , \\ 
	\noalign{\smallskip}
\Ds \widetilde{y} = be^{-c t}u1_{\omega \times \{ 0\}} & \mathrm{on}\ \Sigma_{T}  , \\ 
	\noalign{\smallskip}
\Ds \widetilde{y}_{\mid t=0}=y_0, \quad \widetilde{y}_{\mid t=T}=0 & \mathrm{in}\ \Omega .
	\end{array}
	\right. 
	\end{equation*}
This is equivalent to adding the constant $c > 0 $ to the sequences $\Lambda^{(1)}_2$ and $\Lambda^{(2)}_2$.
The same remark holds for system~\eqref{internal_ap}. 
\end{remark}

For the operator~\eqref{Lbcoupled}, thanks to~\eqref{potential}, one has that its spectrum is given by
	\begin{equation}\label{spectrum}
\sigma (\calA ) = \left\{ \nu_{m, k}^{(1)} := \mu_m + \nu_k^{( 1 )} ,  \quad \nu_{m, k}^{(2)} := \mu_m + \nu_k^{( 2 )}  : (m, k )\in \N^2 \right\}.
	\end{equation}
The associated eigenfunctions of $\calA $ are defined on $\Omega $ by
	\begin{equation*}
\Ds \Phi _{m,k}^{ ( 1 ) } ( x',x ) := \begin{pmatrix} \psi _{m} ( x' ) \phi_k^{( 1 )} (x) \\ 0 \end{pmatrix}, 
\qquad
\Phi_{m,k}^{ ( 2 ) } ( x',x ) := \begin{pmatrix} 0 \\ \psi _{m} ( x' ) \phi_k^{( 2 )} (x) \end{pmatrix},
	\end{equation*}
for any $m,k \ge 1$. Thus, getting back to~\eqref{Ex_PbMoments_bord_abstrait}, it comes that the solution $y$ of~\eqref{boundary_ap} satisfies $y(T)=0$ if and only if
\begin{equation} \label{Ex_PbMoments_bord}
-e^{-\nu_{m,k}^{(j)} T} \left\langle y_0 , \Phi _{m,k}^{ ( j ) } \right\rangle_{H^{-1},H^1_0} 
= \mathfrak{b}_j  \int_0^T \int_\omega v(t,x') e^{-\nu_{m,k}^{(j)} t} \psi_m(x') \, dx' dt
\end{equation}
for any $m,k \ge 1$ and $j: 1 \le j \le 2 $, where we have used the normalization condition~\eqref{H6bis}.

\paragraph{Internal control problem.}

For the internal control problem, we consider the normalization condition
\begin{equation}\label{H6ter}
\int_a^b \left( \phi_k^{(j)} (x) \right)^2 dx = 1
\end{equation}
instead of~\eqref{H6bis}. Then, the solution $y$ of~\eqref{internal_ap} satisfies $y(T)=0$ if and only if
\begin{equation} \label{Ex_PbMoments_interne}
-e^{-\nu_{m,k}^{(j)} T} \left\langle y_0 , \Phi _{m,k}^{ ( j ) } \right\rangle_{L^2} 
= \mathfrak{b}_j \int_0^T \int_\omega \int_a^b v(t,x',x) e^{-\nu_{m,k}^{(j)} t} \psi_m(x') \phi_k^{(j)} (x) \, dx dx' dt
\end{equation}
for any $m,k \ge 1$ and $j: 1 \le j \le 2 $.


\subsection{A boundary controllability problem for a coupled parabolic system} \label{Subsec:Ex_bord}

In this section we analyze the boundary null controllability problem associated with~\eqref{boundary_ap} with a particular focus on the minimal final time $T$ needed to achieve such property.

The one-dimensional version of this question ($d = 1$) has been analyzed in~\cite{Ou:19}-\cite{Ou:19-b}. Let us describe the results obtained in this work. First, in order to solve the null controllability problem associated with~\eqref{boundary_ap} ($d = 1$ and $\Omega = (0, \pi) $), it is necessary to impose the conditions~\eqref{distinct} and
	\begin{equation}\label{fb}
\mathfrak{b}_1 \mathfrak{b}_2 \not= 0 .
	\end{equation}
It is clear that conditions~\eqref{distinct} and~\eqref{fb} are also necessary to solve~\eqref{Ex_PbMoments_bord} in the general case $d > 1$. In fact,~\eqref{distinct} and~\eqref{fb} are equivalent conditions to the corresponding approximate controllability property for the parabolic system associated to the null controllability problem~\eqref{boundary_ap}. Secondly, in~\cite{Ou:19} the author analyzes problem~\eqref{boundary_ap} when $d = 1$ and~\eqref{distinct} and~\eqref{fb} holds. Setting 
	\begin{equation}\label{min_time_boun}
T_0( q ):=\limsup_{k \to +\infty} \frac{-\log \left| \nu_{k}^{(2)} - \nu_{k}^{(1)} \right| }{\nu_{k}^{(1)}} \in \left[ 0,\infty \right] ,  
	\end{equation}
the author proves, using the moment method, that:
\begin{itemize}
\item if $T>T_0 (q)$, for any $y_0 \in H^{-1} (\Omega; \R )$ there exists $u\in L^{2} ( 0,T )$ such that the solution $y$ of problem~\eqref{boundary_ap} (with $d = 1$) satisfies $y(T)=0$. 
\item if $T_0 (q) > 0$ and $0<T<T_0 (q)$, then there exists $y_0 \in H^{-1} (\Omega; \R )$ such that for any $u\in L^{2} ( 0,T )$ the solution $y$ of problem~\eqref{boundary_ap} (with $d = 1$) satisfies $y(T)\neq 0$.
\end{itemize}
Notice that when $\overline q \not= 0$, one has $T_0 (q) = 0$ and, under assumptions~\eqref{distinct} and~\eqref{fb}, null controllability for problem~\eqref{boundary_ap} (with $d = 1$) holds for any time $T>0$. 

However, when $q \in L^2(0, \pi)$ satisfies $\overline q = 0$, the elements of the sequence $\left\{ \nu_k^{( 1 )}, \nu_k^{( 2 )} \right\}_{k \ge 1}$ condense:
	$$
\left| \nu_{k}^{(2)} - \nu_{k}^{(1)} \right| = | \xi_ k | \to 0 
	$$
and it can happen that $T_0(q) > 0$. More precisely, for any $\tau_0 \in [0, +\infty]$, there exists $q \in L^2(0, \pi)$ such that $T_0(q)=\tau_0$.

This result has been extended in~\cite{Ou:19-b} to the case $d>1$ in two particular settings
\begin{itemize}
\item[$\bullet$] first, where $\omega=\Omega_1$ in~\cite[Section 4.1]{Ou:19-b}. In this case, the problem remains roughly one-dimensional and the proof uses the study in the one-dimensional case and the fact that $\{ \psi_m \}_{m \ge 1}$ is a Hilbert basis of $L^2(\Omega_1)$.
\item[$\bullet$] then, the case where $\omega\subset \Omega_1$ in~\cite[Section 4.2]{Ou:19-b} but with extra assumptions ensuring that $T_0(q)=0$. In this case, the proof uses in a fundamental way the cost of null controllability in the one-dimensional setting to develop a Lebeau-Robbiano's iteration scheme as in~\cite{BBGBO:2014}.
\end{itemize}

Our objective is now to generalize this controllability result to the setting described at the beginning of Section~\ref{s5}. In this sense, one has:


\begin{theorem}\label{boundary_result} 
Let $\Omega$ and $\omega$ be as defined at the beginning of Section~\ref{s5}. Let $\mathfrak{b}$ and $\calA$ defined by~\eqref{Lbcoupled} with $q \in L^2 (\Omega )$ satisfying~\eqref{potential} and assume conditions~\eqref{distinct} and~\eqref{fb}. Consider $T_0(q) \in [0, \infty]$ given by~\eqref{min_time_boun}. Let $T>0$. Then, 
\begin{enumerate}
\item if $T>T_0 (q)$, for any $y_0 \in H^{-1} (\Omega; \R^2 )$ there exists $u\in L^{2} ( (0,T) \times \omega )$ such that the solution $y$ of problem~\eqref{boundary_ap} satisfies $y(T)=0$. 
\item if $T<T_0 (q)$, then there exists $y_0 \in H^{-1} (\Omega; \R^2 )$ such that for any $u\in L^{2} ( (0,T) \times \omega )$ the solution $y$ of problem~\eqref{boundary_ap} satisfies $y(T)\neq 0$.
\end{enumerate}
\end{theorem}


The proof of Theorem~\ref{boundary_result}, item~$\mathit{2.}$ is a direct consequence of the results obtained in~\cite[Section 4.1]{Ou:19-b}: if problem~\eqref{boundary_ap} is not controllable with $\omega=\Omega_1$ it cannot be controllable in the more restrictive setting $\omega \subset \Omega_1$. 

Theorem~\ref{boundary_result}, item~$\mathit{1.}$ will be proved applying the moment method. It strongly relies on the biorthogonal family designed in Theorem~\ref{main} applied to the framework described in Section~\ref{Subsec:Ex_prelim}. 

In order to prove it, let us first establish the existence of a biorthogonal family in $L^2 ( ( 0,T ) \times \omega )$ of the sequence $ \calF := \left\{ \widetilde F_{m,k}^{(j)}\right\}_{\substack{ m, k \ge 1\\  j= 1, 2 } }$ given by
	\begin{equation}\label{Fmkbis}
\widetilde F_{m,k}^{(j)} (t, x') := e^{-\nu_{m,k}^{ ( j ) }t} \psi_{m} ( x') , \quad ( t,x' ) \in (0,T) \times \Omega_1 , 
	\end{equation}
with $m, k \ge 1$ and  $j= 1, 2$ ($\nu_{m,k}^{ ( j ) }$ is given in~\eqref{spectrum}). 
One has:


\begin{proposition} \label{Prop:Ex_biortho}
Under the previous notations, let us assume that $q \in L^2 (\Omega)$ and the sequence $\left\{ \nu_k^{(1)}, \nu_k^{(2)} \right\}_{  k \ge 1 }$ satisfies~\eqref{potential} and~\eqref{distinct}. Then, there exists a constant $\calC >0$, only depending on $q$, such that for any $T >0$, the sequence $\left\{ \widetilde F_{m,k}^{(j)}\right\}_{\substack{ m, k \ge 1\\  1 \le j \le 2 } }$ (see~\eqref{Fmkbis}) admits a biorthogonal family $\left\{ \widetilde Q_{m,k}^{(j)}\right\}_{\substack{ m, k \ge 1\\  1 \le j \le 2 } }$ in $L^{2}( ( 0,T ) \times \omega ) $ that satisfies 
	\begin{equation}\label{esti-bound}
\Ds \left\| \widetilde Q_{m,k}^{(j)} \right\|_{L^2 ((0, T) \times \omega)} \le \calC \exp \left(\frac{ \calC }{T } \right) \exp \left( \calC \sqrt{ \nu_{m, k}^{(j)} } \right)  \frac{1}{\left| \nu_k^{(1)} - \nu_k^{(2)} \right| },
	\end{equation}
for any $ m, k \ge 1$ and $j: 1 \le j \le 2 $.
\end{proposition}


\begin{proof}
Let us consider $\calU_2 = \R$ and $\frakC_2 $ the linear operator
	$$
\frakC_2 : \phi \in  H^2 (0, \pi) \cap H_0^1 (0, \pi) \mapsto \frakC_2 \phi =   \phi' (0) \in \R,
	$$
which satisfies $\frakC_2 \in \calL ( H^2 (0, \pi) \cap H_0^1 (0, \pi), \R)$. Let us also consider $(\Lambda_1, \calB_1, \Lambda_2, \calB_2)$ given in Section~\ref{Subsec:Ex_prelim}. 
To apply Theorem~\ref{main}, let us check that $(\Lambda_1, \calB_1, \Lambda_2, \calB_2, \calU_2 , \frakC_2)$ satisfies Assumption~\ref{A1}.

Using Weyl's law, it is classical that the sequence $\left\{ \mu_m \right\}_{m \ge 1}$ satisfies
	$$
\calN_{\Lambda_1}  ( r )  \leq \kappa_1 r^{\frac{d-1}{2}}, \quad \forall r \in (0, \infty).
	$$
for a positive constant $\kappa_1$, only depending on $\Omega_1 \subset \R^{d-1}$ and $d$ (the definition of $\calN_{\Lambda_1} $ is given in~\eqref{counting}). So,~\eqref{H1} holds for $\kappa_1 $ and $\theta_1 = (d-1)/2$.

Let us now check conditions~\eqref{H3}--\eqref{H5} for sequence $\Lambda_2$. Recall that $\Lambda_2 = \Lambda_2^{(1)} \cup \Lambda_2^{(2)}$ with (see~\eqref{asymp0} and~\eqref{asymp})
	\begin{equation*}
\Lambda_2^{(1)} = \left\{ \nu_k^{( 1 )} \right\}_{k \ge 1 } \quad \hbox{and} \quad \Lambda_2^{(1)} =\left\{ \nu_k^{( 2 )} \right\}_{k \ge 1 }.
	\end{equation*}
Using a general result for Sturm-Liouville operators (see for instance~\cite[Theorem IV.1.3]{B:22} we can deduce that, for $i = 1,2$, $\Lambda_2^{(i)} \in \calL (1, \rho_i, 1/2, \kappa_i) $ for appropriate constants $\rho_i,  \kappa_i >0 $ (recall that the class $\calL$ is defined in~\eqref{Lclass}). In our particular case and using~\eqref{asymp0} and~\eqref{asymp}, it is not difficult to see that
	$$
\Lambda_2^{(1)} \in \calL (1, 3, 1/2, 1) \quad \hbox{and} \quad \Lambda_2^{(2)} \in \calL (1, \rho_2, 1/2, \kappa_2) 
	$$
with $\rho_2,  \kappa_2 >0 $ only depending on $q$. Now, from~\cite[Lemma~V.4.20]{B:22}, we deduce
	$$
\Lambda_2 \in \calL (p, \rho, \theta, \kappa) \quad \hbox{with} \quad p = 2, \quad \rho = \min \left( 3 , \rho_2 \right), \quad \theta = 1/2, \quad \kappa = 2 \left( 1 + \kappa_ 2 \right). 
	$$
Using~\cite[Proposition~2.2]{BM:23} we deduce the existence of a countable family $\left\{ G_k \right\}_{k\geq 1}$ of disjoint subsets of $\Lambda_2$ satisfying~\eqref{Gk} and~\eqref{Gk2}. In the same way, we can also rearrange the elements of $\calB_2$  in such a way that we have~\eqref{B2}. Condition~\eqref{H6} is a direct consequence of the normalization condition~\eqref{H6bis}.

Finally, inequality~\eqref{H2} is a consequence of a result of Jerison-Lebeau (see~\cite{JL:90}): If $\Omega_1 \subset \R^{d-1}$ is a bounded domain with boundary $\partial \Omega_1$ regular enough, and $\omega \subset \Omega _{1}$ is an arbitrary nonempty open subset of $\R^{d-1}$, then, there exists a constant $\beta >0$, only depending on $\omega$ and $\Omega_1$, such that for any sequence $\left\{ b_{m}\right\}_{m\geq 1}\subset \C$ and any $\lambda \in \left( 0 ,\infty \right) $, one has 
	\begin{equation} \label{Inegalite_Spectrale:JL}
\int_{\Omega_1 }\left \vert \sum_{\sqrt{\mu _{m}}\leq \lambda }b_{m}\psi _{m}(x' )\right \vert ^{2}\, dx'  \leq e^{\beta \lambda } \int_{\omega }\left \vert \sum_{\sqrt{\mu _{m}}\leq \lambda }b_{m}\psi _{m}(x' )\right \vert ^{2}\, dx' .
	\end{equation}
As already noticed in Remark~\ref{Rem:Inegalite_Spectrale}, applying the previous inequality to 
\begin{equation*} 
b_m = \sum_{k \leq \lambda} \sum_{j=1}^{g_k} b_{m,k}^{(j)} \frakC_2 \phi_k^{(j)}
\end{equation*}
implies the validity of~\eqref{H2} with $\beta >0$ and $ \vartheta = 1/2$. 
We have therefore proved that $(\Lambda_1, \calB_1, \Lambda_2, \calB_2, \calU_2 , \frakC_2)$ satisfies Assumption~\ref{A1}.

We now apply Theorem~\ref{main} to obtain a biorthogonal family $\left\{ \widetilde Q_{m,k}^{(j)}\right\}_{\substack{ m, k \ge 1\\  1 \le j \le 2 } }$ to $\left\{ \widetilde F_{m,k}^{(j)}\right\}_{\substack{ m, k \ge 1\\  1 \leq j \leq 2 } }$ in $L^{2}( ( 0,T ) \times \omega )$ satisfying~\eqref{main-esti}. 
To prove the estimate~\eqref{esti-bound} it thus only remains to estimate $\left( M_k^{-1} \right)_{j,j}$ where $M_k$ is the matrix given in~\eqref{Def:M_k}. 
For any $k \geq 1$, we consider the two possible cases. 

\begin{itemize}
\item[$\bullet$] If $\sharp G_k = 2$ then there exists $\ell \geq 1$ such that $G_k = \left\{ \nu_\ell^{(1)}, \nu_{\ell}^{(2)} \right\}$. Then, from~\eqref{Def:M_k}, using the normalization condition~\eqref{H6bis} we have 
\[
M_k = \begin{pmatrix} 1 & 1 \\ 1 & 1 + \left( \nu_{\ell}^{(1)} - \nu_{\ell}^{(2)} \right)^2 \end{pmatrix}.
\]
Thus, for any $j \in \{ 1, 2\}$, we have
\[
\left( M_k^{-1} \right)_{j,j} \leq \frac{1+\rho^2}{\left| \nu_{\ell}^{(1)} - \nu_{\ell}^{(2)}\right|^2}. 
\]
\item[$\bullet$] If $\sharp G_k = 1$ then there exists $\ell \geq 1$ and $i \in \{1, 2 \}$ such that $G_k = \left\{ \nu_\ell^{(i)} \right\}$. Then, from~\eqref{Def:M_k}, we have $M_k = 1$. 

To obtain an estimate valid in both cases, notice that, from~\eqref{asymp} it comes that there exists $C>0$ depdending on $q$ such that
\begin{equation} \label{Sturm_distance_vp}
\left| \nu_k^{(2)} - \nu_k^{(1)} \right| \leq C, \qquad \forall k \geq 1.
\end{equation}
Thus,
\[
M_k^{-1} = 1 \leq \frac{C}{\left| \nu_{\ell}^{(1)} - \nu_{\ell}^{(2)}\right|^2}.
\]
\end{itemize}
Gathering both cases proves~\eqref{esti-bound} and ends the proof of Proposition~\ref{Prop:Ex_biortho}. 
\end{proof}


We now have all the ingredients to prove Theorem~\ref{boundary_result}.
\begin{proof}[Proof of Theorem~\ref{boundary_result}]
The proof of item~$\mathit{2.}$ follows from the study done in~\cite[Section 4.1]{Ou:19-b}: if problem~\eqref{Ex_PbMoments_bord} is null controllable with a control $u \in L^2((0,T) \times \omega)$ then it is null controllable with a control $u \in L^2((0,T) \times \Omega_1)$.
The latter property does not hold if $T<T_0(q)$. 

The proof of item~$\mathit{1.}$ relies on the moment method. Let $T>T_0(q)$ and $y_0 \in H^{-1}(\Omega, \R^2)$. Let $\left\{ \widetilde Q_{m,k}^{(j)} \right\}_{\substack{m, k \geq 1 \\ 1 \leq j \leq 2}}$ be the biorthogonal family designed in Proposition~\ref{Prop:Ex_biortho}. Let us consider $v$ given by the formal series
\begin{equation} \label{Ex:Def_v}
v : (t,x') \mapsto \sum_{m \geq 1} \sum_{k \geq 1} \sum_{j=1}^2 
\frac{-1}{\mathfrak{b}_j} \left\langle y_0 , \Phi_{m,k}^{(j)} \right\rangle_{H^{-1}, H^1_0} e^{-\nu_{m,k}^{(j)} T} \widetilde Q_{m,k}^{(j)}(t,x'). 
\end{equation}
There exists $C>0$ such that $\| \psi_m \|_{H^1_0(\Omega_1)} \leq C \sqrt{\mu_m}$ for any $m \geq 1$. Applying classical results for Sturm-Liouville operator (see for instance~\cite[Lemma 2.3]{ABM_2018}), taking into account the normalization condition~\eqref{H6bis}, there exists $C>0$ depending on $q$ such that
\[
\left| \phi_k^{(j)}(x) \right|^2 + \frac{1}{\nu_{k}^{(j)}} \left| \left( \phi_k^{(j)} \right)'(x) \right|^2 \leq \frac{C}{\nu_k^{(j)}}, 
\qquad \forall x \in (0, \pi), \quad \forall k \geq 1, \: \forall j : 1 \leq j \leq 2.
\]
Thus, we deduce that 
\[
\left\| \Phi_{m,k}^{(j)} \right\|_{H^1_0}^2 \leq C \left( 1 + \frac{\mu_m}{\nu_{k}^{(j)}} \right), 
\qquad \forall m, k \geq 1, \: \forall j : 1 \leq j \leq 2.
\]
From~\eqref{esti-bound} we have
\[
\Ds \left\| \widetilde Q_{m,k}^{(j)} \right\|_{L^2 ((0, T) \times \omega)} \le \calC \exp \left(\frac{ \calC }{T } \right) \exp \left( \calC \sqrt{ \nu_{m, k}^{(j)} } \right)  \frac{1}{\left| \nu_k^{(1)} - \nu_k^{(2)} \right| }
\]
for any $m, k \geq 1$ and $j: 1 \leq j \leq 2$. 
Recall that $T_0(q)$ is defined by~\eqref{min_time_boun}. Then, since $T> T_0(q)$, it comes that the series~\eqref{Ex:Def_v} converges in $L^2((0,T) \times \omega)$. 

Using the biorthogonality property we directly obtain that $v$ solves the moment problem~\eqref{Ex_PbMoments_bord} which ends the proof of Theorem~\ref{boundary_result}. 
\end{proof}


\subsection{A distributed controllability problem for a coupled parabolic system} \label{Subsec:Ex_interne}

In this subsection we give the adjustments with respect to Section~\ref{Subsec:Ex_bord} to study the null controllability with a distributed control problem given in~\eqref{internal_ap}. Recall that the associated moment problem is given by~\eqref{Ex_PbMoments_interne} and we considered the normalization condition~\eqref{H6ter}. 

The following one-dimensional version of this question ($d = 1$)
\begin{equation} \label{Ex:Pb_1D_interne}
\left\{
\begin{aligned}
& \partial_t y_1 - \partial_{xx} y_1 = \mathfrak{b}_1 \mathbf{1}_{(a,b)} u, \quad \text{in } (0,T) \times (0,\pi),
\\
& \partial_t y_2 - \partial_{xx} y_2 + q y_2 = \mathfrak{b}_2 \mathbf{1}_{(a,b)} u,  \quad \text{in } (0,T) \times (0,\pi),
\\
& y_j(t,0) = y_j(t,\pi) = 0, \qquad t \in (0,T), \: j = 1, 2,
\\
& y_j(0,x) = y_0^j(x), \qquad x \in (0,\pi),  \: j = 1, 2,
\end{aligned}
\right.
\end{equation}
has also been analyzed in~\cite{Ou:19} in the particular case 
\[
\operatorname{Supp} (q) \subset (0,a) 
\qquad \text{or} \qquad
\operatorname{Supp} (q) \subset (b,\pi). 
\]
Then, under conditions~\eqref{distinct} and~\eqref{fb}, she proved that $T_0(q)$ defined by~\eqref{min_time_boun} is the minimal null control time for system~\eqref{internal_ap} from $L^2(\Omega ; \R^2)$ with controls in $L^2((0,T) \times \omega \times (a,b))$. 

In the general setting, we prove the following result. 
\begin{theorem}\label{internal_result} 
Let $\Omega$, $\omega$ be as defined at the beginning of Section~\ref{s5} and let $0 \leq a < b \leq \pi$. Let $\mathfrak{b}$ and $\calA$ defined by~\eqref{Lbcoupled} with $q \in L^2 (\Omega )$ satisfying~\eqref{potential} and assume conditions~\eqref{distinct} and~\eqref{fb}.
For any $k \geq 1$ let
\[
\calG_k :=  \operatorname{Gram}_{L^2(a,b)} \left( \phi_k^{(1)}, \, \phi_k^{(2)} \right)
\]
and set 
	\begin{equation}\label{min_time_internam}
T_0( q , a, b ):=\limsup_{k \to +\infty} \frac{-\log \left| \sqrt{ \det \calG_k + \left| \nu_{k}^{(2)} - \nu_{k}^{(1)} \right|^2} \right| }{\nu_{k}^{(1)}} \in \left[ 0,\infty \right].
	\end{equation}
Then, 
\begin{enumerate}
\item if $T>T_0 (q,a,b)$, for any $y_0 \in H^{-1} (\Omega; \R^2 )$ there exists $u\in L^{2} ( (0,T) \times \omega \times (a,b))$ such that the solution $y$ of problem~\eqref{internal_ap} satisfies $y(T)=0$. 
\item if $T<T_0 (q,a,b)$, then there exists $y_0 \in H^{-1} (\Omega; \R^2 )$ such that for any $u\in L^{2} ( (0,T) \times \omega \times (a,b))$ the solution $y$ of problem~\eqref{internal_ap} satisfies $y(T)\neq 0$.
\end{enumerate}
\end{theorem}


Here again, the proof of item~$\mathit{1.}$ will follow from the moment method and particularly the use of the following biorthogonal family. 
\begin{proposition} \label{Prop:Ex_biortho_interne}
Under the previous notations, let us assume that $q \in L^2 (\Omega)$ and the sequence $\left\{ \nu_k^{(1)}, \nu_k^{(2)} \right\}_{  k \ge 1 }$ satisfies~\eqref{potential} and~\eqref{distinct}. Then, there exists a constant $\calC >0$, only depending on $q$, such that for any $T >0$, there exists a family 
\[
\left\{ \widetilde Q_{m,k}^{(j)}\right\}_{\substack{ m, k \ge 1\\  1 \le j \le 2 } } \subset L^{2}( ( 0,T ) \times \omega \times (a,b) ) 
\] 
that satisfies the biorthogonal property
\begin{equation} \label{Ex:Def_biortho_interne}
\int_0^T \int_\omega \int_a^b \widetilde Q_{m,k}^{(j)}(t,x',x) e^{-\nu_{n,\ell}^{(i)}t} \psi_n(x') \phi_{\ell}^{(i)}(x) dt dx' dx = \delta_{m n} \delta_{k \ell} \delta_{j i}
\end{equation}
for any $ m, k \ge 1$ and $j: 1 \le j \le 2 $ and
	\begin{equation}\label{esti-bound_interne}
\Ds \left\| \widetilde Q_{m,k}^{(j)} \right\|_{L^2 ((0, T) \times \omega \times(a,b))}^2 \le \calC \exp \left(\frac{ \calC }{T } \right) \exp \left( \calC \sqrt{ \nu_{m, k}^{(j)} } \right)  \frac{1}{\det \calG_k + \left| \nu_{k}^{(2)} - \nu_{k}^{(1)} \right|^2},
	\end{equation}
for any $ m, k \ge 1$ and $j: 1 \le j \le 2 $. 
\end{proposition}

\begin{proof}
Let us consider $\calU_2 = L^2(a,b)$ and $\frakC_2 $ the linear operator
	$$
\frakC_2 : \phi \in  H^2 (0, \pi) \cap H_0^1 (0, \pi) \mapsto \frakC_2 \phi =   \mathbf{1}_{(a,b)} \phi. 
	$$
The proof follows the lines of that of Proposition~\ref{Prop:Ex_biortho}. To obtain that $(\Lambda_1, \calB_1, \Lambda_2, \calB_2, \calU_2 , \frakC_2)$ satisfies Assumptions~\ref{A1} it only remains to prove the spectral inequality~\eqref{H2}. As stated in Remark~\ref{Rem:Inegalite_Spectrale} it follows from the spectral inequality~\eqref{Inegalite_Spectrale:JL}. 
Applying it for any fixed $x \in (a,b)$ with
\[
b_m = \sum_{k \leq \lambda} \sum_{j=1}^{g_k} b_{m,k}^{(j)} \phi_k^{(j)}(x)
\]
gives
\[
\int_{\Omega_1} \left| \sum_{\sqrt{\mu_m}, k \leq \lambda} \sum_{j= 1}^{g_k} b_{m,k}^{(j)} \psi_m(x') \phi_k^{(j)}(x) \right|^2 dx' 
\leq e^{\beta \lambda} \int_{\omega} \left| \sum_{\sqrt{\mu_m}, k \leq \lambda} \sum_{j= 1}^{g_k} b_{m,k}^{(j)} \psi_m(x') \phi_k^{(j)}(x) \right|^2 dx'.
\]
Integrating with respect to $x \in (a,b)$ proves~\eqref{H2}. Thus $(\Lambda_1, \calB_1, \Lambda_2, \calB_2, \calU_2 , \frakC_2)$  satisfies Assumption~\ref{A1}.

We can now apply Theorem~\ref{main} to obtain a family $\left\{ \widetilde Q_{m,k}^{(j)}\right\}_{\substack{ m, k \ge 1\\  1 \le j \le 2 } }$ satisfying the biorthogonal family~\eqref{Ex:Def_biortho_interne} and the estimate~\eqref{main-esti}. 
To prove the estimate~\eqref{esti-bound_interne} it thus only remains to estimate $\left( M_k^{-1} \right)_{j,j}$ where $M_k$ is the matrix given in~\eqref{Def:M_k}. 
For any $k \geq 1$, we consider the two possible cases. 

\begin{itemize}
\item[$\bullet$] If $\sharp G_k = 2$ then there exists $\ell \geq 1$ such that $G_k = \left\{ \nu_\ell^{(1)}, \nu_{\ell}^{(2)} \right\}$. Then, from~\eqref{Def:M_k}, using the normalization condition~\eqref{H6ter} we have 
\[
M_k = \begin{pmatrix} 1 & \left\langle \phi_\ell^{(1)}, \phi_\ell^{(2)} \right\rangle_{L^2(a,b)} \\
\left\langle \phi_\ell^{(1)}, \phi_\ell^{(2)} \right\rangle_{L^2(a,b)} & 1 + \left( \nu_{\ell}^{(1)} - \nu_{\ell}^{(2)} \right)^2 \end{pmatrix}.
\]
Explicit computations yield 
\begin{align} \notag
\det M_k &= 1 - \left\langle \phi_\ell^{(1)}, \phi_\ell^{(2)} \right\rangle_{L^2(a,b)}^2 + \left( \nu_{\ell}^{(1)} - \nu_{\ell}^{(2)} \right)^2
\\
&= \det \calG_\ell + \left( \nu_{\ell}^{(1)} - \nu_{\ell}^{(2)} \right)^2
\label{Ex:Calcul_Det_Mk}
\end{align}
and
\begin{equation} \label{Ex:Inverse_Mk_interne}
M_k^{-1} = \frac{1}{\det M_k} \begin{pmatrix} 1 + \left( \nu_{\ell}^{(1)} - \nu_{\ell}^{(2)} \right)^2 & -\left\langle \phi_\ell^{(1)}, \phi_\ell^{(2)} \right\rangle_{L^2(a,b)} \\
-\left\langle \phi_\ell^{(1)}, \phi_\ell^{(2)} \right\rangle_{L^2(a,b)} & 1 \end{pmatrix}.
\end{equation}
Thus, for any $j \in \{ 1, 2\}$, we have
\[
\left( M_k^{-1} \right)_{j,j} \leq \frac{1+\rho^2}{\det \calG_\ell + \left| \nu_{\ell}^{(1)} - \nu_{\ell}^{(2)}\right|^2}. 
\]
\item[$\bullet$] If $\sharp G_k = 1$ then there exists $\ell \geq 1$ and $i \in \{1, 2 \}$ such that $G_k = \{ \nu_\ell^{(i)} \}$. Then, from~\eqref{Def:M_k}, we have $M_k = 1$. 

As in Section~\ref{Subsec:Ex_bord}, let us formulate an estimate valid in both cases. Using estimate~\eqref{Sturm_distance_vp} and the normalization condition~\eqref{H6ter}, it comes that 
\[
\det \calG_k + \left| \nu_{k}^{(1)} - \nu_{k}^{(2)}\right|^2 \leq 1+C^2, \qquad \forall k \geq 1.
\]
Thus we obtain
\[
M_k^{-1} = 1 \leq \frac{1+C^2}{\det \calG_\ell + \left| \nu_{\ell}^{(1)} - \nu_{\ell}^{(2)}\right|^2}.
\]

\end{itemize}
Gathering both cases proves~\eqref{esti-bound_interne} and ends the proof of Proposition~\ref{Prop:Ex_biortho_interne}.
\end{proof}

We now turn to the proof of Theorem~\ref{internal_result}. 

\begin{proof}[Proof of Theorem~\ref{internal_result}]
The proof of item~$\mathit{1.}$ follows the line of the proof of Theorem~\ref{boundary_result} replacing the biorthogonal family coming from Proposition~\ref{Prop:Ex_biortho} with the one coming from Proposition~\ref{Prop:Ex_biortho_interne} and is not detailed. 

We now turn to the proof of item~$\mathit{2}$. Let $T>0$ and assume that for any $y_0 \in H^{-1} (\Omega; \R^2 )$ there exists $u\in L^{2} ( (0,T) \times \omega \times (a,b))$ such that the solution $y$ of problem~\eqref{internal_ap} satisfies $y(T)=0$. As in the proof of Theorem~\ref{internal_result}, if system~\eqref{internal_ap} is null controllable then it is also null controllable with controls in $L^2( (0,T) \times \Omega_1 \times (a,b))$. This implies
\[
T \geq \widetilde T_0(q, a,b)
\]
where $\widetilde T_0(q,a,b)$ is the minimal null control time of the one dimensional problem~\eqref{Ex:Pb_1D_interne} from $L^2(0,\pi)$ with controls in $L^2((0,T) \times (a,b))$. 

From~\cite[Theorems 11 and 18]{BM:23}, considering the particular initial condition $y_1^0 = 0$ and $y_2^0 = \phi_k^{(2)}$ for system~\eqref{Ex:Pb_1D_interne} it comes that
\[
\widetilde T_0(q,a,b) \geq \limsup_{k \to + \infty} \frac{\log \left(M_k^{-1} \right)_{2,2}}{\nu_{k}^{(2)}}. 
\]
where the matrix $M_k^{-1}$ has been computed in~\eqref{Ex:Inverse_Mk_interne}. 
Finally, using the expression of $\det M_k$ given in~\eqref{Ex:Calcul_Det_Mk} and the asymptotic~\eqref{asymp} we obtain
\[
T \geq \limsup_{k \to +\infty} \frac{-\log \left| \sqrt{ \det \calG_k + \left| \nu_{k}^{(2)} - \nu_{k}^{(1)} \right|^2} \right| }{\nu_{k}^{(1)}}
\]
which ends the proof of item~$\mathit{2}$.
\end{proof}


\section{Moment problems associated with geometrically multiple eigenvalues in $\Lambda_2$} \label{s6}

In this section we extend Theorem~\ref{main} to the case where the moment problem involves geometrically multiple eigenvalues.

As it appears in the application of our strategy to explicit examples in Section~\ref{s5}, our assumption on $\calB_2$ is only valid for geometrically simple eigenvalues in the $1D$ variable. However this assumption is not necessary and our strategy also apply with geometrically multiple eigenvalues. The price to pay is the introduction of extra heavier notation. To lighten the article we chose to present this extension and indicate the modifications in the proof in this subsection.

We stick with every assumption except for the assumption concerning $\calB_2$. We now assume that we have $\calB_2 \subset L^2(0, \pi)$, a family of $L^2(0, \pi)$, given by
	\begin{equation}\label{B2_vpm}
\calB_2 := \bigcup_{k \ge 1} \bigcup_{ 1\leq j \leq g_k } B_{k,j}, \quad B_{k,j}= \left\{ \phi_k^{(j,1)}, \dots, \phi_k^{(j,\gamma_{k,j})} \right\}  \subset H^2 (0, \pi) \cap H_0^1 (0, \pi), 
	\end{equation}
and satisfying
	\begin{equation}\label{H6_vpm}
\left\{\frakC_2 \phi_k^{(j,i)}\right\}_{1\leq i \leq \gamma_{k,j}} \hbox{is linearly independant in}\,\, \calU_2 , \quad \forall k \ge 1 \hbox{ and } j : 1 \le j \le g_k.  
	\end{equation}

More precisely, our main assumption is now the following. 
\begin{assumption}  \label{A1_vpm}
We have two positive real sequences $\Lambda_1 $ and $\Lambda_2$, an orthonormal basis $\calB_1$ of $L^2 (\Omega_1)$, a sequence $\calB_2$ of $L^2 (0, \pi)$, a Hilbert space $\calU_2$, and an operator $\frakC_2 \in \calL ( H^2 (0, \pi) \cap H_0^1 (0, \pi), \calU_2)$ such that
	\begin{equation*}
	\left\{
	\begin{array}{l}
\Ds  \Lambda_1  \hbox{ satisfy~\eqref{H1} and with } \kappa_1, \theta_1 \hbox{ and } \vartheta \in (0, 1); \\
	\noalign{\smallskip}
 \Ds \Lambda_2 \in \calL (p, \rho, \theta , \kappa), \ p \in \N, \ \rho, \kappa > 0 \hbox{ and } \theta \in (0, 1), \hbox{ satisfying~\eqref{Gk}, \eqref{Gk2}}; \\
	\noalign{\smallskip}
\Ds \calB_2 \hbox{ is given by~\eqref{B2_vpm} and satisfies~\eqref{H6_vpm}}; 
\hbox{ inequality~\eqref{H2} holds with } \beta >0.
	\end{array}
	\right.
	\end{equation*}
\end{assumption}
Let us consider the sequence 
\[
\calF := \left\{ F_{m,k}^{(j,i)} : m, k \ge 1,\  1 \le j \le g_k, \ 1\leq i \leq \gamma_{k,j} \right \}
\] 
of elements of $\calU_2$ given by
	\begin{equation}\label{Fmk_vpm}
F_{m,k}^{(j,i)} (t, x') := e_{m, k}^{(j,i)} (t) \psi_{m} (x') = e^{-\lambda _{m,k}^{ ( j ) }t} \psi_{m} ( x') \frakC_2 \phi_k^{(j,i)} , \quad ( t,x'   ) \in (0,T) \times \Omega_1 , 
	\end{equation}
We now consider the matrix $M_k$ given by
\begin{equation} \label{Def:M_k_vpm}
M_k = \sum_{\ell=1}^{g_k} \operatorname{Gram}_{\calU_2} \left(
\delta_{k,\ell}^1 \frakC_2 \phi_k^{(1, 1)}, \dots, \delta_{k,\ell}^1 \frakC_2 \phi_k^{(1, \gamma_{k,1})},
\dots,
\delta_{k,\ell}^{g_k} \frakC_2 \phi_k^{(g_k, 1)}, \dots, \delta_{k,\ell}^{g_k} \frakC_2 \phi_k^{(g_k, \gamma_{k,g_k})}
\right)
\end{equation}
where $\delta_{k,\ell}^j$ is given by~\eqref{Def:Delta_k}. We consider the associated renumbering function
\begin{equation} \label{Renumbering}
R : (j,i) \in \N^* \times \N^* \mapsto \gamma_{k,1} + \dots + \gamma_{k,j-1} + i,
\end{equation}
with the convention $\gamma_{k,0} = 0$. 

With these notation we obtain the following theorem.
 \begin{theorem}\label{main_vpm}
Let us assume that $(\Lambda_1, \calB_1, \Lambda_2, \calB_2, \calU_2 , \frakC_2)$ satisfies Assumption~\ref{A1_vpm}. Then, there exists a constant $\calC >0$, only depending on $p$, $\rho$, $\theta$, $\kappa$, $\beta$, $\vartheta$, $\theta_1$ and $\kappa_1$, such that for any $T >0$, the sequence $\calF $ (see~\eqref{Fmk_vpm}) admits a biorthogonal family 
\[
\left\{ Q_{m,k}^{(j,i)} \: : \:  m, k \ge 1,\, 1 \le j \le g_k,\,\,1\leq i \leq\gamma_{k,j} \right \}
\] 
in $L^{2}( ( 0,T ) \times \omega ; \calU_2)$, \textit{i.e.}, such that for any $ m,n \ge 1$, any $k , \ell \geq 1$, any $j : 1 \le j \le g_k $ any $j': 1 \le j' \le g_\ell$, any $i : 1 \le i \le \gamma_{k,j} $ and any $i': 1 \le i' \le \gamma_{\ell,j'}$, we have 
\[
\int_0^T \int_\omega \psU{ Q_{m,k}^{(j,i)}(t,x') }{ F_{n,\ell}^{(j',i')}(t,x') } dx' dt = \delta_{mn} \delta_{k \ell} \delta_{j j'} \delta_{i i'},
\] 
that satisfies
	\begin{equation}\label{main-esti_vpm}
\Ds \left\| Q_{m,k}^{(j,i)} \right\|^2_{L^2 ((0, T) \times \omega; \calU_2) } \le \calC \exp \left(\frac{ \calC }{T^b } + \frac{ \calC} {T^{ \theta '} } \right) \exp \left( \calC \left[ \lambda_{m, k}^{(1)} \right]^{ \frac{ b}{1 + b}} + \calC \left[ \lambda_{m , k }^{(1)} \right]^\theta \right)  \left( M_k^{-1} \right)_{R(j,i),R(j,i)},
	\end{equation}
for any $ m, k \ge 1$, any $j: 1 \le j \le g_k $ and any $i : 1\leq i \leq \gamma_{k,j}$, where $M_k$ is the matrix defined in~\eqref{Def:M_k_vpm}, $\theta '$ and $b$ are given by~\eqref{theta} and \eqref{bcoef} and $R$ is the renumbering function defined in~\eqref{Renumbering}. 
\end{theorem}

\begin{proof}
The proof of this theorem follows the same steps as that of Theorem~\ref{main}.  Let us briefly explain the necessary adjustments. 

\begin{itemize}
\item[$\bullet$] \textbf{Theorem~\ref{tPN}}. Following the lines of the proof of Theorem~\ref{tPN} directly gives
	\begin{equation}\label{f7_vpm}
\Ds  \int_{0}^T  \!\! \int_{ \Omega_1 }  e^{-\frac{\alpha \beta }{t^{ b} }} \normeU{ P_{N} (t,x') }^{2} \, dx'\,dt  \leq  6 \int_0^T   \!\! \int_{ \omega } \normeU{ P_{N}(t,x') }^{2} \, dx'\,dt , 
	\end{equation}
for any $ T \in (0, \tau_0]$, any $N \ge 1$ and any $P_N$ given by
	\begin{equation} \label{PN_vpm}
P_{N} (t, x')  := \sum_{\mu _{m}^\vartheta, k \leq N} \sum_{j= 1 }^{g_k}  \sum_{i= 1 }^{\gamma_{k,j}} a_{m,k}^{ (j,i , N ) } F_{m,k}^{(j,i)} (t, x'), \qquad (t, x' ) \in (0, T) \times \Omega_1. 
	\end{equation}
where $\beta $ is the constant in~\eqref{H2} and $b$ is given by~\eqref{bcoef}. Indeed the key estimate given by Lemma~\ref{l2} still holds with the same proof taking into account the new definition of $g_{m,k}^{ ( N ) }$ given by
\[
g_{m,k}^{ ( N ) } (t) := \sum_{j= 1 }^{g_k} e^{-\lambda _{m,k}^{ (j) }t}  \sum_{i= 1 }^{\gamma_{k,j}} a_{m,k}^{ (j,i,N) } \frakC_2 \phi_k^{(j,i)} \in \calU_2.
\]
\item[$\bullet$] \textbf{Theorem~\ref{tmainrest}}. Replacing the spaces $E_{\eta_\alpha}$ and $E^\omega$ by
	\begin{equation*}
	\begin{array}{c}
\Ds E_{\eta _{\alpha }} =\overline{\mathrm{span} \, \left \{ F_{m, k}^{(j,i)} : k,m\geq 1, \ j : 1 \le j \le g_k, \ i : 1 \leq i \leq \gamma_{k,j} \right\}  }^{L_{\eta _{\alpha }}^{2} ((0, T) \times \Omega_1 ; \calU_2)} , \\
	\noalign{\smallskip}
\Ds \Ds E^\omega =\overline{\mathrm{span}\, \left \{ F_{m, k}^{(j,i)} |_\omega : k,m\geq 1, \ j : 1 \le j \le g_k, \ i : 1 \leq i \leq \gamma_{k,j} \right\}  }^{L^{2} ((0, T) \times \omega ; \calU_2) },
	\end{array}
	\end{equation*}
we obtain that Theorem~\ref{tmainrest} holds without any modification. Notice that neither Theorem~\ref{tPN} nor Theorem~\ref{tmainrest} uses the assumption~\eqref{H6} (here replaced by~\eqref{H6_vpm}). 
\smallskip
\item[$\bullet$] \textbf{Theorem~\ref{main}}. Finally the proof of Theorem~\ref{main} combines two steps: the existence (with estimates) of a biorthogonal family in $L^2((0,T) \times \Omega_1 ; \calU_2)$ and the isomorphism property of the restriction operator coming from Theorem~\ref{tmainrest}. The latter and its use to deduce a biorthogonal family in $L^((0,T) \times \omega ; \calU_2)$ given in Proposition~\ref{Prop:restriction_biortho} remains unchanged. 

The existence with suitable estimates of a biorthogonal family under Assumption~\eqref{A1_vpm} follows the line of Section~\ref{ss42} replacing the use of Proposition~\ref{Prop:Annexe_biortho} by~\cite[Theorem 51]{BM:23}. Here it is necessary to assume~\eqref{H6_vpm}. 
\end{itemize}
\end{proof}
Notice that the use of such biorthogonal families allows to prove that the minimal null control time for system~\eqref{internal_ap} given in Theorem~\ref{internal_result} still holds replacing condition~\eqref{distinct} by the approximate controllability of system~\eqref{internal_ap}.


\appendix
\section{Biorthogonal families in $L^2(0,T ; \calU_2)$}\label{Annexe1}

Recall that we have defined in~\eqref{ekm} the following functions
\[
e_{m, k}^{(j)}  : t \mapsto e^{- \lambda_{m, k}^{( j)} t} \frakC_2 \phi_k^{(j)} \in \calU_2, \quad \forall m, k \geq 1 \hbox{ and } j:  1\le j \le g_k. 
\]
In the setting considered in Assumption~\ref{A1}, using results proved in~\cite{BM:23}, we obtain a biorthogonal family to 
$\left\{ e_{m,k}^{(j)} \right\}_{ \substack{ k \ge 1 \\  1 \le j \le g_k }  }$
in $L^2( 0,T ; \calU_2)$. 
More precisely, we obtain the following result. 
\begin{proposition} \label{Prop:Annexe_biortho}
Let $\Lambda_1 \subset (0,\infty)$. 
Let $\left( \Lambda_2, \calB_2 , \calU_2, \frakC_2 \right)$ satisfying $\Lambda_2 \in \calL (p, \rho, \theta , \kappa)$,~\eqref{B2} and~\eqref{H6}. Let $(G_k)_{k \geq 1}$ be a grouping satisfying~\eqref{Gk} and~\eqref{Gk2}.
For any $k \geq 1$, let $M_k$ be the matrix defined by~\eqref{Def:M_k}.
There exists a positive constant $ \widetilde C_ 1 $ depending on $p$, $\rho$, $\theta$ and $\kappa$ such that for any $T > 0$, for any $m \geq 1$, there exists a biorthogonal family  $\left\{ q_{m,k}^{(j)} \right\}_{\substack{ k \ge 1 \\  1 \le j \le g_k } } $ to $\left\{ e_{m,k}^{(j)} \right\}_{\substack{ k \ge 1 \\  1 \le j \le g_k } } $ in $ L^2(0,T ;\calU_2)$, \textit{i.e.}, such that
	\begin{equation*}
\Ds \int_0^T \psU{ q_{m,k}^{(j)}(t)}{ e^{-(\lambda_\ell^{(i)} +\mu_m)t} \frakC_2 \phi_{\ell}^{(i)} } \, dt = \delta_{k \ell} \delta_{ij}, \quad \forall k , \ell \ge 1 \hbox{ and } 1 \le j \le g_k, \, 1 \le i \le g_{\ell} .
	\end{equation*}
satisfying the following estimate
	\begin{equation}\label{H7}
\left\| q_{m,k}^{(j)}  \right\|_{L^2(0,T ; \, \calU_2)}^2 \leq \widetilde C_1 \exp \left( \frac{\widetilde C_1} {T^{ \theta '} } \right) e^{\widetilde C_1 \left[ \lambda_{k}^{(1)}+\mu_m \right]^\theta } \left( M_k^{-1} \right)_{j,j} , \quad \forall k \ge 1 \hbox{ and } 1 \le j \le g_k ,
	\end{equation}
with $\theta ' $ given by~\eqref{theta}. 

Moreover, there exists a positive constant $\widetilde C_2$ depending on $p$ and $\inf \Lambda_2$ such that for any $T>0$, for any $m \geq 1$, any biorthogonal family  $\left\{ q_{m,k}^{(j)} \right\}_{\substack{ k \ge 1 \\  1 \le j \le g_k } } $ to $\left\{ e_{m,k}^{(j)} \right\}_{\substack{ k \ge 1 \\  1 \le j \le g_k } } $ in $ L^2(0,T ;\calU_2)$ satisfies
\[
\left\| q_{m,k}^{(j)}  \right\|_{L^2(0,T ; \, \calU_2)}^2 \geq \widetilde C_2 \left( M_k^{-1} \right)_{j,j} , \quad \forall k \ge 1 \hbox{ and } 1 \le j \le g_k. 
\]
\end{proposition}

This proposition is exactly~\cite[Theorem 51]{BM:23}. The only assumption to check to apply~\cite[Theorem 51]{BM:23} is that 
\[
\Lambda_2^{(m)} := \mu_m + \Lambda_2 = \left\{ \lambda_k^{(j)} + \mu_m \right\}_{  \substack{ k \ge 1 \\  1 \le j \le g_k } } 
\in \calL (p, \rho, \theta , \kappa)
\]
for any $m \geq 1$ which we proved in the proof of Corollary~\ref{cblock}. 

Notice that in the particular case $\calU_2 = \R$, then we are dealing with classical biorthogonal families to time exponentials and, in this setting, Proposition~\ref{Prop:Annexe_biortho} is a consequence of~\cite{GBO:19}.


\section{A Lebeau-Robbiano construction} \label{Annexe2}

In this appendix we revisit the classical Lebeau-Robbiano strategy from the point of view of biorthogonal families. Using the restriction operator of Section~\ref{s3}, we prove in Theorem~\ref{tmainrest3} that the spectral inequality implies the existence of biorthogonal families with estimates allowing to recover null controllability in arbitrary time without any geometrical extra condition on the space domain. 

In~\cite{M:10} L.~Miller was interested in an adaptation of the Lebeau-Robbiano strategy for the proof of an inequality of observability of
heat-like semigroups. In particular, he considered the Dirichlet-Laplace operator in $L^{2}(\Omega )$, where $\Omega $ is a sufficiently smooth bounded domain of $\mathbb{R}^{d}$ with $d \ge 1$. Let us denote by $\left\{ e^{tA} \right\}_{t\geq 0}$, the semigroup generated by this operator in $L^{2}(\Omega )$ and $\left\{ \mu _{m},\psi _{m} \right\}_{m\geq 1}$, its eigenelements ($ \| \psi_m \|_{L^2(\Omega)} = 1$).  Using that $\left\{ \psi _{m} \right\}_{m\geq 1}$ satisfy~\eqref{Inegalite_Spectrale:JL}, with $ \Omega_{1}$ replaced by $\Omega $ and $\omega \subset \Omega $ an arbitrary nonempty open subset, i.e., using
	\begin{equation} \label{NSpI}
\int_{\Omega }\left \vert \sum_{\sqrt{\mu _{m}}\leq \lambda }b_{m}\psi _{m}(x' )\right \vert ^{2}\, dx'  \leq e^{\beta \lambda } \int_{\omega }\left \vert \sum_{\sqrt{\mu _{m}}\leq \lambda }b_{m}\psi _{m}(x' )\right \vert ^{2}\, dx' ,
	\end{equation}
for any $\{ b_m \}_{m \ge 1} \in \ell^2$, L.~Miller proved that, for all $T>0,$ the observability inequality
	\begin{equation}\label{IOchaleurLR}
\left \Vert e^{tA}f \right\Vert _{L^{2}(\Omega )}^{2}\leq \calK (T,\omega )\int_{0}^{T} \! \! \int_{\omega }\left \vert e^{tA}f\right \vert ^{2}, \quad \forall f \in L^2 (\Omega),
	\end{equation}
holds for a positive constant $\calK (T,\omega )$ satisfying
	\begin{equation} \label{CTHLR}
\limsup_{T \to 0 }\left( T\ln (\calK (T,\omega ))\right) \leq 2\beta ^{2}. 
	\end{equation}
It is well-known that the observability inequality~\eqref{IOchaleurLR} is equivalent to the null controllability property for the heat equation at time $T > 0$ together with an estimate of the associated control cost: for all $y_{0}\in L^{2}(\Omega )$, there exists $u\in L^{2}\left( Q_T \right) $ ($Q_T := (0,T) \times \Omega$) satisfying 
	\begin{equation}\label{coutuLR}
\Vert u \Vert _{L^{2}((0,T)\times \omega )}^{2}\leq \calK (T,\omega )\Vert y_{0}\Vert _{L^{2}(\Omega )}^{2},
	\end{equation}
with $ \calK (T,\omega )$ the constant in~\eqref{IOchaleurLR}, and such that the problem 
	\begin{equation}\label{controlchaleurd}
	\left \{ 
	\begin{array}{ll}
\Ds \partial _{t}y-\Delta y=1_{\omega }u,  & \hbox{in } Q_T , \\ 
	\noalign{\smallskip}
\Ds y=0, & \hbox{on } \Sigma_T:= (0,T) \times \partial \Omega , \\ 
	\noalign{\smallskip} 
\Ds y(0,\cdot ) =y_{0},  \quad  y(T,\cdot )=0, & \hbox{in } \Omega ,
	\end{array}
	\right. 
	\end{equation}
admits a weak solution $y $.

The purpose of this appendix is to prove the null controllability result~\eqref{controlchaleurd}, with controls $u$ satisfying~\eqref{coutuLR} for a constant $\calK (T,\omega ) $ fulfiling~\eqref{CTHLR}. To be precise, we will prove:


\begin{theorem}\label{mainLR}
Let us assume that $\Lambda _{1}=\{ \mu _{m}\}_{m\geq }$ satisfies~\eqref{H1} for $\theta_1 > 0$ and that $\mathcal{B}_{1}=\{ \psi _{m}\}_{m\geq 1}$ is an orthonormal basis of $L^{2}(\Omega )$ that satisfies \eqref{NSpI}.  Let $T>0$. Then, there exists a constant $\calK (T, \omega) > 0$ satisfying~\eqref{CTHLR} and such that the null controllability problem~\eqref{controlchaleurd} has a solution $u \in L^2 ((0,T) \times \omega)$ satisfying~\eqref{coutuLR}.
\end{theorem}

To this end, we will use a different approach to that of~\cite{M:10}: we will solve~\eqref{controlchaleurd} by solving the associated moment problem. This will entail the construction of a biorthogonal family in $L^{2}((0,T)\times \omega )$ to the sequence $\left\{ F_m \right\}_{m\geq 1}$, given by
	\begin{equation}\label{FmLR}
F_{m}(t,x)=e^{-\mu _{m}t}\, \psi _{m}(x),\quad (t,x)\in Q_T , \quad m\geq 1,   
	\end{equation}
with an explicit estimate with respect to $T > 0$ of the norms of the elements of this family.

In order to construct an appropriate biorthogonal family in $L^{2}((0,T)\times \omega )$ to $\left\{ F_m \right\}_{m\geq 1}$ we will use the restriction argument of the previous sections. In fact, this restriction argument will allow us to construct the unique optimal biorthogonal family to $\{ F_{m} \}_{m \ge 1}$ in $L^{2}((0,T)\times \omega )$.

As a first step, we begin by the following result:


\begin{proposition}\label{Prop_epsLR} 
For any $\varepsilon \in \left( 0,T\right) $, there exists a sequence $\left\{f_{m}^{\varepsilon } \right\}_{m \geq 1}$ biorthogonal  in $L^{2}(Q_T)$ to $\left\{F_{m} \right\}_{m\geq 1}$, $F_m $ given by~\eqref{FmLR}, satisfying
	\begin{equation}\label{fepsestimLR}
	\left \{ 
	\begin{array}{l}
\mathrm{Supp} \, \left( f_{m}^{\varepsilon }\right) \subset \left[ \varepsilon ,T\right] \times \overline \Omega , \\ 
	\noalign{\smallskip}
\Ds \left\| f_{m}^{\varepsilon } \right\|_{L^{2}(Q_T)}^{2}\leq c\left( \dfrac{1}{T-\varepsilon }+\mu _{m}\right) e^{2\varepsilon \mu _{m}},%
	\end{array}
	\right.  
	\end{equation}
for any $ m \ge 1$, with $c=2e$.
\end{proposition}


\begin{proof}
Let $\tau >0$. Using the orthonormality of the sequence $\{ \psi_{m}\}_{m\geq 1}$ in $L^{2}(\Omega )$, a natural choice of a biorthogonal
family $\{f^\tau_{m}\}_{m\geq 1}$ to $\{F_{m}\}_{m\geq 1}$ in $L^{2}\left( ( 0,\tau ) \times \Omega \right) $ is a sequence of the form 
	\begin{equation*}
g^\tau_{m} (t,x) = \frac 1{C^\tau_m} e^{- \mu_m t} \psi _{m}(x),\quad \mathrm{with} \ C_m^\tau = \int_{0}^{\tau } e^{-2\mu_{m}t}\,dt  ,\quad m\geq 1.
	\end{equation*}
Now, for any $\varepsilon \in \left( 0,T\right)$, we choose $\tau =T-\varepsilon $ and we define
	\begin{equation*}
f_{m}^{\varepsilon }(t,x) =
	\left\{
	\begin{array}{ll}
\Ds 0 &\mathrm{if} \ t \in (0, \varepsilon], \quad x \in \Omega , \\
	\noalign{\smallskip}
\Ds e^{\varepsilon \mu _{m}}g_{m}^{T-\varepsilon }( t-\varepsilon, x ) &\mathrm{if} \ t\in  ( \varepsilon ,T ) , \quad x \in \Omega, 
	\end{array}
	\right.
	\end{equation*}
for any $m \ge 1$. It is not difficult to check that the sequence $\left\{ f_{m}^{\varepsilon }\right\}_{m\geq 1}$ is biorthogonal to $\left\{ F_{m}\right\}_{m\geq 1}$ in $L^{2}(Q_T) $ and satisfies the first condition in~\eqref{fepsestimLR}. On the other hand,
	\begin{equation*}
\left\| f_{m}^{\varepsilon } \right\|_{L^{2}(Q_T )}^{2} = \frac{e^{2\varepsilon \mu _{m}}}{C_m^{T - \varepsilon}}  , \quad \forall m\geq 1.
	\end{equation*}

If $2\mu _{m} ( T-\varepsilon  ) \leq 1$, then
	\begin{equation*}
C_m^{T - \varepsilon} = \int_{0}^{T-\varepsilon }e^{-2\mu _{m}t} \, dt\geq \int_{0}^{T-\varepsilon }e^{-1}\, dt = e^{-1}\left( T-\varepsilon \right) , \quad \forall m\geq 1.
	\end{equation*}
If $2\mu _{m} ( T-\varepsilon  ) >1$, then
	\begin{equation*}
C_m^{T - \varepsilon} = \int_{0}^{T-\varepsilon }e^{-2\mu _{m}t} \, dt \geq \int_{0}^{\frac{1}{2\mu _{m} }} e^{-2\mu _{m}t} \, dt\geq \frac{e^{-1}}{2\mu _{m}} , \quad \forall m\geq 1.
	\end{equation*}
The two previous inequalities together with the expression of $\left\| f_{m}^{\varepsilon } \right\|_{L^{2}(Q_T )}$ proves the second condition in~\eqref{fepsestimLR}). This ends the proof of the result.
\end{proof}

For $N\geq 1$, let us define 
	\begin{equation}  \label{P_Nouvert}
P_{N} ( t,x ) := \sum_{\sqrt{\mu _{m}}\leq N}a_{m}^{( N ) }F_{m} ( t,x ) =\sum_{\sqrt{\mu _{m}}\leq N}a_{m}^{(N)}e^{-\mu
_{m}t}\psi _{m} ( x ) ,
	\end{equation}
where $a_{m}^{( N ) } \in \R$, $1 \le m \le N$. The second step in our approach is:


\begin{proposition}\label{Prop_restrictouvertLR}
For all $\alpha > 2\beta $ ($\beta $ is the constant in~\eqref{NSpI}), any integer $N\geq 1$ and any $P_N$ (see~\eqref{P_Nouvert}), one has
	\begin{equation}\label{Eq_3ouvertLR}
\intdouble_{Q_T } e^{-\frac{ \alpha \beta}{t}} \left \vert P_{N} ( t,x ) \right \vert^{2}\, dx\, dt \leq \left( 3 \left\| P_{N} \right \|_{L^2 ((0,T) \times \omega )}^{2} + \mathcal{M}  ( \alpha ,T ) \intdouble_{Q_T } e^{-\frac{ \alpha \beta}{t}} \left \vert P_{N} ( t,x ) \right \vert^{2}\, dx\, dt \right) ,  
	\end{equation}
where 
	\begin{equation*}
	\left\{
	\begin{array}{l}
\Ds \mathcal{M} ( \alpha ,T ) = \widehat{\mathcal{C}} \, \frac{\chi ( T,\alpha  ) }{\left( \alpha ^{2} - 2\alpha \beta \right)^{2\theta _{1}+1} } e^{-\frac{\alpha ^{2}-2\alpha \beta }{2T}},  \\
	\noalign{\smallskip}
\Ds \chi  ( T,\alpha  ) = \left( \frac{1}{T}+1\right)\left[ T^{\theta _{1} + 1} \left( 1 + T^{\theta _{1}+1} \right) + \alpha^{2\theta_1+2} \right],
 	\end{array}
 	\right.
	\end{equation*}
$\theta_1 = {d}/{2}$ and $\widehat{\mathcal{C}} $ is a positive constant only depending on $\theta_1 $. 

\end{proposition}


\begin{proof}
The proof follows the steps of the proofs of Lemmas~\ref{lrest1} and~\ref{lrest2} ($b = 1$). Given $ N \ge 1$, $\alpha$ and $T >0$, we will assume $T > \alpha / N  $. See Remark~\ref{r4} when $T \le \alpha / N  $.

Recall that for all $N\geq 1$ and $P_N $ given by~\eqref{P_Nouvert}, one has 
	\begin{equation*}
\int_{\Omega }\left \vert P_{N} ( t,x ) \right \vert ^{2}\, dx = \sum_{\sqrt{\mu _{m}}\leq N} \left| a_{m}^{N} \right|^{2} e^{-2\mu _{m}t}, \quad \forall t \in (0, T).
	\end{equation*}
The spectral inequality~\eqref{NSpI} gives (see the proof of Lemma~\ref{lrest1}) 
	\begin{equation} \label{Eq_4ouvertLR}
\int_{0}^{\frac{\alpha }{N}} \! \! \int_{\Omega }e^{-\frac{\alpha \beta }{t}}\left \vert P_{N} ( t,x ) \right \vert ^{2}\, dx\,dt \leq \int_{0}^{\frac{ \alpha }{N}} \! \! \int_{\omega }\left \vert P_{N} ( t,x ) \right \vert^{2}\, dx\,dt. 
	\end{equation}

Reasoning as in the proof of Lemma~\ref{lrest2} (see~\eqref{Eq_2} with $b = 1$ and $\vartheta = 1/2 $), on the interval $\left( \alpha / N  ,T\right)  $, one has: 
	\begin{equation} \label{maj_1}
e^{-\frac{ \alpha \beta }{t} }\int_{\Omega }\left \vert P_{N} ( t,x ) \right \vert ^{2}\, dx \leq 3\left( \int_{\omega }\left \vert P_{N} ( t,x ) \right \vert ^{2}\, dx + \sum_{\frac{\alpha }{t} < \sqrt{\mu _{m}}\leq N}|a_{m}^{(N)}|^{2}e^{-2\mu _{m}t}\right) ,
	\end{equation}
for any $t \in \left( \alpha / N  ,T\right) $. Now, using Proposition~\ref{Prop_epsLR} with $\varepsilon \in (0, T)$, we obtain 
	\begin{equation*}
	\begin{array}{l}
\Ds \left| a_{m}^{N} \right|^2 = \left| \int_0^T  \! \! \int_{\Omega }  P_{N} (t,x) f_{m}^{\varepsilon } (t, x) \, dx\, dt \right|^2 = \left| \int_\varepsilon^T  \! \! \int_{\Omega }  e^{- \frac{\alpha \beta}{2t} }P_{N} (t, x) e^{ \frac{\alpha \beta}{2t} } f_{m}^{\varepsilon } (t, x) \, dx\, dt \right|^2 \\
	\noalign{\smallskip}
\Ds \phantom{ \left| a_{m}^{N} \right|^2} \leq c e^{ \frac{\alpha \beta }{\varepsilon }}\left( \dfrac{1}{T-\varepsilon }+\mu _{m}\right) e^{2\varepsilon \mu _{m}} \intdouble_{Q_T } e^{-\frac{ \alpha \beta}{t}} \left \vert P_{N} ( t,x ) \right \vert^{2}\, dx\, dt. 
	\end{array}
	\end{equation*}
Let us choose $\varepsilon = t/2 \in (0, T)$. Thus, the previous inequality implies
	\begin{equation} \label{jarto}
	\begin{array}{l}
\Ds \sum_{\frac{\alpha }{t}<\sqrt{\mu _{m}}\leq N} \left| a_{m}^{(N)} \right|^{2} e^{-2\mu _{m}t} \leq c e^{\frac{2 \alpha \beta }{t } } \sum_{\frac{\alpha }{t}<\sqrt{\mu _{m}}} \left(\frac 2{2T-t}+\mu _{m}\right)e^{-\mu _{m}t} \intdouble_{Q_T } e^{-\frac{ \alpha \beta}{t}} \left \vert P_{N} ( t,x ) \right \vert^{2}\, dx\, dt, \\
	\noalign{\smallskip}
\Ds \phantom{\sum_{\frac{\alpha }{t}<\sqrt{\mu _{m}}\leq N} \left| a_{m}^{(N)} \right|^{2} e^{-2\mu _{m}t} } \leq e^{\frac{2 \alpha \beta }{t } } S ( t,\alpha ) \intdouble_{Q_T } e^{-\frac{ \alpha \beta}{t}} \left \vert P_{N} ( t,x ) \right \vert^{2}\, dx\, dt, 
	\end{array}
	\end{equation}
for any $t \in \left( \alpha / N  ,T\right) $, with 
	\begin{equation*}
S ( t,\alpha  ) = 2 c \sum_{\frac{\alpha }{t}<\sqrt{\mu _{m}}} \left(\frac 1T +\mu _{m}\right)e^{-\mu _{m}t}. 
	\end{equation*}

In what follows $\widehat{\mathcal{C}} $ will denote a generic positive constant only depending on $\theta_1$ and whose value can change from one line to the next. 

We can bound the series appearing in $S(t, \alpha)$ using Lemma~\ref{l1}. Indeed, the sequence $\left \{ \mu _{m}\right \} _{m\geq 1}$ satisfies~\eqref{H1} for $\kappa_1>0$, only depending on $\Omega \subset \R^{d}$ and $d$, and $\theta_1 = d / 2$ (Weyl's law). From inequality~\eqref{sum} applied to $q=\theta _{1}$, $\gamma = \frac{\alpha ^{2}}{t^{2}}>0$ and $\sigma =t>0$, we can write 
\begin{equation*}
	\mathcal{S}_{1}(t) :=\sum_{\frac{\alpha }{t}<\sqrt{\mu _{m}}} e^{-\mu_{m}t} \leq \widehat{\mathcal{C}} \, \frac{t^{\theta _{1}}+\alpha ^{2\theta_{1}}}{t^{2\theta _{1}}} e^{-\frac{\alpha ^{2}}{t}}, \quad \forall t \in(0,T). 
	\end{equation*}
An adaptation of the proof of Lemma~\ref{l1} leads to the existence of a new positive constant $\widehat{\mathcal{C}}$ such that
	\begin{equation*}
\mathcal{S}_{2}(t):= \sum_{\frac{\alpha }{t}<\sqrt{\mu _{m}}} \mu_{m} e^{-\mu _{m}t} \leq \widehat{\mathcal{C}} \, \frac{t^{\theta_{1}+1} + \alpha ^{2\theta _{1}+2}}{t^{2\theta _{1}+2}} e^{-\frac{\alpha ^{2}}{t}} , \quad \forall t \in(0,T).
	\end{equation*}
Then, for a new constant $\widehat{\mathcal{C}} >0$  we deduce
	$$
	\begin{array}{l}
\Ds S(t ,\alpha) \le  \widehat{\mathcal{C}} \, \left( \frac{1}{T} +1\right) \left( \frac{t^{\theta _{1}+2}+ \alpha ^{2\theta_{1}}t^{2} + t^{\theta _{1}+1}+\alpha ^{2\theta _{1}+2}}{t^{2\theta _{1}+2}} \right) e^{-\frac{\alpha ^{2}  }{t}} \\
	\noalign{\smallskip}
\phantom{S(t ,\alpha)} \Ds \le \widehat{\mathcal{C}} \, \left( \frac{1}{T} +1\right) \left( \frac{T^{\theta _{1}+2}+ \alpha ^{2\theta_{1}}T^{2} + T^{\theta _{1}+1}+\alpha ^{2\theta _{1}+2}}{t^{2\theta _{1}+2}} \right) e^{-\frac{\alpha ^{2}  }{t}} , \quad \forall t \in(0,T). 
	\end{array}
	$$
In order to get a simpler estimate of $S (t, \alpha)$, we will use Young's inequality as follows: 
	$$
	\left\{
	\begin{array}{l}
\Ds T^a = 1 \cdot T^a \le \frac{b-a}b 1^{\frac b{b - a}} + \frac ab T^b \le 1 + T^b , \quad \forall a,b: 0 < a < b, \\
	\noalign{\smallskip}
\Ds \alpha^{2\theta_{1}} T^{2} \le \frac{\theta_1}{\theta_1 + 1} \alpha^{2 \theta_1 + 2 } + \frac{1}{\theta_1 + 1} T^{2 \theta_1 + 2 } \le \alpha^{2 \theta_1 + 2 } + T^{2 \theta_1 + 2 },
	\end{array}
	\right.
	$$
for any $ T >0$. Thus,
	$$
	\left\{
	\begin{array}{l}
\Ds {T^{\theta _{1}+2}+ \alpha ^{2\theta_{1}}T^{2} + T^{\theta _{1}+1}+\alpha ^{2\theta _{1}+2}} \le T^{\theta_1 + 1} \left(T + T^{\theta_1 + 1} +1 \right) + 2\alpha ^{2\theta _{1}+2}  \\
	\noalign{\smallskip}
\phantom{T^{\theta _{1}+2}+ \alpha ^{2\theta_{1}}T^{2} + T^{\theta _{1}+1}+\alpha ^{2\theta _{1}+2}} \Ds \le2 \left( T^{\theta_1 + 1} \left(1 +T^{\theta_1 + 1} \right) + \alpha ^{2\theta _{1}+2}\right). 
	\end{array}
	\right.
	$$

Coming to the last estimate of $S (t, \alpha )$ and~\eqref{jarto}, we obtain
	\begin{equation}\label{somme_maj}
\Ds \sum_{\frac{\alpha }{t}<\sqrt{\mu _{m}}\leq N} \left| a_{m}^{(N)} \right|^{2} e^{-2\mu _{m}t}  \leq  \widehat{\mathcal{C}} \, \chi  ( T , \alpha ) h(\alpha ,t) \intdouble_{Q_T } e^{-\frac{ \alpha \beta}{t}} \left \vert P_{N} ( t,x ) \right \vert^{2}\, dx\, dt, 
	\end{equation}
for any $t \in  (  \alpha / N , T  )$, where $\widehat{\mathcal{C}}  > 0$ is a new constant, $ \chi  ( T , \alpha )$ is given in the statement of Proposition~\ref{Prop_restrictouvertLR} and 
	\begin{equation*}
h(\alpha ,t):=t^{-\left( 2\theta _{1}+2\right) }e^{-\frac{\alpha^{2}-2\alpha \beta }{t}}, \quad t \in (0,T).
	\end{equation*}
The function $h$ is clearly bounded on $\left( 0,T\right) $ if $\alpha >2\beta $. Therefore, 
	\begin{equation*}
	\begin{array}{l}
\Ds \int_{0}^{T}h(\alpha ,t)\,dt \leq 2\max_{t\in  ( 0,T ) }\left(  \frac{t^{-2\theta _{1}}e^{-\frac{\alpha ^{2}-2\alpha \beta }{2t}}}{\alpha^{2}-2\alpha \beta }\right) \int_{0}^{T}\left( \frac{\alpha ^{2}-2\alpha \beta }{2t^{2}}e^{-\frac{\alpha ^{2}-2\alpha \beta }{2t}}\right) \, dt \\
	\noalign{\smallskip}
\Ds \phantom{ \int_{0}^{T}h(\alpha ,t)\,dt} \leq \widehat{\mathcal{C}} \, \frac{e^{-\frac{\alpha ^{2}-2\alpha \beta }{ 2T}}}{\left( \alpha ^{2}-2\alpha \beta \right) ^{2\theta _{1}+1}},
	\end{array}
	\end{equation*}
with $\widehat{\mathcal{C}} > 0$. Going back to~\eqref{somme_maj}, we get
	\begin{equation*}
\int_{\frac \alpha N}^{T} \sum_{\frac{\alpha }{t}<\sqrt{\mu _{m}}\leq N} \left| a_{m}^{(N)} \right|^{2} e^{-2\mu _{m}t} \, dt \le \mathcal{M}  ( \alpha ,T )  \intdouble_{Q_T } e^{-\frac{ \alpha \beta}{t}} \left \vert P_{N} ( t,x ) \right \vert^{2}\, dx\, dt, 
	\end{equation*}
where $\mathcal{M}  ( \alpha ,T )$  is given in the statement of Proposition~\ref{Prop_restrictouvertLR}.

To summarize, if $\alpha >2\beta $, after integrating~\eqref{maj_1} on $\left( \frac{\alpha }{N},T\right) $, we have
	\begin{equation*}
\int_{\frac{\alpha }{N}}^{T}\int_{\Omega }e^{-\frac{\beta \alpha }{t}}\left \vert P_{N} ( t,x ) \right \vert ^{2}dxdt\leq \left( 3 \left\| P_{N} \right \|_{L^2 ((\frac \alpha N,T) \times \omega )}^{2} + \mathcal{M}  ( \alpha ,T ) \intdouble_{Q_T } e^{-\frac{ \alpha \beta}{t}} \left \vert P_{N} ( t,x ) \right \vert^{2}\, dx\, dt \right) .
	\end{equation*}
Adding the previous inequality and inequality~\eqref{Eq_4ouvertLR}) we get~\eqref{Eq_3ouvertLR}. This ends the proof. 
\end{proof}


Let us continue with our reasoning. The following result is our third step (see Theorem~\ref{tPN}):


\begin{proposition}\label{Prop_alpha_0}
For all $T>0$ there exists $\alpha_0(T,\beta)>0$ satisfying
	\begin{equation}\label{limalpha_0LR}
\lim_{T\to 0^+} \alpha_0(T,\beta) = 2\beta,
	\end{equation}
and such that any $N\geq 1$, any $P_N$ given by~\eqref{P_Nouvert} and for all $\alpha \geq \alpha _{0} ( T,\beta  ) $, one has
	\begin{equation}\label{IsomLR}
\intdouble_{Q_T} e^{-\frac{ \alpha \beta }{t}} \left \vert P_{N} ( t,x ) \right \vert ^{2} \, dx\, dt \leq 6\int_{0}^{T} \!\!\int_{\omega } \left \vert P_{N} ( t,x ) \right \vert ^{2}\,dx \, dt.
	\end{equation}
\end{proposition}


\begin{proof}
Let us first take $\alpha > 2 \beta$ and consider the expressions of $\calM (\alpha, T)$ and $\chi (T, \alpha)$ in the statement of Proposition~\ref{Prop_restrictouvertLR}. If we take $\alpha = 2 \beta + \sqrt{ \delta T}$, with $\delta \ge 1$ to be determined, one has $\alpha
^{2} - 2\alpha \beta = \left( 2 \beta + \sqrt{ \delta T} \right)\sqrt{ \delta T} $ and we can write
	\begin{equation*}
	\begin{array}{l}
\Ds \calM (\alpha, T) =  \widehat{\mathcal{C}} \,  (T + 1) \frac{ \frac{ T^{\theta _{1} + 1} \left( 1 + T^{\theta _{1}+1} \right) }{ \left( 2 \beta + \sqrt{ \delta T} \right)^{2 \theta_1 +1 } } + 2 \beta + \sqrt{ \delta T} }{ T \left( \delta T\right)^{ \theta_1 + \frac 12} } e^{- \frac \delta 2} e^{- \beta \sqrt{\frac \delta T}} \\
	\noalign{\smallskip}
\Ds \phantom{ \calM (\alpha, T) }
\le  \widehat{\mathcal{C}} \,  \frac{T + 1}{T^{\theta_1 + \frac 32}} \left[ T^{\frac 12 } \left( 1 + T^{\theta_1 + 1} \right) + 2 \beta + T^{\frac 12 } \right] e^{- \frac \delta 2} := B(\beta , T ) \, e^{- \frac \delta 2}  , \quad \forall T > 0.
	\end{array}
	\end{equation*}
Recall that $\widehat \calC > 0$ is the constant provided by Proposition~\ref{Prop_restrictouvertLR}, only depending on $\theta_1$. 

From the previous inequality, it is clear that taking $\alpha_0 (T, \beta) = 2 \beta +  \sqrt{T \delta_0 (T, \beta)} $ with
	\begin{equation*}
\delta_0(T, \beta) = \max\left\{1,2 \ln\left(2 B(\beta , T ) \right)\right\}
	\end{equation*} 
one has~\eqref{limalpha_0LR}. In addition, $\calM (\alpha, T) \le 1/2$, for any $\alpha \ge \alpha_0 (T, \beta)$ and inequality~\eqref{Eq_3ouvertLR} implies~\eqref{IsomLR}. This finalizes the proof. 
\end{proof}


Let $\eta _{\alpha }$ and $L_{\eta _{\alpha }}^{2}((0,T)\times \Omega )$ the function and the space defined in~\eqref{eta} and~\eqref{L2etaalpha} with $\Omega _{1}$ replaced by $\Omega $. As in~\eqref{eta_alpha}, we define 
	\begin{equation}\label{eta_alpha2}
	\left\{ 
	\begin{array}{l}
E_{\eta _{\alpha }}=\overline{\mathrm{span}\, \left \{ F_{m}:m\geq 1\right \} }^{L_{\eta_\alpha }^{2}(Q_T )}, \\ 
	\noalign{\smallskip}
E^{\omega }=\overline{\mathrm{span}\, \left \{ F_{m}|_{\omega }: m \ge 1 \right \}}^{L^{2}((0,T)\times \omega )}.
	\end{array}
	\right.  
	\end{equation}
where the function $F_m$ is given in~\eqref{FmLR}. As a consequence of Propositions~\ref{Prop_epsLR}, \ref{Prop_restrictouvertLR} and~\ref{Prop_alpha_0}, we can prove a result as Theorem~\ref{tmainrest} in our framework. One has:


%
\begin{theorem}\label{tmainrest2}  
Let us assume that $\Lambda _{1}=\{ \mu _{m}\}_{m\geq }$ satisfies~\eqref{H1} and that $\mathcal{B}_{1}=\{ \psi _{m}\}_{m\geq 1}$ is an orthonormal basis of $L^{2}(\Omega )$ fulfilling~\eqref{NSpI}.  Let $T>0$ and let us consider the constant $\alpha_0(T,\beta)$ provided by Proposition~\ref{Prop_alpha_0}. Then, for all $\alpha \geq \alpha _{0} ( T,\beta  )$, the operator  
	\begin{align}
\mathcal{R}_{\omega } : L_{\eta _{\alpha }}^{2} (Q_T )& \rightarrow L^{2}((0,T)\times \omega ) \notag  \\
\varphi & \mapsto \mathcal{R}_{\omega }(\varphi )=\varphi |_{\omega }  \notag
	\end{align}
satisfies 
	\begin{equation*}
\left\Vert \varphi \right\Vert^2 _{L^2_{\eta_\alpha}(Q_T)} \leq 7\left \Vert \mathcal{R}_{\omega }(\varphi )\right \Vert _{L^{2}((0,T)\times \omega )}^{2}\leq 7 \left\Vert \varphi \right\Vert^2 _{L^2_{\eta_\alpha}(Q_T)},\quad \forall \varphi \in E_{\eta _{\alpha }}, 
	\end{equation*}
and, therefore, $\mathcal{R}_{\omega }\in \mathcal{L}\left( E_{\eta _{\alpha }},E^{\omega }\right) $ is an isomorphism.
\end{theorem}

The proof is analogous to the proof of Theorem~\ref{tmainrest} and will be omitted.

Let us now prove a result which plays the role of Theorem~\ref{main} in our framework.  One has:


\begin{theorem} \label{tmainrest3} 
Under the assumptions of Theorem~\ref{tmainrest2}, there exists a positive constant $\widehat \calC$ such that, for any $T > 0$, the family $\{F_{m}\}_{m\geq 1}$ (see~\eqref{FmLR}) has a unique biorthogonal family $\left\{Q_{m}\right\}_{m\geq 1} \subset E^\omega$ in $L^{2}((0,T)\times \omega )$ that satisfies 
	\begin{equation}\label{f123}
\left\Vert Q_{m}  \right\Vert _{L^{2}((0,T)\times \omega )}^{2} \leq \widehat \calC \left( \dfrac{1}{T-\varepsilon }+\mu _{m}\right) e^{\frac{\alpha \beta }{\varepsilon } + 2\varepsilon \mu _{m}}, \quad \forall m\geq 1, 
	\end{equation}
for any $\varepsilon \in (0, T) $ and $\alpha \ge \alpha_ 0 (T, \beta)$ ($\alpha_0 (T, \beta)$ is the constant provided by Proposition~\ref{Prop_alpha_0}).
\end{theorem}


\begin{proof}
The proof is a consequence of Theorem~\ref{tmainrest2}. Let us take $T>0$, $\alpha \geq \alpha _{0}(T,\beta)$ and $\varepsilon \in  ( 0,T ) $. We consider, 
	\begin{equation*}
	\left\{ 
	\begin{array}{l}
\Ds \widetilde f_{m}^\varepsilon (t,x) = e^{ \frac{\eta_ \alpha (x) }t } f_{m}^\varepsilon (t,x) , \quad \forall (t, x) \in Q_T,  \\
	\noalign{\smallskip}
\Ds Q_{m}:= \left( \mathcal{R}_{\omega }^{-1}\right) ^{\star } \mathcal{P}_{\eta _{\alpha }} \widetilde f_{m}^\varepsilon \in E^\omega, \quad m\geq 1,
	\end{array}
	\right.
	\end{equation*}
where he sequence $\left \{ f_{m}^{\varepsilon } \right \}_{m\geq 1}$ is given by Proposition~\ref{Prop_epsLR} and satisfies~\ref{fepsestimLR}, $E_{\eta _{\alpha }}$ in \eqref{eta_alpha2}, $\mathcal{R}_{\omega }$ is the restriction operator defined in Theorem~\ref {tmainrest2} and $\mathcal{P}_{\eta _{\alpha }}$ is the orthogonal projection from $L_{\eta _{\alpha }}^{2}( Q_T )$ to $E_{\eta _{\alpha }}$. As in the proof of Theorem~\ref{main} (see Section~\ref{ss42}), the sequence $\{ Q_{m} \}_{m \ge 1}$ belongs to $E^\omega$. This will imply the uniqueness (and therefore, the independence with respect to the parameter $\varepsilon$) and the optimality of this family.

Observe that $\widetilde f_{m}^\varepsilon \in L_{\eta _{\alpha }}^{2}( Q_T )$ and, for any $m \ge 1$,
	\begin{equation*}
\left\| \widetilde f_{m}^\varepsilon \right\|_{ L_{\eta _{\alpha }}^{2}( Q_T )}^2 = \int_\varepsilon^T \! \! \int_{\Omega } e^{\frac{\eta_ \alpha (x) }t } \left| f_{m}^\varepsilon (t,x)\right|^2 dx \, dt \le c\left( \dfrac{1}{T-\varepsilon }+\mu _{m}\right) e^{\frac{\alpha \beta }{\varepsilon } + 2\varepsilon \mu _{m}}, 
	\end{equation*}
thanks to~\eqref{fepsestimLR}. On the other hand, 
	\begin{align*}
\Ds \delta_{nm} & = \intdouble_{Q_T} F_n (t,x) f_{m}^\varepsilon (t,x) \, dx \, dt = \left( F_n, \widetilde f_{m}^\varepsilon \right)_{ L_{\eta _{\alpha }}^{2}( Q_T )} = \left( F_n, \mathcal{P}_{\eta _{\alpha }} \left(\widetilde f_{m}^\varepsilon \right) \right)_{ E_{\eta_\alpha }} \\
	\noalign{\smallskip}
\Ds & =  \left( \mathcal{R}_{\omega }^{-1} \mathcal{R}_{\omega } F_n, \mathcal{P}_{\eta _{\alpha }} \left(\widetilde f_{m}^\varepsilon \right) \right)_{ E_{\eta_\alpha }} = \left( \mathcal{R}_{\omega } F_n, \left( \mathcal{R}_{\omega }^{-1}\right) ^{\star }  \mathcal{P}_{\eta _{\alpha }} \left(\widetilde f_{m}^\varepsilon \right) \right)_{ L^2 ((0, T) \times \omega )} \\
	\noalign{\smallskip}
\Ds & = \int_0^T \! \! \int_\omega F_n (t, x) Q_m (t, x) \,dx \, dt, \quad \forall n,m \ge 1.  
	\end{align*}
Therefore, the family $\left\{Q_{m} \right\}_{m\geq 1}$ is biorthogonal to $\{F_{m}\}_{m\geq 1}$ in $L^{2}((0,T)\times \omega )$. Finally, 
	\begin{align*}
\left\| Q_{m} \right \|_{ L^2((0,T)\times\omega)}^2 = & \left\| \left( \mathcal{R}_{\omega }^{-1}\right) ^{\star } \mathcal{P}_{\eta _{\alpha }} \widetilde f_{m}^\varepsilon \right\|_{ L^2((0,T)\times\omega)}^2
\leq 7 \left\| \calP_{\eta _{\alpha }} \widetilde f_{m}^\varepsilon \right\|_{ L^2_{\eta_\alpha}(Q_T)}^2
\\
&\leq 7 \left\| \widetilde f_{m}^\varepsilon \right\|_{ L^2_{\eta_\alpha}(Q_T)}^2, \quad \forall m \ge 1. 
	\end{align*}
From this inequality we deduce~\eqref{f123}. This ends the proof. 
\end{proof}

We can now prove Theorem~\ref{mainLR}

\begin{proof}[Proof of Theorem~\ref{mainLR}]
 The null controllability problem~\eqref{controlchaleurd} in $L^2 (\Omega)$ is equivalent to a moment problem in $L^2 ((0,T) \times \omega)$ for the family $\left\{ F_m \right\}_{m \ge 1}$ given by~\eqref{FmLR}. This moment problem can be solved by means of the biorthogonal family $\left\{Q_{m} \right\}_{m\geq 1}$ provided by Theorem~\ref{tmainrest3} with $\alpha = \alpha_0(T, \beta) > 0$ . Thus, let us take  the family $\left\{Q_{m}  \right\}_{m\geq 1}$. An explicit solution to problem~\eqref{controlchaleurd} ($y_0 \in L^2 (\Omega)$) is 
	\begin{equation*}
u (t, x) =-\sum_{m\geq 1}e^{-\mu _{m}T}\left ( y_{0},\psi _{m}\right )_{L^2 (\Omega )} Q_{m}  (T -t, x), \quad \forall (t, x) \in (0,T)  \times \omega .
	\end{equation*}
Therefore (see~\eqref{f123}), for all $\eps\in (0,T)$
	\begin{align*}
\Vert u\Vert _{L^{2}((0,T)\times \omega )}^{2} &\leq \left( \sum_{m\geq 1}e^{-2\mu _{m}T} \left\Vert Q_{m} \right \Vert _{L^{2}((0,T)\times \omega )}^{2}\right) \left \Vert y_{0}\right \Vert _{L^{2}\left( \Omega \right) }^{2} 
	\\
&\leq \widehat \calC e^{\frac{\alpha_0 \beta }{\varepsilon }} \left( \sum_{m\geq 1}\left( \dfrac{1}{T-\varepsilon }+\mu _{m}\right) e^{-2\mu _{m} ( T-\varepsilon  ) }\right) \left \Vert y_{0}\right \Vert _{L^{2} ( \Omega  ) }^{2}  
	\\
& := \widehat \calC e^{\frac{\alpha_0 \beta }{\varepsilon }} S(\varepsilon, T) \left \Vert y_{0}\right \Vert _{L^{2} ( \Omega  ) }^{2}, \quad \forall \eps \in (0,T).
	\end{align*}

Using condition~\eqref{H1}, it is possible to estimate $S(\varepsilon, T)$ ($\widehat \calC $ is a generic constant only depending on $\theta_1$):
	\begin{align*}
 S ( \varepsilon , T )  &= \int_{0}^{\infty }\left( \dfrac{1}{T-\varepsilon }+x \right) e^{-2 ( T-\varepsilon  ) x} \, d\mathcal{N}_{\Lambda _{1}} ( x ) 
 	\\
&\leq  \kappa_1 \int_{0}^{\infty }\left( x^{\theta _{1}} + 2 ( T-\varepsilon  ) x^{\theta _{1}+1} \right) e^{-2 ( T-\varepsilon  ) x} \, dx.
 	\\
& \le \widehat \calC \left( \frac{1}{2 ( T-\varepsilon  ) }\right)^{\theta_{1}+1} \int_{0}^{\infty }\left( \xi ^{\theta _{1}} + \xi ^{\theta_{1}+1}\right) e^{-\xi }\, d\xi = \frac{ \widehat \calC }{(T - \varepsilon)^{\theta_1 +1 }}
	\end{align*}

Coming back to $\Vert u\Vert _{L^{2}((0,T)\times \omega )}$, we deduce
	\begin{equation*}
\Vert u\Vert _{L^{2}((0,T)\times \omega )}^{2} \leq \widehat \calC e^{\frac{\alpha_0 \beta }{\varepsilon }} \frac{ 1  }{(T - \varepsilon)^{\theta_1 +1 }} \left \Vert y_{0}\right \Vert _{L^{2} ( \Omega  ) }^{2},\quad \forall \eps\in (0,T) . 
	\end{equation*}
Let us now choose $\varepsilon \in (0,T)$ in order to minimize the right member of the previous estimate. To this end, we consider the function $g$ defined on $ ( 0,T ) $ by: 
	\begin{equation*}
g ( \varepsilon  ) =e^{\frac{\alpha_0 \beta }{\varepsilon }} ( T-\varepsilon  )^{-\sigma }, \quad \forall \varepsilon \in  ( 0,T ) , \quad \sigma =\theta _{1}+1. 
	\end{equation*}
This function achieves its minimum in $(0, T)$ at point
	\begin{equation*}
\varepsilon_0 (T, \beta) := \frac{\sqrt{\alpha_0 ^{2}\beta ^{2}+4T\alpha_0 \sigma \beta } - \alpha_0 \beta }{2\sigma  }=\frac{2T\alpha_0 \beta}{\sqrt{\alpha_0 ^{2}\beta ^{2}+4\sigma T\alpha_0 \beta }+\alpha_0 \beta } \in (0, T). 
	\end{equation*}
This expression implies the estimate~\eqref{coutuLR} for the constant ($\sigma =\theta _{1}+1$)
	\begin{equation*}
\calK (T, \beta) := \widehat \calC e^{\frac{\alpha_0 \beta }{\varepsilon_0 (T, \beta) }} \frac{ 1  }{(T - \varepsilon_0 (T, \beta) )^{\theta_1 +1 }}. 
	\end{equation*}

Finally, let us check~\eqref{CTHLR}. For this purpose, we will use property~\eqref{limalpha_0LR}. One has,
	\begin{equation*}
 \lim_{T \to 0^+} \varepsilon_0 (T, \beta) = T, \quad \lim_{T \to 0^+} \frac{T - \varepsilon_0 (T, \beta)}{T^2} = \frac{\theta_1 + 1}{2 \beta^2}, 
	\end{equation*}
and, then,
	\begin{align*}
\limsup_{T \to 0^+} T \ln \left( \calK (T, \beta) \right) & = \lim_{T \to 0^+}  \frac{\alpha_0 \beta T}{\varepsilon_0 (T, \beta) } - (\theta_1 +1 ) \lim_{T \to 0^+} T \ln (T - \varepsilon_0 (T, \beta) ) 
	\\
& = 2 \beta^2. 
	\end{align*}
This proves~\eqref{limalpha_0LR} and ends the proof of Theorem~\ref{mainLR}. 
\end{proof}



\begin{thebibliography}{99}
\bibitem{ABM_2018} 
\textsc{D.~Allonsius, F.~Boyer, M.~Morancey}, 
\emph{Spectral analysis of discrete elliptic operators and applications in control theory}, 
Numer. Math.~\textbf{140} (2018), no.~4, pp.~857--911. 


\bibitem{ABGD:14} 
\textsc{F.~Ammar Khodja, A.~Benabdallah, M.~Gonz\'{a}lez-Burgos, L.~de Teresa}, 
\emph{Minimal time for the null controllability of parabolic systems: the effect of the condensation index of complex sequences}, 
J.~Funct.\ Anal.~\textbf{267} (2014), no.~7, pp.~2077--2151.


\bibitem{AKBGBdT_JMAA16}
\textsc{F.~Ammar~Khodja, A.~Benabdallah, M.~Gonz\'alez-Burgos, and L.~de~Teresa},
\emph{New phenomena for the null controllability of parabolic systems: minimal time and geometrical dependence}, J.~Math.~Anal.~Appl.~\textbf{444} (2016), no.~2, pp.~1071--1113. 
 

\bibitem{BBGBO:2014}
\textsc{A.~Benabdallah, F.~Boyer, M.~Gonz{\'{a}}lez-Burgos, G.~Olive}, 
\emph{Sharp estimates of the one-dimensional boundary control cost for parabolic systems and application to the $N$-dimensional boundary null controllability in cylindrical domains}, 
SIAM J.~Control and Optim.~\textbf{52} (2014), no.~5, p.~2970-3001.


\bibitem{BBM} 
\textsc{A.~Benabdallah, F.~Boyer, M.~Morancey}, 
\emph{A block moment method to handle spectral condensation phenomenon in parabolic control problems},
Ann.~H.~Lebesgue \textbf{3} (2020), pp.~717-793.


\bibitem{B:22}
\textsc{F.~Boyer}, \emph{Controllability of linear parabolic equations and systems}, Master, France. 2022. 
https://hal.science/hal-02470625


\bibitem{BO:24}
\textsc{F.~Boyer and G.~Olive}, 
\emph{Boundary null-controllability of some multi-dimensional linear parabolic systems by the moment method}, 
to appear in Annales de l'Institut Fourier. 


\bibitem{BM:23}
\textsc{F.~Boyer and M.~Morancey}, 
\emph{Analysis of non scalar control problems for parabolic systems by the block moment method}, 
C.~R.~Math.~Acad.~Sci.\ Paris \textbf{361} (2023), pp.~1191--1248.


\bibitem{CMV:20}
\textsc{P.~Cannarsa, P.~Martinez and J.~Vancostenoble}, 
\emph{The cost of controlling strongly degenerate parabolic equations}, 
ESAIM Control Optim.~Calc.~Var.~\textbf{26} (2020), Paper No.~2, 50 pp.



\bibitem{CMV:21}
\textsc{P.~Cannarsa, P.~Martinez and J.~Vancostenoble}, 
\emph{Sharp estimate of the cost of controllability for a degenerate parabolic equation with interior degeneracy}, 
Minimax Theory Appl.~\textbf{6} (2021), no.~2, pp.~251--280.


\bibitem{Duprez_Tmin2017}
\textsc{M.~Duprez}, 
\emph{Controllability of a $2\times 2$ parabolic system by one force with space-dependent coupling term of order one}, 
ESAIM Control Optim.~Calc.~Var.~\textbf{23} (2017), no.~4, pp.~1473--1498.


\bibitem{Dolecki:73} 
\textsc{S.~Dolecki}, 
\emph{Observability for the one-dimensional heat equation}, 
Studia Math.~\textbf{48} (1973), pp.~291--305.

\bibitem{Fatto} 
\textsc{H.O.~Fattorini, }
\emph{Boundary control of temperature distributions in a parallelepipedon},
 SIAM J.~Control~\textbf{13} (1975), pp.~1--13.

\bibitem{FR1:71} 
\textsc{H.O.~Fattorini, D.~L.~Russell}, 
\emph{Exact controllability theorems for linear parabolic equations in one space dimension}, 
Arch.~Rational Mech.~Anal.~\textbf{43} (1971), pp.~272--292.

\bibitem{FR2:74} 
\textsc{H.O.~Fattorini, D. L. Russell}, 
\emph{Uniform bounds on biorthogonal functions for real exponentials with an application to the control theory of parabolic equations}, Quart.~Appl.\ Math.~\textbf{32} (1974/75), pp.~45--69.


\bibitem{FCGBdT} 
\textsc{E.~Fern\'andez-Cara, M.~Gonz\'alez-Burgos, L.~de Teresa}, 
\emph{Boundary controllability of parabolic coupled equations}. 
J.~Funct.~Anal.~\textbf{259} (2010), no.~7, pp.~1720--1758.

\bibitem{FI:96} 
\textsc{A.~Fursikov, O.~Yu.~Imanuvilov}, \emph{Controllability of Evolution Equations}, 
Lecture Notes Ser., \textbf{34}, Seoul National University, Research Institute of Mathematics, Global Analysis Research Center, Seoul, 1996.


\bibitem{GBO:19} 
\textsc{M.~Gonz\'{a}lez-Burgos, L.~Ouaili}, 
\emph{Sharp estimates for biorthogonal families to exponential functions associated to complex sequences without gap conditions},
Evol.~Equ.~Control Theory \textbf{13} (2024), no.~1, pp.~215--279.


\bibitem{JL:90}
\textsc{D.~Jerison, G.~Lebeau}, 
\emph{Nodal sets of sums of eigenfunctions} in \textquotedblleft Harmonic Analysis and Partial
Differential Equations\textquotedblright, Chicago Lectures in Math, University of Chicago Press, Chicago, IL, 1999, pp.~223--239.


\bibitem{K:95} 
\textsc{A.~Kirsch}, 
\emph{An Introduction to the Mathematical Theory of Inverse Problems}, 
Second edition, Appl. Math. Sci.,~\textbf{120}, Springer, New York, 2011. 


\bibitem{LR:95} 
\textsc{G.~Lebeau, L.~Robbiano}, 
\emph{Contr\^{o}le exact de l'\'{e}quation de la chaleur}, 
Comm.\ Partial Differential Equations \textbf{20} (1995), no.~1-2, pp.~335--356.


\bibitem{M:10} 
\textsc{L.~Miller}, 
\emph{A direct Lebeau-Robbiano strategy for the observability of heat-like semigroups}, 
Discrete Contin.~Dyn.\ Syst.~Ser.~B \textbf{14} (2010), no.~4, pp.~1465--1485.


\bibitem{Morancey24} 
\textsc{M.~Morancey}, 
\emph{Some remarks on the moment method for null controllability of parabolic equations in higher dimension}; \textit{In preparation}.  

\bibitem{Ou:19} 
\textsc{L.~Ouaili}, 
\emph{Minimal time of null controllability of two parabolic equations}, 
Math.~Control Relat.~Fields \textbf{10} (2020), no.~1, pp.~89--112.

\bibitem{Ou:19-b} 
\textsc{L.~Ouaili}, 
\emph{Contr\^olabilit\'e de Quelques Syst\`emes Paraboliques}, Ph.D thesis, Aix-Marseille University, 2020. 


\bibitem{EHS:15} 
\textsc{E.~H.~Samb}, 
\emph{Internal null-controllability of the $N$-dimensional heat equation in cylindrical domains}, 
C.~R.~Acad.~Sci.~Paris~\textbf{353} (2015), no.~10, pp.~925--930.


\bibitem{Schwartz} 
\textsc{L.~Schwartz}, \emph{\'Etude des Sommes d'Exponentielles R\'eelles}, 
Actualit\'es Sci.~Ind., no.~959. Hermann et Cie., Paris, 1943.

\end{thebibliography}
\end{document}